\documentclass[11pt,a4paper]{article}

\usepackage{amsmath,amssymb}
\usepackage{color}
\newcommand{\R}{\mathbb{R}}

\newcommand{\Z}{\mathbb{Z}}
\def\cA{{\mathcal A}}
\def\cC{{\mathcal C}}

\def\cS{{\mathcal S}}
\newcommand{\ee}{\varepsilon}
\renewcommand{\div}{{\rm div}\,}
\newcommand{\curl}{{\rm curl}\,}

\newcommand{\Frac}{\displaystyle \frac}

\newcommand{\Sum}{\displaystyle \sum}

\def\d{\partial}
\def\ddl{\dot \Delta_l}
\def\ddj{\dot \Delta_j}
\def\ddq{\dot \Delta_q}
\def\ddk{\dot \Delta_k}
\def\tilde{\widetilde}
\def\hat{\widehat}

\newcommand{\n}{\nabla}
\newcommand{\fd}{\frac{d}{2}}

\newcommand{\p}{\partial}

\newcommand{\qe}{q_\ee}
\newcommand{\ue}{u_\ee}

\newtheorem{thm}{Theorem}
\newtheorem{lem}{Lemma}

\newtheorem{prop}{Proposition}
\newtheorem{defi}{Definition}

\newtheorem{rem}{Remark}

\title{On a Lagrangian method for the convergence from a non-local to a local Korteweg capillary fluid model}

\author{Fr\'ed\'eric Charve\footnote{Universit\'e Paris-Est Cr\'eteil, Laboratoire d'Analyse et de Math\'ematiques Appliqu\'ees (UMR 8050), 61 Avenue du G\'en\'eral de Gaulle, 94 010 Cr\'eteil Cedex (France). E-mail: frederic.charve@u-pec.fr}, Boris Haspot  \thanks{Ceremade UMR CNRS 7534
Universit\'e de Paris  Dauphine,
Place du Mar\'echal DeLattre De Tassigny
75775 Paris Cedex 16. E-mail: haspot@ceremade.dauphine.fr}}

\date{}
\begin{document}

\maketitle

\begin{abstract} In the present article we are interested in further investigations for the barotropic compressible Navier-Stokes system endowed with a non-local capillarity we studied in \cite{CH}. Thanks to an accurate study of the associated linear system using a Lagrangian change of coordinates, we provide more precise energy estimates in terms of hybrid Besov spaces naturally depending on a threshold frequency ( which is determined in function of the physical parameter) distinguishing the low and the high regimes. It allows us in particular to prove  the convergence of the solutions from the non-local to the local Korteweg system. Another mathematical interest of this article is the study of the effect of the Lagrangian change on the non-local capillary term.
\end{abstract}
\section{Introduction}
\subsection{Presentation of the system}
The local and non-local Korteweg systems aim to study the dynamics of a liquid-vapour mixture in the diffuse interface approach (DI), where the phase changes are seen through the variations of the density. These systems are based upon the compressible Navier-Stokes system with a Van der Waals state law for ideal fluids, and endowed with a capillary tensor modelling the behaviour at the interfaces between the phases. This capillary term was introduced in the DI approach in order to obtain physically relevant solutions by penalizing the high variations of the density.

We refer to \cite{CH} for a physical presentation of the diffuse interface model, and of the local and non-local Korteweg systems. Let us recall that the local model of the capillary term was introduced by Korteweg and the non-local model was introduced by Van der Waals and renewed by F. Coquel, D. Diehl, C. Merkle and C. Rohde (for an in-depth presentation of the capillary models, we refer to \cite{Rohdehdr} and \cite{5CR}).

Let $\rho$ and $u$ denote the density and the velocity of a compressible viscous fluid. As usual, $\rho$ is a non-negative function and $u$ is a vector-valued function defined on $\R^d$. In the sequel we will denote by $\cA$ the following diffusion operator
$$
\cA u= \mu \Delta u+ (\lambda+\mu)\nabla \div u, \quad \mbox{with} \quad \mu>0 \quad \mbox{and} \quad \nu=\lambda+ 2 \mu >0.
$$
The Navier-Stokes equations for compressible fluids endowed with internal capillarity read:
$$
\begin{cases}
\begin{aligned}
&\d_t\rho+\div (\rho u)=0,\\
&\d_t (\rho u)+\div (\rho u\otimes u)-\cA u+\nabla(P(\rho))=\kappa\rho\nabla D[\rho].\\
\end{aligned}
\end{cases}
$$
Let us mention that the capillary coefficient $\kappa$ may depend on $\rho$ but in this article only the constant case is considered. In the local Korteweg system $(NSK)$, the capillary term $D[\rho]$ is given by (see \cite{3DS}):
$$D[\rho]=\Delta\rho,$$
and, in the non-local Korteweg system $(NSRW)$ (introduced in its modern form by C. Rohde in \cite{5Ro} and also \cite{5CR}, see Van der Waals \cite{VW} for the original works), if $\phi$ is an interaction potential which satisfies the following conditions
$$
 (|.|+|.|^2)\phi(.)\in L^1(\R^d)\mbox{, }\quad\int_{\R^d}\phi(x)dx=1,\quad\phi\mbox{ even, and }\phi\geq0,
$$
$D[\rho]$ is a non-local term:
$$D[\rho]=\phi*\rho-\rho.$$

If we compute the Fourier transform of the capillary terms, we obtain $(\hat{\phi}(\xi)-1) \hat{\rho}(\xi)$ in the non-local model, and $-|\xi|^2 \hat{\rho}(\xi)$ in the local model.

We are interested in the closedness of the solutions of these models when $\hat{\phi}(\xi)$ is formally "close" to $1-|\xi|^2$. In \cite{CH}, we approximated the local Korteweg model $(NSK)$ with a non-local model such as system $(NSRW)$ where we chosed a specific function $\phi_{\ee}$ in the capillarity tensor. In this paper we will once more consider the following non-local system:
$$
\begin{cases}
\begin{aligned}
&\d_t\rho_\ee+\div (\rho_\ee \ue)=0,\\
&\d_t (\rho_\ee \ue)+\div (\rho_\ee u\otimes \ue)-\cA \ue+\nabla(P(\rho_\ee))=\rho_\ee\frac{\kappa}{\ee^2}\nabla(\phi_\ee*\rho_\ee-\rho_\ee),\\
\end{aligned}
\end{cases}
\leqno{(NSRW_\ee)}
$$
where
$$\phi_{\ee}=\frac{1}{\ee^d} \phi(\frac{x}{\ee}) \quad \mbox{with} \quad \phi(x)=\frac{1}{(2 \pi)^d} e^{-\frac{|x|^2}{4}}$$
For a fixed $\xi$ the Fourier transform of $\phi_\ee$ is $\hat{\phi_\ee}(\xi)=e^{-\ee^2|\xi|^2}$, and when $\ee$ is small, $\Frac{\hat{\phi_\ee}(\xi)-1}{\ee^2}$ is close to $-|\xi|^2$.

We will consider a density which is close to an equilibrium state $\overline{\rho}$ and we will introduce the change of function $\rho= \overline{\rho}(1+q)$.  For simplicity we take $\overline{\rho}=1$. The previous systems become:
$$
\begin{cases}
\begin{aligned}
&\d_t q+ u.\nabla q+ (1+q)\div u=0,\\
&\d_t u+ u.\nabla u -\cA u+P'(1).\nabla q-\kappa\nabla \Delta q= K(q).\nabla q- I(q) \cA u,\\
\end{aligned}
\end{cases}
\leqno{(K)}
$$
and
$$
\begin{cases}
\begin{aligned}
&\d_t \qe+ \ue.\nabla \qe+ (1+\qe)\div \ue=0,\\
&\d_t \ue+ \ue.\nabla \ue -\cA \ue+P'(1).\nabla \qe-\frac{\kappa}{\ee^2}\nabla(\phi_\ee*\qe-\qe)=K(\qe).\nabla \qe- I(\qe) \cA \ue,\\
\end{aligned}
\end{cases}
\leqno{(RW_\ee)}
\label{Rexpr}
$$
where $K$ and $I$ are real-valued functions defined on $\R$ given by:
$$
K(q)=\left(P'(1)-\frac{P'(1+q)}{1+q}\right) \quad \mbox{and} \quad I(q)=\frac{q}{q+1}.
$$

\subsection{Existence results}

Let us now recall some results concerning the local and non-local Korteweg systems. As for the compressible Navier-Stokes system (we refer to  \cite{Dbook}, \cite{Dinv}, \cite{arma}, \cite{CD}) both of these systems have been studied in the context of the existence of global strong solutions with small initial data in critical spaces for the scaling of the equations. For example, concerning the strong solutions, we refer to \cite{DD}, \cite{Has7} (and \cite{Has2} in the non isothermal case) for a study of $(NSK)$ system, and to \cite{Has1} for $(NSRW)$. In \cite{Has8}, we show the existence of global strong solution with large initial data on the rotational part when we add a friction term.
\\
Let us mention that the well-posedness of the compressible Euler Korteweg system (when $\mu=\lambda=0$) has been studied in the case of variable capillary coefficient by Benzoni, Danchin and Descombes in \cite{BDD1}.

The solutions of the compressible Navier Stokes or non-local Korteweg systems have the same behaviour. Namely the density regularity is separated by a frequency threshold: in low frequencies, the solution is subject to a heat-type smoothing, and in the high frequencies, there is only a damping effect due to the term of pressure (modulo that the pressure is at least locally increasing with respect to the density). The solution of the local Korteweg system is more regular: for all frequencies, we have a parabolic regularization on the density (see \cite{DD,Has7}). We refer to the appendix for the definitions of the Besov spaces introduced in the following results. 

\begin{thm} (\cite{DD})
\sl{ Assume that $P'(1)>0$, $\min(\mu, \nu)>0$ where $\nu\overset{def}{=}2 \mu+\lambda$, that the initial density fluctuation $q_0$ belongs to $\dot{B}_{2,1}^{\fd-1}\cap \dot{B}_{2,1}^{\fd}$, and that the initial velocity $u_0$ is in $(\dot{B}_{2,1}^{\fd-1})^d$. Then there exist constants $\eta_K>0$ and $C>0$ depending on $\kappa$, $\mu$, $\nu$, $P'(1)$ and $d$ such that if:
$$
\|q_0\|_{\dot{B}_{2,1}^{\fd-1} \cap \dot{B}_{2,1}^{\fd}} +\|u_0\|_{\dot{B}_{2,1}^{\fd-1}}\leq \eta_K
$$
then system $(K)$ has a unique global solution $(\rho, u)$ such that the density fluctuation and the velocity satisfy:
$$
\begin{cases}
q \in \cC(\R_+, \dot{B}_{2,1}^{\fd-1} \cap \dot{B}_{2,1}^{\fd}) \cap L^1(\R_+, \dot{B}_{2,1}^{\fd+1} \cap \dot{B}_{2,1}^{\fd+2}),\\
u\in \cC(\R_+, \dot{B}_{2,1}^{\fd-1})^d \cap  L^1(\R_+, \dot{B}_{2,1}^{\fd+1})^d.
\end{cases}
$$
Moreover the norm of $(q,u)$ in this space is estimated by the initial norm $C (\|q_0\|_{\dot{B}_{2,1}^{\fd-1} \cap \dot{B}_{2,1}^{\fd}} +\|u_0\|_{\dot{B}_{2,1}^{\fd-1}})$.}
\end{thm}
\begin{rem}
 \sl{Further in this article R. Danchin and B. Desjardins provide a Fourier study of the linearized system and observe different behaviours whether the quantity $\nu^2- 4\kappa$ is positive, negative of zero. In all cases they obtain parabolic regularization.}
\end{rem}
\begin{rem}
 \sl{In \cite{Has6}, the second author obtains some generalizations of \cite{DD} in as much as with a specific choice on the capillarity and the viscosity $\kappa(\rho)=\frac{1}{\rho}$, $\mu(\rho)=\rho$, we obtain the existence of global strong solution with $u_{0}\in B^{\frac{d}{2}-1}_{2,\infty}$ and $\ln(1+q_{0})\in B^{\frac{d}{2}}_{2,\infty}\cap B^{\frac{d}{2}-1}_{2,\infty}$ (let us point out that in this case we can work with discontinuous initial density). We also refer to \cite{Has8} for the existence of global strong solution with large initial data on the irrotational part when we added a friction term. In the sequel as we shall deal with a Lagrangian change of coordinate which requires a Lipschitz control on the velocity, we need to work in the framework of the functional space introduced in \cite{DD}. The second reason is that the existence of global strong solution for the system $(NSRW)$ with initial data in $B^{\frac{d}{2}}_{2,\infty}\cap B^{\frac{d}{2}-1}_{2,\infty}$ remains open and is probably false in general.}
  \end{rem}
As for the compressible Navier-Stokes system, in the $(NSRW)$ model the density fluctuation has two distinct behaviours in some low and high frequencies, separated by a frequency threshold. This naturally leads to the definition of the hybrid Besov spaces, involving two different regularities for low and high frequencies, introduced in \cite{CH} and defined for $l_\ee=[\frac{1}{2}\log_2(\frac{\gamma}{C_0\ee^2})-1]$ ($\gamma$ is a constant and $C_0=\frac{8}{3}$) and $s,t\in \R$ by :
\begin{equation}
 \|q\|_{\dot{B}_{\ee}^{s,t}} \overset{def}{=} \Sum_{l\leq l_\ee} 2^{ls} \|\ddl q\|_{L^2}+ \Sum_{l> l_\ee} \frac{1}{\ee^2} 2^{lt} \|\ddl q\|_{L^2}.
\label{normhybride1}
\end{equation}

\begin{defi}
\sl{(\cite{CH}) The space $E_{\ee}^s$ is the set of functions $(q,u)$ in
$$
\left(\cC_b(\R_+, \dot{B}_{2,1}^{s-1}\cap \dot{B}_{2,1}^s)\cap L^1(\R_+, \dot{B}_\ee^{s+1,s}\cap \dot{B}_\ee^{s+2,s})\right) \times
\left(\cC_b(\R_+, \dot{B}_{2,1}^{s-1})\cap L^1(\R_+, \dot{B}_{2,1}^{s+1})\right)^d
$$
endowed with the norm
\begin{multline}
\|(q,u)\|_{E_{\ee}^s} \overset{def}{=} \|u\|_{L^{\infty} \dot{B}_{2,1}^{s-1}}+ \|q\|_{L^{\infty} \dot{B}_{2,1}^{s-1}}+ \|q\|_{L^{\infty} \dot{B}_{2,1}^{s}}\\
+\|u\|_{L^1 \dot{B}_{2,1}^{s+1}}+ \|q\|_{L^1 \dot{B}_{\ee}^{s+1,s}}+ \|q\|_{L^1 \dot{B}_{\ee}^{s+2,s}}.
\end{multline}}
\end{defi}
We first state the global well-posedness for system $(RW_\ee)$ with uniform estimates with respect to $\ee$ (see \cite{CH}):
\begin{thm} (\cite{CH})
\sl{Let $\ee>0$ and assume that $\min(\mu,\nu=2\mu+\lambda)>0$. There exist two positive constants $\eta_R$ and $C$ only depending on $d$, $\kappa$, $\mu$, $\lambda$ and $P'(1)$ such that if $q_0\in \dot{B}_{2,1}^{\fd-1}\cap \dot{B}_{2,1}^{\fd}$, $u_0 \in \dot{B}_{2,1}^{\fd-1}$ and
$$
\|q_0\|_{\dot{B}_{2,1}^{\fd-1} \cap \dot{B}_{2,1}^{\fd}} +\|u_0\|_{\dot{B}_{2,1}^{\fd-1}}\leq \eta_R
$$
then system $(RW_\ee)$ has a unique global solution $(\rho_\ee, u_\ee)$ with $(q_\ee, u_\ee)\in E_{\ee}^{\fd}$ such that:
$$
\|(q_\ee,u_\ee)\|_{E_{\ee}^\fd} \leq C (\|q_0\|_{\dot{B}_{2,1}^{\fd-1} \cap \dot{B}_{2,1}^{\fd}} +\|u_0\|_{\dot{B}_{2,1}^{\fd-1}}).
$$}
\label{thexist}
\end{thm}
\begin{rem}
 \sl{Note that in the low frequency regime ($l\leq l_\ee\sim -\frac{1}{2}\log \ee$), the parabolic regularization for $q_\ee$ is the same as for the Korteweg system, indeed the low frequencies of $q_\ee$ are in $L_t^{1}(\dot{B}_{2,1}^{\fd+1} \cap \dot{B}_{2,1}^{\fd+2})$}.
\end{rem}

The main result in \cite{CH} is the following: when the initial data are small enough (so that we have global solutions for $(K)$ and $(RW_\ee)$) the solution of $(RW_\ee)$ goes to the solution of $(K)$ when $\ee$ goes to zero.

\begin{thm} (\cite{CH})
\sl{Assume that $\min(\mu,2\mu+\lambda)>0$, $P'(1)>0$ and that $q_0\in \dot{B}_{2,1}^{\fd-1}\cap \dot{B}_{2,1}^{\fd}$, $u_0 \in \dot{B}_{2,1}^{\fd-1}$. There exists $0<\eta\leq \min(\eta_K, \eta_R)$ such that if
$$
\|q_0\|_{\dot{B}_{2,1}^{\fd-1} \cap \dot{B}_{2,1}^{\fd}} +\|u_0\|_{\dot{B}_{2,1}^{\fd-1}}\leq \eta,
$$
then systems $(K)$ and $(RW_\ee)$ both have global solutions and $\|(\qe-q, \ue-u)\|_{E_\ee^{\fd}}$ tends to zero as $\ee$ goes to zero. Moreover, with the same notations as before, there exists a constant $C=C(\eta, \kappa, P'(1))>0$ such that for all $\alpha\in ]0, 1[$ (if $d=2$) or $\alpha\in ]0,1]$ (if $d\geq 3$), and for all $t\in\R_+$,
$$
 \|(\qe-q, \ue-u)\|_{E_\ee^{\fd-\alpha}} \leq C \ee^{\alpha},
$$}
\label{thcv}
\end{thm}
This results relies on the following estimates:
\begin{prop} (\cite{CH}, Proposition $1$)
\sl{Let $\ee>0$, $s\in \R$, $I=[0,T[$ or $[0, +\infty[$ and $v\in L^1(I,\dot{B}_{2,1}^{\fd+1}) \cap L^2 (I,\dot{B}_{2,1}^{\fd})$. Assume that $(q,u)$ is a solution of System $(LR_\ee)$ (see below) defined on $I$. There exists a constant $C>0$ depending on $d$, $s$, $\mu$, $\nu$, $p$, $\kappa$ such that for all $t\in I$,
\begin{multline}
 \|u\|_{\tilde{L}_t^{\infty} \dot{B}_{2,1}^{s-1}}+ \|q\|_{\tilde{L}_t^{\infty} \dot{B}_{2,1}^{s-1}}+ \|q\|_{\tilde{L}_t^{\infty} \dot{B}_{2,1}^{s}}+ \|u\|_{\tilde{L}_t^1 \dot{B}_{2,1}^{s+1}}+ \|q\|_{\tilde{L}_t^1 \dot{B}_{\ee}^{s+1,s}}+ \|q\|_{\tilde{L}_t^1 \dot{B}_{\ee}^{s+2,s}}\\
\leq C e^{C\int_0^t (\|\nabla v(\tau)\|_{\dot{B}_{2,1}^\fd}+ \|v(\tau)\|_{\dot{B}_{2,1}^\fd}^2)d\tau} \Big(\|u_0\|_{\dot{B}_{2,1}^{s-1}}+ \|q_0\|_{\dot{B}_{2,1}^{s-1}} + \|q_0\|_{\dot{B}_{2,1}^{s}}\\
+ \|F\|_{\tilde{L}_t^1 \dot{B}_{2,1}^{s-1}}+ \|F\|_{\tilde{L}_t^1 \dot{B}_{2,1}^{s}}+ \|G\|_{\tilde{L}_t^1 \dot{B}_{2,1}^{s-1}}\Big).
\end{multline}}
\label{oldapriori}
\end{prop}
\begin{rem}
 \sl{ Let us point out that the main difficulty consists in obtaining the previous accurate estimates in Besov space  depending on the parameter $\ee$ for the linear system  $(LR_\ee)$. Indeed it is quite tricky to deal with the convection terms (let us point out that it is absolutely necessary to integrate the convection term $v\cdot \n q$ in  $(LR_\ee)$ in order not to loose regularity on the density in the remainder term $F$, indeed $(RW_\ee)$ does not provide any regularizing effects on the density in high frequencies), in \cite{CH} we use energy methods and symmetrizers.}
 \end{rem}
The goal of this paper is to propose a more robust method which allows us to precisely keep track of the dependance with respect to the physical coefficients: viscosity and capillarity.
 \subsection{Statement of the results}

Classically, as in the study of compressible Navier-Stokes systems-type in critical spaces (see \cite{Dinv, CD, arma}), for proving the theorems \ref{thexist} and \ref{thcv} (see \cite{CH} section 2) the key point  consists in obtaining a priori estimates on the following advected linear system ($\ee>0$ is fixed and for more simplicity we write $(q,u)$ instead of $(\qe,\ue)$):

$$
\begin{cases}
\begin{aligned}
&\d_t q+ v.\nabla q+ \div u= F,\\
&\d_t u+ v.\nabla u -\cA u+ p\nabla q-\frac{\kappa}{\ee^2} \nabla(\phi_\ee*q-q)= G.\\
\end{aligned}
\end{cases}
\leqno{(LR_\ee)}
$$
With
$$\cA u= \mu \Delta u+ (\lambda+\mu)\nabla \div u.$$

As mentioned before, in the present article we are interested in obtaining, from a different point of view, these energy estimates by using a Lagrangian change of coordinate. Our method provides a more precise dependency of the various constants with respect to $p=P'(1)$ and the well-known ratio $\nu^2/4\kappa$ (that also appears in the local system, and in any evanescent viscosity-capillarity limit. In the 1D case, we refer to \cite{CH1} for the study of this limit (when $\frac{\nu^2}{4\kappa}=O(1)$). This result is proven via the introduction of an effective velocity).

\begin{rem}
 \sl{In Proposition \ref{normhybride2} we provide two equivalent (and more accurate) expressions for the hybrid norm $\|.\|_{\dot{B}_{\ee}^{s,s+2}}$ and we will only use them:
\begin{multline}
 \|f\|_{\dot{B}_\ee^{s+2,s}} \sim \|\frac{\phi_\ee*f-f}{\ee^2}\|_{\dot{B}_{2,1}^s} \sim \sum_{j\in \Z} \frac{1-e^{-c\ee^2 2^{2j}}}{\ee^2}  2^{js} \|\ddj f\|_{L^2}\\
\sim \sum_{j\in \Z} \min(\frac{1}{\ee^2}, 2^{2j}) 2^{js} \|\ddj f\|_{L^2}
\end{multline}
}
\end{rem}
Consequently the non-local capillary term $\frac{\phi_\ee*\nabla q_\ee-\nabla q_\ee}{\ee^2}$ has in fact the same regularity as the local capillary term $\nabla \Delta q$ : both of them belong to $L_t^1 (\dot{B}_{2,1}^{\fd-2}\cap \dot{B}_{2,1}^{\fd-1})$.\\

Let us now give the main result of the present article, which is a sharper version of Proposition \ref{oldapriori} (see the previous section):
\begin{thm} Let $\ee>0$, $-\fd+1<s<\fd+1$, $I=[0,T[$ or $[0, +\infty[$ and $v\in L^1(I,\dot{B}_{2,1}^{\fd+1}) \cap L^2 (I,\dot{B}_{2,1}^{\fd})$. Assume that $(q,u)$ is a solution of System $(LR_\ee)$ defined on $I$. There exists $\ee_0>0$, a constant $C>0$ depending on $d$, $s$ such that if $\ee\leq \ee_0$, for all $t\in I$ (denoting $\nu=\mu+2\lambda$ and $\nu_0=\min(\nu, \mu)$),
\begin{multline}
 \|u\|_{\tilde{L}_t^{\infty} \dot{B}_{2,1}^{s-1}}+ \|q\|_{\tilde{L}_t^{\infty} \dot{B}_{2,1}^{s-1}}+ \nu\|q\|_{\tilde{L}_t^{\infty} \dot{B}_{2,1}^{s}}+ \nu_0\|u\|_{\tilde{L}_t^1 \dot{B}_{2,1}^{s+1}}+ \nu\|q\|_{\tilde{L}_t^1 \dot{B}_{\ee}^{s+1,s-1}}+ \nu^2\|q\|_{\tilde{L}_t^1 \dot{B}_{\ee}^{s+2,s}}\\
\leq C_{p,\frac{\nu^2}{4\kappa}} e^{\displaystyle{C_{p,\frac{\nu^2}{4\kappa}} C_{visc}\int_0^t (\|\nabla v(\tau)\|_{\dot{B}_{2,1}^\fd}+ \|v(\tau)\|_{\dot{B}_{2,1}^\fd}^2)d\tau}}\\
\times\Big(\|u_0\|_{\dot{B}_{2,1}^{s-1}}+ \|q_0\|_{\dot{B}_{2,1}^{s-1}} + \nu\|q_0\|_{\dot{B}_{2,1}^{s}} +\|F\|_{\tilde{L}_t^1 \dot{B}_{2,1}^{s-1}}+ \nu\|F\|_{\tilde{L}_t^1 \dot{B}_{2,1}^{s}}+ \|G\|_{\tilde{L}_t^1 \dot{B}_{2,1}^{s-1}}\Big).
\label{estimapriori}
\end{multline}
\label{apriori}
where
$$
\begin{cases}
 \displaystyle{C_{p,\frac{\nu^2}{4\kappa}}=C \max(\sqrt{p}, \frac{1}{\sqrt{p}}) \max(\frac{4\kappa}{\nu^2}, (\frac{\nu^2}{4\kappa})^2),}\\
\displaystyle{C_{visc}=\frac{1+|\lambda+\mu|+\mu+\nu}{\nu_0}+\max(1, \frac{1}{\nu^3}).}
\end{cases}
$$
\end{thm}
\begin{rem}
 \sl{The viscous coefficient $C_{visc}$ satisfies:
$$C_{visc}=
\begin{cases}
 \frac{1+2\nu}{\mu}+\max(1, \frac{1}{\nu^3}) & \mbox{If } \lambda+\mu>0,\\
\frac{1+2\mu}{\nu}+\max(1, \frac{1}{\nu^3}) & \mbox{If } \lambda+\mu\leq0
\end{cases}
$$
and when both viscosities are small, we simply have $C_{visc}\leq \max(1,\frac{1}{\nu_0^3})$.
}
\end{rem}

\begin{rem}
\sl{The norm $\|q\|_{\tilde{L}_t^1 \dot{B}_{\ee}^{s+1,s-1}}$ in the left-hand side has to be compared to the term $\|q\|_{\tilde{L}_t^1 \dot{B}_{\ee}^{s+1,s}}$ from proposition \ref{oldapriori}. Thanks to the last term in the left-hand side, there is no loss of regularity. We only chose to write it this way in reference to the more meaningful equivalent expression from \eqref{equivmin}.
}
\end{rem}

\begin{rem}
\label{rqref}
\sl{Let us give a few comments about the advantages of the Lagrangian method compared to the symetrizers techniques used in \cite{CH}:
\begin{itemize}
\item First we are able with this method to provide accurate estimates tracking the physical coefficients and especially the influence of the ratio $\frac{\nu^2}{4\kappa}$ (we also refer to \cite{DD} for the importance of this ratio in the local Korteweg system). Indeed it plays an important role when considering the vanishing viscosity-capillarity process. In the one dimensional case we proved in \cite{CH1} the global convergence of the classical Korteweg solutions to the global weak entropy solutions of the compressible Euler system when $\kappa= \nu^2$ and $\nu$ goes to zero. Let us also mention the works of Lax and Levermore (\cite{LL}) who consider the non-viscous case and prove a vanishing capillarity process from the KdV equation to the Burger equation in the context of dispersive shock solutions (the main tool is the inverse scattering theory). Lefloch also studies the KdV equation with viscosity, and shows that the previous ratio is critical for the convergence towards a weak entropy solution or to a dispersive shock solution (notice that this also seems to be observed in numerical simulations). In particular we expect our linear a priori estimates to be useful to study the vanishing process. This is the object of a future work.
\item Seconds, the coefficients of the last two terms in the left-hand side of \eqref{estimapriori} are $\nu$ and $\nu^2$, which in the setting $\kappa\sim \nu^2$ with $\nu$ small, are larger than the $\nu^3$ obtained in \cite{CH} (see the end of section $2.1$ therein).  
\item Another important feature of the present paper is that this method allows us to deal with Besov spaces constructed on general $L^r$ spaces with $r\neq 2$ in the spirit of \cite{CD}. Indeed, such results cannot be obtained with symmetrizers methods, which are by nature based on energy methods and scalar products in $L^2$. We refer to the appendix for results in this direction and to \cite{CH} for the compressible Navier-Stokes system.
\end{itemize}
}
\end{rem}

The paper is structured in the following way: instead of using energy methods and symmetrizers (see \cite{Dinv, CH}), we will first consider the linear system $(LR_\ee)$ without transport (this is the object of the second section, where we precisely study the different frequency regimes). In section $3$ after presenting the Langrangian change of coordinates, we provide equivalent expressions for the hybrid Besov norm, and perform the Lagrangian change of variable, as introduced by T. Hmidi in \cite{TH1} (in high frequency regime, transport terms prevent any direct use of the linear estimates).

As in \cite{Dlagrangien, CD, TH2, TH3, TH4} we get estimates on the advected linear system and for this we need to bound additional external force terms generated by the change of variable. This is done thanks to estimates dealing with the action of a lagrangian change of variables on frequency truncation, such as the one proved by Vishik (see \cite{Vishik}), we also refer to \cite{Dbook}. In this part we focus on the main difficulty of the present article which consists in estimating the commutator of our non-local capillary operator under the Lagrangian change of variables. We end this section by giving an extension of our estimates allowing results in Besov spaces defined on $L^r$ with $r\neq 2$. The appendix is devoted to give the main tools of the Littlewood-Paley theory and recall classical results on the Lagrangian flow.

\begin{rem}
 \sl{ Let us mention that another approach would consist in introducing an effective velocity as in \cite{arma} in order to diagonalize the system in a certain way and to cancel out the coupling between density and velocity. This is the object of another article.
}
\end{rem}

\begin{rem}
\sl{For general considerations about non-local operators we can refer to the work of Rohde and Yong (see \cite{Rohderadiation}) and the recent paper of Alibaud, Cifani and Jakobsen (see \cite{Alibaud}).
}
\end{rem}

\section{Linear estimates}

The aim of this section is to obtain linear estimates for the following system:
$$
\begin{cases}
\begin{aligned}
&\d_t q+ \div u= F,\\
&\d_t u-\cA u+ p\nabla q-\frac{\kappa}{\ee^2} \nabla(\phi_\ee*q-q)= G.\\
\end{aligned}
\end{cases}
\leqno{(L_\ee)}
$$
With
$$\cA u= \mu \Delta u+ (\lambda+\mu)\nabla \div u.$$
Let us state the frequency-localized result that we will use in this article:
\begin{prop}
 \sl{Let $\ee>0$, $s\in \R$, $I=[0,T[$ or $[0, +\infty[$. Assume that $(q,u)$ is a solution of System $(L_\ee)$ defined on $I$. There exists $\ee_0>0$, a constant $C>0$ depending on $d$, $s$, $c_0$ and $C_0$ such that if $\ee\leq \ee_0$, for all $t\in I$ (as usual $\nu_0=\min(\nu, \mu)$), and for all $j\in \Z$,
\begin{multline}
 \|\ddj u\|_{L_t^\infty L^2}+\nu_0 2^{2j} \|\ddj u\|_{L_t^1 L^2}+(1+\nu 2^j)\left(\|\ddj q\|_{L_t^\infty L^2} +\nu\min(\frac{1}{\ee^2}, 2^{2j})\|\ddj q\|_{L_t^1 L^2}\right)\\
\leq C \max(\sqrt{p}, \frac{1}{\sqrt{p}}) \max(\frac{4\kappa}{\nu^2}, (\frac{\nu^2}{4\kappa})^2)\\
\times \left( (1+\nu 2^j)\|\ddj q_0\|_{L^2} +\|\ddj u_0\|_{L^2} + (1+\nu 2^j)\|\ddj F\|_{L_t^1 L^2} +\|\ddj G\|_{L_t^1 L^2}\right)
\end{multline}
}
\label{estimlinloc}
\end{prop}

\subsection{Study of the eigenvalues}

As in \cite{Dinv} or \cite{CD} we first introduce the Helmholtz decomposition of $u$. If the pseudo-differential operator $\Lambda$ is defined by $\Lambda f=\mathcal{F}^{-1}(|.|\hat{f}(.))$, we set:
\begin{equation}
\begin{cases}
 v=\Lambda^{-1} \div u,\\
w=\Lambda^{-1} \curl u\\
\end{cases} 
\label{eqhelmoltz}
\end{equation}
then $u=-\Lambda^{-1}\nabla v+\Lambda^{-1}\div w$ and the system turns into:
$$
\begin{cases}
\begin{aligned}
&\d_t q+\Lambda v=F,\\
&\d_t v-\nu \Delta v -p \Lambda q +\frac{\kappa}{\ee^2}\Lambda(\phi_\ee*q-q)=\Lambda^{-1} \div G,\\
&\d_t w-\mu \Delta w=\Lambda^{-1} \curl G.
\end{aligned}
\end{cases}
\leqno{(L_\ee')}
$$
The last equation is a heat equation, easily dealt thanks to classical heat estimates in Besov spaces (we refer to \cite{Dbook} chapter 2), so we can focus on the first two lines and compute the eigenvalues and eigenvectors of the matrix associated to the Fourier transform of this new system:
$$
\begin{cases}
\begin{aligned}
&\d_t \hat{q}+|\xi| \hat{v}=\hat{F},\\
&\d_t \hat{v}+\nu |\xi|^2 \hat{v} -p |\xi| \hat{q} +\frac{\kappa}{\ee^2}|\xi|(e^{-\ee^2 |\xi|^2}-1)\hat{q}=\frac{i}{|\xi|}\xi.\hat{G}.\\
\end{aligned}
\end{cases}
$$
As the external forces appear through homogeneous pseudo-differential operators of degree zero, we can compute the estimates in the case $F=G=0$ and deduce the general case from the Duhamel formula. So we finally study the following system:
$$
\d_t\left(\begin{array}{c}\hat{q}\\ \hat{v}\end{array}\right)
=A(\xi)\left(\begin{array}{c}\hat{q}\\\hat{v}\end{array}\right)\quad\hbox{with}\quad
A(\xi):= \left(\begin{array}{cc}0&-|\xi|\\p|\xi|+\frac{\kappa}{\ee^2}|\xi|(1-e^{-\ee^2|\xi|^2})&-\nu|\xi|^2\end{array}\right).
$$
The discriminant of the characteristic polynomial of the matrix is:
$$
\Delta=|\xi|^2\left(\nu^2|\xi|^2-4\left(p+\frac{\kappa}{\ee^2}(1-e^{-\ee^2|\xi|^2})\right)\right),
$$
and thanks to the variations of function $f_\ee:x\mapsto \nu^2 x-4\left(p+\frac{\kappa}{\ee^2}(1-e^{-\ee^2x})\right)$, we obtain the existence of a unique threshold $x_\ee>0$ such that
$$
\Delta(\xi)\begin{cases}
<0 \mbox{ if } |\xi|^2<x_\ee,\\
>0 \mbox{ if } |\xi|^2>x_\ee.
\end{cases}
$$
We emphasize that when $\frac{\nu^2}{4K}\geq 1$, $f_\ee$ is an increasing function on $\R_+$, and when $\frac{\nu^2}{4K}< 1$, $f_\ee$ is decreasing in $[0, \frac{1}{\ee^2}\log(\frac{4K}{\nu^2})]$ and then increasing.
\begin{prop}
 \sl{Under the same assumptions, we have:
$$
x_\ee \underset{\ee\rightarrow 0}{\sim}
\begin{cases}
\vspace{0.2cm}
\displaystyle{\frac{4p}{\nu^2-4\kappa}} & \mbox{ if }\frac{\nu^2}{4\kappa}>1,\\
\vspace{0.2cm}
\displaystyle{\frac{1}{\ee}\sqrt{\frac{2p}{\kappa}}} & \mbox{ if }\frac{\nu^2}{4\kappa}=1,\\
\displaystyle{\frac{a_{\kappa, \nu}}{\ee^2}} & \mbox{ if }\frac{\nu^2}{4\kappa}<1,\\
\end{cases}
$$
where $a_{\kappa, \nu}=a(\frac{\nu^2}{4\kappa})$ is the unique positive root of function $x\mapsto 1-\frac{4\kappa}{\nu^2}\frac{1-e^{-x}}{x}$. Moreover we have $\frac{4\kappa}{\nu^2}-1< a_{\kappa, \nu}< \frac{4\kappa}{\nu^2}$.
}
\end{prop}
\textbf{Proof:} First, as the function $h:x\mapsto \frac{1-e^{-x}}{x}$ decreases from $[0, \infty[$ to $]0,1]$, we easily prove that the following function:
\begin{equation}
 g_\ee(x)=\frac{f_\ee(x)}{\nu^2 x}=1-\frac{4p}{\nu^2 x}-\frac{4\kappa}{\nu^2}\frac{1-e^{-\ee^2x}}{\ee^2x}
\label{gepsilon}
\end{equation}
is an increasing function from $]0,\infty[$ to $]-\infty, 1[$, and for a fixed $x> 0$, it is increasing with respect to $\ee$ ($\forall x>0$, $\forall 0<\ee<\ee'$, $g_\ee(x)<g_{\ee'}(x)$), which implies that $\ee\mapsto x_\ee$ increases when $\ee$ decreases to zero.
\\
Then, computing $g_\ee(\frac{4\kappa}{\nu^2\ee^2})=e^{-\frac{4\kappa}{\nu^2}}-\frac{p\ee^2}{\kappa}$, which is positive when $\ee$ is small enough, as $g_\ee$ is increasing, we obtain that $x_\ee<\frac{4\kappa}{\nu^2\ee^2}$ when $\ee$ is small enough. Computing $g_\ee(\frac{4p}{\nu^2})=-\frac{4\kappa}{\nu^2}h(\ee^2\frac{4p}{\nu^2})<0$ we easily get that $x_\ee>\frac{4p}{\nu^2}$. Finally we obtain that:
\begin{itemize}
 \item If $\frac{\nu^2}{4\kappa}>1$, as for all $x$, $0<h(x)\leq 1$, we have $1-\frac{4\kappa}{\nu^2}-\frac{4p}{\nu^2 x_\ee}<g_\ee(x_\ee)=0\leq 1-\frac{4p}{\nu^2 x_\ee}$ which implies that the sequence $x_\ee$ is bounded (and therefore convergent because monotonous), so that $\ee^2 x_\ee$ goes to zero, and thanks to the fact that $f_\ee(x_\ee)=0$ and $h(x)\underset{x\rightarrow 0}{\rightarrow}1$, we obtain that $x_\ee\underset{\ee\rightarrow 0}{\rightarrow} \frac{4p}{\nu^2-4\kappa}$.
\item If $\frac{\nu^2}{4\kappa}=1$, guided by the value of $g_\ee(\frac{C}{\ee})=1-\frac{4p\ee}{\nu^2 C}-h(\ee C)$, and thanks to function study and Taylor expansions, we get that:
$$
\sqrt{\frac{2p}{\kappa}}\frac{1+\ee^2}{\ee}<x_\ee< \sqrt{\frac{2p}{\kappa}}\frac{1+\sqrt{\ee}}{\ee}.
$$
\item If $\frac{\nu^2}{4\kappa}<1$, guided by the value of $g_\ee(\frac{C}{\ee^2})$, we obtain that for sufficently small $\ee$:
$$
\frac{a(\frac{\nu^2}{4\kappa})}{\ee^2}<x_\ee< \frac{a(\frac{\nu^2}{4\kappa})+\ee}{\ee^2},
$$
where $a(\frac{\nu^2}{4\kappa})$ is the unique positive root of function $x\mapsto 1-\frac{4\kappa}{\nu^2}\frac{1-e^{-x}}{x}$ and satisfies $\frac{4\kappa}{\nu^2}-1< a_{\kappa, \nu}< \frac{4\kappa}{\nu^2}$. $\blacksquare$
\end{itemize}
Then we can compute the expressions of $\hat{q}$ and $\hat{v}$:
\\

\textbf{-For the low frequencies ($\Delta<0$)}, when $|\xi|<\sqrt{x_\ee}$, we have:
$$
\begin{cases}
\vspace{0.2cm}
 \hat{q}(\xi)=\frac{1}{2}\left((1+\frac{i}{S(\xi)})e^{t\lambda_+} +(1-\frac{i}{S(\xi)})e^{t\lambda_-}\right)\hat{q_0}(\xi)-i\frac{e^{t\lambda_+}-e^{t\lambda_-}}{\nu|\xi|S(\xi)}\hat{v_0}(\xi),\\
\hat{v}(\xi)=i\left(p+\frac{\kappa}{\ee^2}(1-e^{-\ee^2|\xi|^2})\right)\frac{e^{t\lambda_+}-e^{t\lambda_-}}{\nu|\xi|S(\xi)} \hat{q_0}(\xi)+\frac{1}{2}\left((1-\frac{i}{S(\xi)})e^{t\lambda_+} +(1+\frac{i}{S(\xi)})e^{t\lambda_-}\right)\hat{v_0}(\xi),
\end{cases}
$$
with:
\begin{equation}
S(\xi)=\sqrt{-g_\ee(|\xi|^2)}=\sqrt{\frac{4}{\nu^2 |\xi|^2}\left(p+\frac{\kappa}{\ee^2}(1-e^{-\ee^2|\xi|^2})\right)-1}
\label{formulelow}
\end{equation}
and
$$
\lambda_{\pm}=-\frac{\nu|\xi|^2}{2}(1\pm i S(\xi)).
$$

\textbf{-For the high frequencies ($\Delta>0$)}, when $|\xi|>\sqrt{x_\ee}$, we have:
$$
\begin{cases}
\vspace{0.2cm}
 \hat{q}(\xi)=\frac{1}{2}\left((1-\frac{1}{R(\xi)})e^{t\lambda_+} +(1+\frac{1}{R(\xi)})e^{t\lambda_-}\right)\hat{q_0}(\xi)+\frac{e^{t\lambda_+}-e^{t\lambda_-}}{\nu|\xi|R(\xi)}\hat{v_0}(\xi),\\
\hat{v}(\xi)=-\left(p+\frac{\kappa}{\ee^2}(1-e^{-\ee^2|\xi|^2})\right)\frac{e^{t\lambda_+}-e^{t\lambda_-}}{\nu|\xi|R(\xi)} \hat{q_0}(\xi)+\frac{1}{2}\left((1+\frac{1}{R(\xi)})e^{t\lambda_+} +(1-\frac{1}{R(\xi)})e^{t\lambda_-}\right)\hat{v_0}(\xi),
\end{cases}
$$
with:
\begin{equation}
R(\xi)=\sqrt{g_\ee(|\xi|^2)}=\sqrt{1-\frac{4}{\nu^2 |\xi|^2}\left(p+\frac{\kappa}{\ee^2}(1-e^{-\ee^2|\xi|^2})\right)}
\label{formulehigh}
\end{equation}
and
$$
\lambda_{\pm}=-\frac{\nu|\xi|^2}{2}(1\pm R(\xi)).
$$
\begin{rem}
 \sl{It will be crucial for the time integration to observe that
$$
p+\frac{\kappa}{\ee^2}(1-e^{-\ee^2|\xi|^2})=\frac{\nu^2|\xi|^2}{4}(1-R(\xi))(1+R(\xi)).
$$
}
\label{remvelocity}
\end{rem}

\begin{rem}
 \sl{In the low frequency regime, both eigenvalues provide parabolic heat regularization. In the high frequency regime, when $|\xi|$ is large, $\lambda_+ \sim -\nu |\xi|^2$ (which generates parabolic regularization), but $\lambda_-\sim -(p+\frac{\kappa}{\ee^2})$, that only provides a damping. This is the same behaviour as for the compressible Navier-Stokes system, and we refer to \cite{Dinv}, \cite{Dbook} and \cite{CD} (for a precise computation of the Fourier transform of the solutions).
}
\end{rem}

\begin{rem}
 \sl{These results have to be compared to the same study for the local Korteweg system, in this case the matrix becomes:
$$
\left(\begin{array}{cc}0&-|\xi|\\|\xi|(p+\kappa |\xi|^2)&-\nu|\xi|^2\end{array}\right).
$$
And the discriminant of the characteristic polynomial, $\Delta=|\xi|^2\left((\nu^2-4\kappa)|\xi|^2-4p\right)$, satisfies:
\begin{itemize}
 \item If $\frac{\nu^2}{4\kappa}\leq1$, $\Delta\leq0$ and $\lambda_\pm=\frac{-\nu|\xi|^2\pm i|\xi|\sqrt{4p-(\nu^2-4\kappa)|\xi|^2}}{2}$ (parabolic regularization everywhere),
\item If $\frac{\nu^2}{4\kappa}>1$, $\Delta<0$ then $\Delta>0$ with threshold $\frac{4p}{\nu^2-4\kappa}$, and in the high frequency regime, when $|\xi|$ is large, $\lambda_-\sim -\frac{\nu |\xi|^2}{2}(1-\sqrt{1-\frac{4\kappa}{\nu^2}})$ which, contrarily to the previous cases, generates parabolic regularization.
\end{itemize}
}
\end{rem}

\subsection{Pointwise estimates}

\subsubsection{Thresholds}

As seen previously, we recall that $x_\ee$ satisfies:
$$
x_\ee \geq
\begin{cases}
\vspace{0.2cm}
\displaystyle{\frac{4p}{\nu^2}} & \mbox{ if }\frac{\nu^2}{4\kappa}>1,\\
\vspace{0.2cm}
\displaystyle{\frac{1}{\ee}\sqrt{\frac{2p}{\kappa}}} & \mbox{ if }\frac{\nu^2}{4\kappa}=1,\\
\displaystyle{\frac{a(\frac{\nu^2}{4\kappa})}{\ee^2}} & \mbox{ if }\frac{\nu^2}{4\kappa}<1.\\
\end{cases}
$$
But in the following we will be able to get a frequency threshold of size $1/\ee^2$ in each of the three cases, for that we will not directly use $x_\ee$ as the announced threshold but a larger term, $y_\ee$, defined the following way (we refer to \eqref{gepsilon} for the definition of $g_\ee$):
\begin{equation}
y_\ee=
\begin{cases}
\vspace{0.2cm}
\displaystyle{g_\ee^{-1}(\frac{1}{4}) \mbox{ if } M<\frac{3}{4},}\\
\displaystyle{g_\ee^{-1}(1-\frac{1}{2M}) \mbox{ if } M\geq\frac{3}{4},}
      \end{cases}
\mbox{where }M=\frac{\nu^2}{4\kappa}.
 \end{equation}
As $g_\ee$ is an increasing function from $]0, \infty[$ to $]-\infty,1[$, we have $x_\ee<y_\ee$. The following property shows that $y_\ee$ is large enough:
\begin{prop}
 \sl{Under the previous estimates, there exist two constants $0<\gamma_1<\gamma_2$ and $\ee_0$ such that for all $0<\ee<\ee_0$,
$$
x_\ee<\frac{\gamma_1}{\ee^2}\leq y_\ee \leq \frac{\gamma_2}{\ee^2}.
$$
Moreover, if $M=\frac{\nu^2}{4\kappa}\geq \frac{3}{4}$, then $\gamma_1$ and $\gamma_1$ are universal constants, and if $M=\frac{\nu^2}{4\kappa}< \frac{3}{4}$ then we have:
$$
\frac{1}{3}\leq \frac{1}{M}-1<a(M)<\gamma_1<\gamma_2<\frac{3}{M}.
$$
}
\label{propgamma}
\end{prop}
\textbf{Proof:} Here it will be more efficient to compare $\frac{\nu^2}{4\kappa}$ to $\frac{3}{4}$ instead of $1$:

\textbf{First case:} If $M=\frac{\nu^2}{4\kappa}\geq \frac{3}{4}$ we define $\gamma_1$ and $\gamma_2$ by:
$$
h(\gamma_1)=\frac{1}{2}\mbox{ and }h(\gamma_2)=\frac{1}{4}.
$$
As $g_\ee(\frac{C}{\ee^2})=1-\frac{4p}{\nu^2 C}\ee^2-\frac{1}{M}h(C)$, we obtain that
$$
g_\ee(\frac{\gamma_1}{\ee^2})=1-\frac{1}{2M}-\frac{4p}{\nu^2 \gamma_1}\ee^2<1-\frac{1}{2M}=g_\ee(y_\ee)
$$
and if $\ee$ is small enough (depending on $p,\nu,\gamma_2$ and $M$),
$$
g_\ee(\frac{\gamma_2}{\ee^2})=1-\frac{1}{4M}-\frac{4p}{\nu^2 \gamma_2}\ee^2>g_\ee(y_\ee)
$$

\textbf{Second case:} If $M=\frac{\nu^2}{4\kappa}<\frac{3}{4}$, $\gamma_1$ and $\gamma_2$ are defined by:
$$
h(\gamma_1)=\frac{3}{4}M\mbox{ and }h(\gamma_2)=\frac{1}{2}M.
$$
Like previously, we compute:
$$
g_\ee(\frac{\gamma_1}{\ee^2})=1-\frac{1}{4}-\frac{4p}{\nu^2 \gamma_1}\ee^2<\frac{1}{4}=g_\ee(y_\ee),
$$
and if $\ee$ is small enough (depending on $p,\nu$ and $\gamma_2$),
$$
g_\ee(\frac{\gamma_2}{\ee^2})=1-\frac{1}{2}-\frac{4p}{\nu^2 \gamma_2}\ee^2>g_\ee(y_\ee),
$$
In both cases we conclude using that function $g_\ee$ is increasing.

In the second case, to get the lower estimate we use that $h(\gamma_1)<M=h(a(M))$ to obtain that (remind that the function $h$ is decreasing) $a(M)<\gamma_1$. We conclude thanks to the fact that $a(M) \geq \frac{1}{M}-1\geq \frac{1}{3}$. $\blacksquare$

\subsubsection{Estimates}

Now that we have defined the frequency thresholds $x_\ee$ and $y_\ee$ we can state the main result of this section:

\begin{prop}
 \sl{Under the previous notations, there exists a constant $C$, such that for all $j\in \Z$ and all $\xi\in 2^j \cC$ where $\cC$ is the annulus $\{\xi\in \R^d, c_0=\frac{3}{4}\leq |\xi|\leq C_0=\frac{8}{3}\}$, we have the following estimates (we denote by $f_j=\ddj f$ and we refer to the appendix for details on the Littlewood-Paley theory):
\begin{itemize}
 \item If $|\xi|<\sqrt{x_\ee}$:
$$
\begin{cases}
\vspace{0.2cm}
 (1+\nu 2^j)|\hat{q_j}(\xi)|\leq C e^{-\frac{\nu t c_0^2 2^{2j}}{4}} \left((1+\nu 2^j)|\hat{q_{0,j}}(\xi)|+(1+\frac{1}{\sqrt{p}})|\hat{v_{0,j}}(\xi)|\right),\\
|\hat{v_j}(\xi)|\leq C e^{-\frac{\nu t c_0^2 2^{2j}}{4}} \left((1+\nu 2^j)(1+\sqrt{p})(1+\frac{4\kappa}{\nu^2})|\hat{q_{0,j}}(\xi)|+|\hat{v_{0,j}}(\xi)|\right).
\end{cases}
$$
\item If $\sqrt{x_\ee}<|\xi|<\sqrt{y_\ee}$:
$$
\begin{cases}
\vspace{0.2cm}
 (1+\nu 2^j)|\hat{q_j}(\xi)|\leq \frac{C}{1-m} e^{-\frac{\nu t c_0^2 2^{2j}}{4}(1-m)} \left((1+\nu 2^j)|\hat{q_{0,j}}(\xi)|+(1+\frac{1}{\sqrt{p}})|\hat{v_{0,j}}(\xi)|\right),\\
|\hat{v_j}(\xi)|\leq \frac{C}{1-m} e^{-\frac{\nu t c_0^2 2^{2j}}{4}(1-m)} \left((\nu 2^j)|\hat{q_{0,j}}(\xi)|+|\hat{v_{0,j}}(\xi)|\right),
\end{cases}
$$
where $m=\frac{1}{2}$ if $M=\frac{\nu^2}{4\kappa}<\frac{3}{4}$, $m=\sqrt{1-\frac{1}{2M}}$ if $M\geq \frac{3}{4}$.
 \item If $|\xi|>\sqrt{y_\ee}$:
$$
\begin{cases}
\vspace{0.2cm}
 (1+\nu 2^j)|\hat{q_j}(\xi)|\leq C \left(e^{-\frac{\nu t |\xi|^2}{2}} +e^{-\frac{\kappa}{\nu \ee^2}(1-e^{-\gamma_1})t}\right) \left((1+\nu 2^j)|\hat{q_{0,j}}(\xi)|+(1+\frac{1}{\sqrt{p}})|\hat{v_{0,j}}(\xi)|\right),\\
|\hat{v_j}(\xi)|\leq C \left(e^{-\frac{\nu t c_0^2 2^{2j}}{4}} +\big(1-\sqrt{g_\ee(c_0^2 2^{2j})}\big)e^{-\frac{\nu t c_0^2 2^{2j}}{2}\big(1-\sqrt{g_\ee(C_0^2 2^{2j})}\big)} \right)\left(\nu 2^j |\hat{q_{0,j}}(\xi)|+|\hat{v_{0,j}}(\xi)|\right).
\end{cases}
$$
\end{itemize}
\label{Estimxi}
}
\end{prop}
\textbf{Proof:} We split the study according to the two frequency thresholds:

\textbf{Low frequencies:} assume that $\xi\in 2^j \cC$ with $|\xi|<\sqrt{x_\ee}$, then the discriminant $\Delta<0$ and the system has two non-real conjugated eigenvalues, thanks to (\ref{formulelow}), we can write:
$$
\begin{cases}
\vspace{0.2cm}
 \hat{q_j}(\xi)=\frac{1}{2}\left((e^{t\lambda_+}+e^{t\lambda_-}) +i\frac{e^{t\lambda_+}-e^{t\lambda_-}}{S(\xi)})\right)\hat{q}_{0,j}(\xi)-i\frac{e^{t\lambda_+}-e^{t\lambda_-}}{\nu|\xi|S(\xi)}\hat{v}_{0,j}(\xi),\\
\hat{v_j}(\xi)=i\left(p+\frac{\kappa}{\ee^2}(1-e^{-\ee^2|\xi|^2})\right)\frac{e^{t\lambda_+}-e^{t\lambda_-}}{\nu|\xi|S(\xi)} \hat{q}_{0,j}(\xi)+\frac{1}{2}\left((e^{t\lambda_+}+e^{t\lambda_-}) +i\frac{e^{t\lambda_-}-e^{t\lambda_+}}{S(\xi)}\right)\hat{v}_{0,j}(\xi),
\end{cases}
$$
with:
$$
S(\xi)=\sqrt{-g_\ee(|\xi|^2)}=\sqrt{\frac{4}{\nu^2 |\xi|^2}\left(p+\frac{\kappa}{\ee^2}(1-e^{-\ee^2|\xi|^2})\right)-1}
$$
and
$$
\lambda_{\pm}=-\frac{\nu|\xi|^2}{2}(1\pm i S(\xi)).
$$
Here we have to cope with two difficulties:
\begin{enumerate}
 \item neutralize $S(\xi)^{-1}$ for frequencies close to $\sqrt{x_\ee}$, as $S(\xi)$ goes to zero as $|\xi|$ goes to $\sqrt{x_\ee}$,
\item  if we are not careful, the first term in the expression of the velocity will provide either $\ee^{-2}$, either $2^{-2j}$ which would make our estimates useless.
\end{enumerate}
The first point will be adressed by considering the following block:
$$
\frac{e^{t\lambda_+}-e^{t\lambda_-}}{S(\xi)} =-e^{-\frac{\nu t |\xi|^2}{2}} \frac{e^{i\frac{\nu t|\xi|^2}{2} S(\xi)}-e^{-i\frac{\nu t|\xi|^2}{2} S(\xi)}}{2i\frac{\nu t|\xi|^2}{2}S(\xi)} i\nu t|\xi|^2 =-i\nu t|\xi|^2 e^{-\frac{\nu t |\xi|^2}{2}} \frac{\sin(\frac{\nu t|\xi|^2}{2}S(\xi))}{\frac{\nu t|\xi|^2}{2}S(\xi)}.
$$
As for all $x\geq 0$, $xe^{-x}\leq \frac{2}{e}e^{-x/2}$, we obtain that:
\begin{equation}
|\frac{e^{t\lambda_+}-e^{t\lambda_-}}{S(\xi)}|\leq \frac{4}{e} e^{-\frac{\nu t |\xi|^2}{4}}.
\label{blockestim}
\end{equation}
We then deduce:
\begin{equation}
|\hat{q_j}(\xi)|\leq C e^{-\frac{\nu t c_0^2 2^{2j}}{4}} \left(|\hat{q_{0,j}}(\xi)|+\frac{1}{\nu |\xi|}|\hat{v_{0,j}}(\xi)|\right).
\label{estimqlow1}
\end{equation}
And finally:
\begin{equation}
\nu 2^j |\hat{q_j}(\xi)|\leq C e^{-\frac{\nu t c_0^2 2^{2j}}{4}} \left(\nu 2^j|\hat{q_{0,j}}(\xi)|+|\hat{v_{0,j}}(\xi)|\right).
\label{estimqlow2}
\end{equation}
For the velocity, we can also use \eqref{blockestim}, except for the first term (that provides large coefficients described as difficulty 2 in the previous page.):
$$
A=i\left(p+\frac{\kappa}{\ee^2}(1-e^{-\ee^2|\xi|^2})\right)\frac{e^{t\lambda_+}-e^{t\lambda_-}}{\nu|\xi|S(\xi)}.
$$
We can observe that $S(\xi)\geq \sqrt{\frac{4p}{\nu^2 |\xi|^2}-1}=\frac{1}{\nu |\xi|}\sqrt{4p-\nu^2|\xi|^2}$, which makes sense only when $|\xi|\leq 2\sqrt{p}/\nu$. Therefore we split the study into two cases:
\begin{itemize}
 \item If $|\xi|\leq \sqrt{2p}/\nu$ then $\nu |\xi|S(\xi)\geq \sqrt{4p-\nu^2 |\xi|^2}\geq \sqrt{2p}$ and as $1-e^{-x}\leq x$, we get:
$$
|A|\leq \frac{p+\kappa|\xi|^2|}{\sqrt{2p}}(|e^{t\lambda_+}|+|e^{t\lambda_-}|)\leq \frac{p+\kappa\frac{2p}{\nu^2}}{\sqrt{2p}}2e^{-\frac{\nu t |\xi|^2}{2}}.
$$
Using this together with \eqref{blockestim}, we get:
$$
|\hat{v_j}(\xi)|\leq C e^{-\frac{\nu t c_0^2 2^{2j}}{2}} \left(\sqrt{p}(1+\frac{4\kappa}{\nu^2})|\hat{q_{0,j}}(\xi)|+|\hat{v_{0,j}}(\xi)|\right).
$$
going back to the density we can estimate the last term, and in addition to (\ref{estimqlow2}), we obtain:
$$
|\hat{q_j}(\xi)|\leq C e^{-\frac{\nu t c_0^2 2^{2j}}{4}} \left(|\hat{q_{0,j}}(\xi)| +|\hat{v_{0,j}}(\xi)|\right).
$$
\item If $\sqrt{2p}/\nu<|\xi|\leq \sqrt{x_\ee}$, then we have $g_\ee(\frac{2p}{\nu^2})<g_\ee(|\xi|^2)< g_\ee(x_\ee)$, which implies:
$$
1<\frac{4}{\nu^2 |\xi|^2} \left(p+\frac{\kappa}{\ee^2}(1-e^{-\ee^2|\xi|^2})\right)< 2(1+\frac{4\kappa}{\nu^2})
$$
and we obtain for the velocity:
$$
|\hat{v_j}(\xi)|\leq C e^{-\frac{\nu t c_0^2 2^{2j}}{4}} \left((1+\frac{4\kappa}{\nu^2})\nu2^j|\hat{q_{0,j}}(\xi)|+|\hat{v_{0,j}}(\xi)|\right).
$$
and for the density fluctuation, thanks to (\ref{estimqlow1}), as $|\xi|>\sqrt{2p}/\nu$, we also have:
$$
|\hat{q_j}(\xi)|\leq C e^{-\frac{\nu t c_0^2 2^{2j}}{4}} \left(|\hat{q_{0,j}}(\xi)|+\frac{1}{\sqrt{2p}}|\hat{v_{0,j}}(\xi)|\right).
$$
\end{itemize}
Finally, gathering these estimates gives the first point of the proposition \ref{Estimxi}.

\begin{rem}
 \sl{Note that $g_\ee(\frac{2p}{\nu^2})<g_\ee(\frac{4p}{\nu^2})<0=g_\ee({x_\ee})$ so $\frac{4p}{\nu^2}<x_\ee$ in any case.
}
\end{rem}
\begin{rem}
 \sl{We emphasize that in the case of the compressible Navier-Stokes system (see \cite{CD}), such a distinction was not necessary as the analogous of $S(\xi)$ was $\sqrt{\frac{4p}{\nu^2 |\xi|^2}-1}$ and the threshold was located at $|\xi|^2=\frac{2p}{\nu^2}$.
}
\end{rem}

\textbf{High frequencies:} assume that $\xi\in 2^j \cC$ with $|\xi|>\sqrt{y_\ee}>\sqrt{x_\ee}$, then the discriminant $\Delta>0$ and the system has two real eigenvalues. Thanks to (\ref{formulehigh}) we write the localized quantities:
$$
\begin{cases}
\vspace{0.2cm}
 \hat{q_j}(\xi)=\frac{1}{2}\left((1-\frac{1}{R(\xi)})e^{t\lambda_+} +(1+\frac{1}{R(\xi)})e^{t\lambda_-}\right)\hat{q_{0,j}}(\xi)+\frac{e^{t\lambda_+}-e^{t\lambda_-}}{\nu|\xi|R(\xi)}\hat{v_{0,j}}(\xi),\\
\hat{v_j}(\xi)=-\left(p+\frac{\kappa}{\ee^2}(1-e^{-\ee^2|\xi|^2})\right)\frac{e^{t\lambda_+}-e^{t\lambda_-}}{\nu|\xi|R(\xi)} \hat{q_{0,j}}(\xi)+\frac{1}{2}\left((1+\frac{1}{R(\xi)})e^{t\lambda_+} +(1-\frac{1}{R(\xi)})e^{t\lambda_-}\right)\hat{v_{0,j}}(\xi),
\end{cases}
$$
with:
\begin{equation}
R(\xi)=\sqrt{g_\ee(|\xi|^2)}=\sqrt{1-\frac{4}{\nu^2 |\xi|^2}\left(p+\frac{\kappa}{\ee^2}(1-e^{-\ee^2|\xi|^2})\right)}
\end{equation}
and
$$
\lambda_{\pm}=-\frac{\nu|\xi|^2}{2}(1\pm R(\xi)).
$$
We have $e^{t\lambda_+}\leq e^{-\frac{\nu t|\xi|^2}{2}}$ and as explained, when $|\xi|$ is large we cannot hope to get parabolic regularization from $\lambda_-$. But the definition of $y_\ee$ implies that in this case
$$
g_\ee(|\xi|^2)\geq g_\ee(y_\ee)=
\begin{cases}
\vspace{0.2cm}
 \frac{1}{4}\mbox{ if }M=\frac{\nu^2}{4\kappa}<\frac{3}{4},\\
1-\frac{1}{2M}\geq \frac{1}{3}\mbox{ if }M\geq\frac{3}{4},
\end{cases}
$$
therefore $\frac{1}{2}\leq R(\xi)\leq 1$, $1+R(\xi)\in[1,2]$ and $1-R(\xi)\in[0,\frac{1}{2}]$. More precisely, as $|\xi|^2\geq y_\ee\geq \frac{\gamma_1}{\ee^2}$, we have:
$$
1-R(\xi)=\frac{1-R(\xi)^2}{1+R(\xi)}=\frac{\frac{4}{\nu^2 |\xi|^2}\left(p+\frac{\kappa}{\ee^2}(1-e^{-\ee^2|\xi|^2})\right)}{1+R(\xi)}\geq \frac{2}{\nu^2 |\xi|^2} \frac{\kappa}{\ee^2}(1-e^{-\gamma_1}).
$$
So that we can bound $\lambda_-=-\frac{\nu |\xi|^2}{2}(1-R(\xi))$ and $e^{t\lambda_-}\leq e^{-t\frac{\kappa}{\nu\ee^2}(1-e^{-\gamma_1})}$, which allows to estimate the density fluctuation the following way:
\begin{multline}
|\hat{q_j}(\xi)|\leq\frac{1}{2}\left((1+\frac{1}{R(\xi)})(e^{t\lambda_+} +e^{t\lambda_-})\right)|\hat{q_{0,j}}(\xi)|+\frac{e^{t\lambda_+}+e^{t\lambda_-}}{\nu|\xi|R(\xi)}|\hat{v_{0,j}}(\xi)|,\\
\leq C(e^{-\frac{\nu t|\xi|^2}{2}} +e^{-t\frac{\kappa}{\nu\ee^2}(1-e^{-\gamma_1})})\left(|\hat{q_{0,j}}(\xi)|+\frac{1}{\nu |\xi|}|\hat{v_{0,j}}(\xi)|\right).
\end{multline}
Using successively that $|\xi|\geq c_0 2^j$ and $|\xi|\geq \sqrt{x_\ee}\geq \frac{\sqrt{2p}}{\nu}$ yields the estimate given in proposition \ref{Estimxi}.
\\
Concerning the velocity, we have two main difficulties in this frequency domain:
\begin{itemize}
 \item Obtain parabolic regularization for the velocity, even when the "damping eigenvalue" $\lambda_-$ is involved.
\item As previously, estimating like before the first term provides coefficients $\ee^{-2}$ (unbounded in $\ee$) or $2^{2j}$ (too many derivatives).
\end{itemize}
In order to neutralize these large coefficients, we rewrite the velocity as in Remark \ref{remvelocity}:
$$
|\hat{v_j}(\xi)|\leq \mathbb{A}+\mathbb{B}+\mathbb{C}\mbox{ with }
\begin{cases}
\vspace{0.2cm}
\mathbb{A}=\frac{\nu |\xi|}{4}(1-R(\xi))(1+R(\xi)) \frac{e^{t\lambda_-}-e^{t\lambda_+}}{R(\xi)} |\hat{q_{0,j}}(\xi)|,\\
\vspace{0.2cm}
\mathbb{B}=\frac{1}{2}(1+\frac{1}{R(\xi)})e^{t\lambda_+} |\hat{v_{0,j}}(\xi)|,\\
\mathbb{C}=\frac{1}{2}(\frac{1}{R(\xi)}-1)e^{t\lambda_-} |\hat{v_{0,j}}(\xi)|.
\end{cases}
$$
We can easily estimate $\mathbb{B}$:
\begin{equation}
\mathbb{B}\leq C e^{-\frac{\nu t c_0^2 2^{2j}}{2}}|\hat{v_{0,j}}(\xi)|.
\label{estB}
\end{equation}
The other two terms involve $\lambda_-$ and have to be handled carefully as we wish to get $L^1$ in time estimates:
$$
\mathbb{C}\leq C (1-R(\xi))e^{-\frac{\nu t c_0^2 2^{2j}}{2}(1-R(\xi))}|\hat{v_{0,j}}(\xi)|.
$$
Using the fact that function $g_\ee$ is increasing, as $c_0 2^j\leq |\xi|\leq C_0 2^j$, we have:
$$
\sqrt{g_\ee(c_0^2 2^{2j})} \leq R(\xi)=\sqrt{g_\ee(|\xi|^2)}\leq \sqrt{g_\ee(C_0^2 2^{2j})}
$$
and
$$
1-\sqrt{g_\ee(C_0^2 2^{2j})} \leq 1-R(\xi)\leq 1-\sqrt{g_\ee(c_0^2 2^{2j})},
$$
so that we can write:
\begin{equation}
\mathbb{C}\leq C\left(1-\sqrt{g_\ee(c_0^2 2^{2j})}\right) e^{-\frac{\nu t c_0^2 2^{2j}}{2}\left(1-\sqrt{g_\ee(C_0^2 2^{2j})}\right)}|\hat{v_{0,j}}(\xi)|.
\label{estC}
\end{equation}
We estimate the last term by:
$$
\mathbb{A}\leq \mathbb{A}_1+\mathbb{A}_2 \mbox{ with }
\begin{cases}
\mathbb{A}_1=\frac{\nu |\xi|}{4}(1-R(\xi))(1+R(\xi)) \displaystyle{\frac{e^{-\frac{\nu t |\xi|^2}{2}(1-R(\xi))}}{R(\xi)}} |\hat{q_{0,j}}(\xi)|,\\
\mathbb{A}_2=\frac{\nu |\xi|}{4}(1-R(\xi))(1+R(\xi)) \displaystyle{\frac{e^{-\frac{\nu t |\xi|^2}{2}(1+R(\xi))}}{R(\xi)} |\hat{q_{0,j}}}(\xi)|,\\
\end{cases}
$$
$\mathbb{A}_2$ is easily dealt:
$$
\mathbb{A}_2 \leq C \nu 2^j e^{-\frac{\nu t c_0^2 2^{2j}}{2}} |\hat{q_{0,j}}(\xi)|,
$$
and finally we get:
$$
\mathbb{A}_1 \leq C \nu 2^j \left(1-\sqrt{g_\ee(c_0^2 2^{2j})}\right) e^{-\frac{\nu t c_0^2 2^{2j}}{2}\left(1-\sqrt{g_\ee(C_0^2 2^{2j})}\right)}|\hat{q_{0,j}}(\xi)|.
$$
Gathering these estimates together with estimates (\ref{estB}) and (\ref{estC}) leads to the announced result in this case.

\begin{rem}
 \sl{Though small, terms like $1-\sqrt{g_\ee(C_0^2 2^{2j})}$ will allow to get $L^1$ in time estimates with parabolic regularization.
}
\end{rem}

\textbf{Medium frequencies:} assume that $\xi\in 2^j \cC$ with $\sqrt{x_\ee}<|\xi|<\sqrt{y_\ee}$. We use the same formula as in the high frequency case. Due to the fact that $R(\xi)$ goes to zero when $|\xi|$ goes to $\sqrt{x_\ee}$, we write the localized density and velocity as in the low frequency case:
$$
\begin{cases}
\vspace{0.2cm}
 \hat{q_j}(\xi)=\frac{1}{2}\left(e^{t\lambda_+}+e^{t\lambda_-}+\displaystyle{\frac{e^{t\lambda_-}-e^{t\lambda_+}}{R(\xi)}}\right)\hat{q_{0,j}}(\xi)+\displaystyle{\frac{e^{t\lambda_+}-e^{t\lambda_-}}{\nu|\xi|R(\xi)}}\hat{v_{0,j}}(\xi),\\
\hat{v_j}(\xi)=\frac{\nu |\xi|}{4}(1-R(\xi))(1+R(\xi)) \displaystyle{\frac{e^{t\lambda_+}-e^{t\lambda_-}}{R(\xi)}} \hat{q_{0,j}}(\xi)+\frac{1}{2} \left(e^{t\lambda_+}+e^{t\lambda_-}+\displaystyle{\frac{e^{t\lambda_+}-e^{t\lambda_-}}{R(\xi)}}\right)\hat{v_{0,j}}(\xi),
\end{cases}
$$
In this frequency domain, we define $m$ such that:
$$
0=\sqrt{g_\ee(x_\ee)}<R(\xi)=\sqrt{g_\ee(|\xi|^2)}< m\overset{\mbox{def}}{=}\sqrt{g_\ee(y_\ee)}=
\begin{cases}
\vspace{0.2cm}
 \frac{1}{2}\mbox{ if }M=\frac{\nu^2}{4\kappa}<\frac{3}{4},\\
\sqrt{1-\frac{1}{2M}}\mbox{ if }M\geq\frac{3}{4},
\end{cases}
$$
As $m<1$ we have
\begin{equation}
 \begin{cases}
  1<1+R(\xi)<2\\
1-m<1-R(\xi)<1
 \end{cases}
\label{m}
\end{equation}
The aim is to get parabolic regularization and neutralize the possibly vanishing $R(\xi)$ near $\sqrt{x_\ee}$, by considering the following blocks:
$$
|\frac{e^{t\lambda_+}-e^{t\lambda_-}}{R(\xi)}|=\frac{e^{t\lambda_-}-e^{t\lambda_+}}{R(\xi)}=e^{t\lambda_-}\frac{1-e^{t(\lambda_+-\lambda_-)}}{R(\xi)}.
$$
As $\lambda_{\pm}=-\frac{\nu|\xi|^2}{2}(1\pm R(\xi))$ (see \ref{formulehigh}), we have: $\lambda_+-\lambda_-=-\nu |\xi|^2 R(\xi)$ and then (thanks to the variations of function $h$ and to (\ref{m})):
$$
|\frac{e^{t\lambda_+}-e^{t\lambda_-}}{R(\xi)}|=\nu t |\xi|^2 e^{t\lambda_-}\frac{1-e^{-\nu t |\xi|^2 R(\xi)}}{\nu t |\xi|^2 R(\xi)}\leq \nu t |\xi|^2 e^{t\lambda_-}\leq \nu t |\xi|^2 e^{-\frac{\nu t |\xi|^2}{2}(1-m)}.
$$
Finally, using once again that $xe^{-x}\leq \frac{2}{e}e^{-x/2}$, we obtain:
\begin{equation}
|\frac{e^{t\lambda_+}-e^{t\lambda_-}}{R(\xi)}|\leq \frac{C}{1-m} e^{-\frac{\nu t |\xi|^2}{4}(1-m)}.
\end{equation}
From this we easily get that:
\begin{equation}
 |\hat{q_j}(\xi)|\leq \frac{C}{1-m} e^{-\frac{\nu t c_0^2 2^{2j}}{4}(1-m)} \left(|\hat{q_{0,j}}(\xi)|+\frac{1}{\nu |\xi|}|\hat{v_{0,j}}(\xi)|\right).
\end{equation}
Using successively that $\xi\in 2^j \cC$ and $|\xi|\geq \sqrt{x_\ee}\geq \frac{\sqrt{2p}}{\nu}$ the estimate of proposition \ref{Estimxi}. As $\lambda_+ \geq \lambda_-$, with the same method, we obtain the corresponding estimate for the velocity. This ends the proof of the proposition. $\blacksquare$

\subsection{Time estimates}

As in the case of the compressible Navier-Stokes system (see \cite{CD} section $3.1$), due to the choice $c_0=3/4$ and $C_0=8/3$ (see the appendix), we can observe that there exist at most two indices $\underline{j}_0 =\overline{j}_0-1$ or $\underline{j}_0 =\overline{j}_0$  such that $\sqrt{y_\ee} \in 2^j[c_0, C_0]$ for $j\in\{\underline{j}_0,\overline{j}_0\}$.

The aim of this part is to prove the following result, which implies Proposition \ref{estimlinloc}.

\begin{prop}
 \sl{Under the same assumptions as in Proposition \ref{estimlinloc}, there exists a constant $C$ such that for all $j\in \Z$ (denoting $M=\frac{\nu^2}{4\kappa}$):
\begin{itemize}
 \item For all $j\leq \overline{j}_0$,
\begin{multline}
\|v_j\|_{L_t^\infty L^2}+ \nu 2^{2j}\|v_j\|_{L_t^1 L^2} +(1+\nu2^j)\left(\|q_j\|_{L_t^\infty L^2}+ \nu 2^{2j}\|q_j\|_{L_t^1 L^2}\right) \leq \\
C \max(\frac{1}{M},M^2) \left((1+\nu 2^j)(1+\sqrt{p})\|q_{0,j}\|_{L^2} +(1+\frac{1}{\sqrt{p}}) \|v_{0,j}\|_{L^2}\right),\\
\label{estimlowint}
\end{multline}
\item For all $j>\overline{j}_0$,
\begin{multline}
\|v_j\|_{L_t^\infty L^2}+ \nu 2^{2j}\|v_j\|_{L_t^1 L^2} +(1+\nu2^j)\left(\|q_j\|_{L_t^\infty L^2}+ \frac{\nu}{\ee^2}\|q_j\|_{L_t^1 L^2}\right) \leq \\
C \max(1,M) \left((1+\nu 2^j)\|q_{0,j}\|_{L^2} +(1+\frac{1}{\sqrt{p}}) \|v_{0,j}\|_{L^2}\right).\\
\label{estimhighint}
\end{multline}
\end{itemize}
}
\label{estimlinseuil}
\end{prop}
\begin{rem}
\sl{Notice that due to the separated results in low and high frequencies, the coefficients here are more precise than in Theorem \ref{apriori}}.
\end{rem}

\textbf{Proof:} As previously, three cases have to be considered:
\\

\textbf{High frequencies:} if $j\geq \overline{j}_0+1$ then $\sqrt{y_\ee}\notin2^j[c_0,C_0]$ and for all $\xi\in 2^j \cC$, we have $\frac{\gamma_1}{\ee^2}<\sqrt{y_\ee}<c_0 2^j \leq |\xi|\leq C_0 2^j$ so that we only use the high frequency case from the previous proposition. Integrating with respect to $\xi$ (we recall that the frequencies are localized in $2^j \cC$), we obtain thanks to the Plancherel formula:
$$
\begin{cases}
\vspace{0.2cm}
 (1+\nu 2^j)\|q_j\|_{L^2}\leq C \left(e^{-\frac{\nu t \gamma_1}{2\ee^2}} +e^{-\frac{\kappa}{\nu \ee^2}(1-e^{-\gamma_1})t}\right) \left((1+\nu 2^j)\|q_{0,j}\|_{L^2}+(1+\frac{1}{\sqrt{p}})\|v_{0,j}\|_{L^2}\right),\\
\|v_j\|_{L^2}\leq C \left(e^{-\frac{\nu t c_0^2 2^{2j}}{4}} +\big(1-\sqrt{g_\ee(c_0^2 2^{2j})}\big)e^{-\frac{\nu t c_0^2 2^{2j}}{2}\big(1-\sqrt{g_\ee(C_0^2 2^{2j})}\big)} \right)\left(\nu 2^j \|q_{0,j}\|_{L^2}+\|v_{0,j}\|_{L^2}\right).
\end{cases}
$$
This immediately implies that:
$$
\begin{cases}
\vspace{0.2cm}
 (1+\nu 2^j)\|q_j\|_{L_t^\infty L^2}\leq C \left((1+\nu 2^j)\|q_{0,j}\|_{L^2}+(1+\frac{1}{\sqrt{p}})\|v_{0,j}\|_{L^2}\right),\\
\|v_j\|_{L_t^\infty L^2}\leq C \left(\nu 2^j \|q_{0,j}\|_{L^2}+\|v_{0,j}\|_{L^2}\right).
\end{cases}
$$
The $L_t^1$-estimates require a little more work:
$$
\begin{cases}
\vspace{0.2cm}
 (1+\nu 2^j)\|q_j\|_{L_t^1 L^2}\leq C \left(\frac{2\ee^2}{\nu \gamma_1} +\frac{\nu \ee^2}{\kappa(1-e^{-\gamma_1})}\right) \left((1+\nu 2^j)\|q_{0,j}\|_{L^2}+(1+\frac{1}{\sqrt{p}})\|v_{0,j}\|_{L^2}\right),\\
\|v_j\|_{L_t^1 L^2}\leq C \left(\frac{4}{\nu c_0^2 2^{2j}} +\frac{2}{\nu c_0^2 2^{2j}} \frac{1-\sqrt{g_\ee(c_0^2 2^{2j})}}{1-\sqrt{g_\ee(C_0^2 2^{2j})}}\right) \left(\nu 2^j \|q_{0,j}\|_{L^2}+\|v_{0,j}\|_{L^2}\right).
\end{cases}
$$
Thanks to the definition of $\gamma_1$ (see Proposition \ref{propgamma}) we can write:
$$
M\frac{\gamma_1}{1-e^{-\gamma_1}}=
\begin{cases}
 2M\mbox{ if }M=\frac{\nu^2}{4\kappa}\geq \frac{3}{4},\\
\frac{4}{3}\mbox{ if }M=\frac{\nu^2}{4\kappa}<\frac{3}{4},
\end{cases}
$$
so that:
$$
\frac{2\ee^2}{\nu \gamma_1} +\frac{\nu \ee^2}{\kappa(1-e^{-\gamma_1})}= \frac{\ee^2}{\nu \gamma_1}\left( 2+ \frac{\nu^2}{\kappa} \frac{\gamma_1}{1-e^{-\gamma_1}}\right) \leq\frac{\ee^2}{\nu \gamma_1}(2+\frac{M}{4}).
$$
Proposition \ref{propgamma} also implies that $1/\gamma_1$ is bounded from above by a universal constant, so that we obtain:
\begin{equation}
 (1+\nu 2^j)\|q_j\|_{L_t^1 L^2}\leq C \frac{\ee^2}{\nu}\max(1,M) \left((1+\nu 2^j)\|q_{0,j}\|_{L^2}+(1+\frac{1}{\sqrt{p}})\|v_{0,j}\|_{L^2}\right).
\end{equation}
We now turn to the velocity: as function $g_\ee$ in increasing, we get that $g_\ee(c_0^2 2^{2j})\leq g_\ee(C_0^2 2^{2j})$ and then $\frac{1-\sqrt{g_\ee(c_0^2 2^{2j})}}{1-\sqrt{g_\ee(C_0^2 2^{2j})}}>1$. This term has to be bounded uniformly in $\ee$ and $j$. Thanks to the expression and variations of $g_\ee$ we first write:
\begin{multline}
 \frac{1-\sqrt{g_\ee(c_0^2 2^{2j})}}{1-\sqrt{g_\ee(C_0^2 2^{2j})}}= \frac{1-g_\ee(c_0^2 2^{2j})}{1-g_\ee(C_0^2 2^{2j})}\times \frac{1+\sqrt{g_\ee(C_0^2 2^{2j})}}{1+\sqrt{g_\ee(c_0^2 2^{2j})}}\\
\leq 4 \frac{C_0^2}{c_0^2} \frac{\displaystyle{\frac{p}{\kappa}+ \frac{1-e^{-\ee^2 c_0^2 2^{2j}}}{\ee^2}}}{\displaystyle{\frac{p}{\kappa}+ \frac{1-e^{-\ee^2 C_0^2 2^{2j}}}{\ee^2}}} \leq 4 \frac{C_0^2}{c_0^2},
\end{multline}
as function $x\mapsto 1-e^{-x}$ is increasing. This yields:
$$
\|v_j\|_{L_t^\infty L^2}+ \nu 2^{2j}\|v_j\|_{L_t^1 L^2} \leq C \left(\nu 2^j \|q_{0,j}\|_{L^2}+\|v_{0,j}\|_{L^2}\right),
$$
so that estimate (\ref{estimhighint}) immediately follows.
\\

\textbf{low frequencies:} when $j\leq \underline{j}_0-1$ we know that $\sqrt{y_\ee}\notin [c_0 2^j, C_0 2^j]$ and here we will have to consider both cases $|\xi|<x_\ee$ and $x_\ee<|\xi|<y_\ee$:
$$
\|q_j\|_{L^2} \leq C \left(\|\hat{q_j} \textbf{1}_{\{|\xi|<\sqrt{x_\ee}\}}\|_{L^2} +\|\hat{q_j} \textbf{1}_{\{\sqrt{x_\ee}<|\xi|<\sqrt{y_\ee}\}}\|_{L^2} \right).
$$
Using the medium and low frequencies estimates from Proposition \ref{Estimxi} provides:
\begin{multline}
 (1+\nu 2^j) \|q_j\|_{L^2} \leq C \left(e^{-\frac{\nu t c_0^2 2^{2j}}{4}}+ \frac{1}{1-m}e^{-\frac{\nu t c_0^2 2^{2j}}{4}(1-m)} \right)\\
\times \left((1+\nu 2^j)\|q_{0,j}\|_{L^2}+(1+\frac{1}{\sqrt{p}})\|v_{0,j}\|_{L^2}\right),
\end{multline}
which implies:
\begin{multline}
 (1+\nu 2^j) \left(\|q_j\|_{L_t^\infty L^2}+ \nu 2^{2j}\|q_j\|_{L_t^1 L^2}\right) \leq\\
 \frac{C}{(1-m)^2} \left((1+\nu 2^j)\|q_{0,j}\|_{L^2}+(1+\frac{1}{\sqrt{p}})\|v_{0,j}\|_{L^2}\right).
\end{multline}
Let us recall that from Proposition \ref{propgamma} we have:
$$
1\leq \frac{1}{1-m} \leq
\begin{cases}
 2\mbox{ if } M=\frac{\nu^2}{4\kappa} \leq \frac{3}{4},\\
4M\mbox{ if } M>\frac{3}{4}
\end{cases}
\leq 4\max(1,M).
$$
The same can be done to the velocity and gives:
\begin{multline}
 \|v_j\|_{L_t^\infty L^2}+ \nu 2^{2j}\|v_j\|_{L_t^1 L^2}\leq\\
C\max(1,M^2) \left((1+\nu 2^j)(1+\sqrt{p})(1+\frac{1}{M})|\hat{q_{0,j}}(\xi)|+|\hat{v_{0,j}}(\xi)|\right).
\end{multline}
Which implies \ref{estimlowint} in the case $j\leq \underline{j}_0-1$.
\\

\textbf{Threshold frequencies:} in the case $j\in\{\underline{j}_0, \overline{j}_0\}$, we know that $\sqrt{y_\ee}\in 2^j[c_0,C_0]$ which may also be true for $\sqrt{x_\ee}$. Therefore we are forced to write (even if the first term of the right-hand side is zero when $M>1$...)
$$
\|q_j\|_{L^2} \leq C \left(\|\hat{q_j} \textbf{1}_{\{|\xi|<\sqrt{x_\ee}\}}\|_{L^2} +\|\hat{q_j} \textbf{1}_{\{\sqrt{x_\ee}<|\xi|<\sqrt{y_\ee}\}}\|_{L^2} +\|\hat{q_j} \textbf{1}_{\{|\xi|>\sqrt{y_\ee}\}}\|_{L^2} \right).
$$
Using all three estimates from Proposition \ref{Estimxi}, together with the specific bounds
$$
\frac{\gamma_1}{C_0^2}< \ee^2 2^{2j} < \frac{\gamma_2}{c_0^2},
$$
we end up with:
\begin{multline}
 (1+\nu 2^j) \|q_j\|_{L^2} \leq C \left( \frac{1}{1-m}e^{-\frac{\nu t c_0^2 2^{2j}}{4}(1-m)} + e^{-\frac{\kappa}{\nu}\frac{1-e^{-\gamma_1}}{\gamma_2}c_0^2 2^{2j}t}\right)\\
\times \left((1+\nu 2^j)\|q_{0,j}\|_{L^2}+(1+\frac{1}{\sqrt{p}})\|v_{0,j}\|_{L^2}\right),
\end{multline}
The rest of the proof follows the lines of the low frequencies case, except that for the time integral, we get:
\begin{multline}
 \nu 2^{2j} (1+\nu 2^j)\|q_j\|_{L_t^1 L^2}\leq\\
C \left(\frac{1}{(1-m)^2} +\frac{\nu^2}{4\kappa}\frac{\gamma_2 }{1-e^{-\gamma_1}}\right) \left((1+\nu 2^j)\|q_{0,j}\|_{L^2}+(1+\frac{1}{\sqrt{p}})|v_{0,j}|_{L^2}\right),
\end{multline}
From Proposition \ref{propgamma}, we have:
$$
\frac{\gamma_2 }{1-e^{-\gamma_1}}\leq
\begin{cases}
2\frac{\gamma_2}{\gamma_1}\mbox{ if }M\geq \frac{3}{4},\\ 
\frac{4}{3M}\frac{\gamma_2}{\gamma_1}\mbox{ if }M<\frac{3}{4},
\end{cases}
\mbox{ and }
\frac{\gamma_2}{\gamma_1}\leq
\begin{cases}
C\mbox{ if }M\geq \frac{3}{4},\\ 
\frac{3}{1-M}\leq 12 \mbox{ if }M<\frac{3}{4}.
\end{cases}
$$
which allows to write:
$$
\frac{1}{(1-m)^2} +\frac{\nu^2}{4\kappa}\frac{\gamma_2 }{1-e^{-\gamma_1}}\leq C\left(\max(1,M^2)+ M\max(1, \frac{1}{M})\right)\leq C \max(1,M^2).
$$
Doing the same frequency truncations for the velocity achieves the proof of the proposition. $\blacksquare$.

\subsection{End of the proof of Proposition \ref{estimlinloc}}

Let us recall that in \ref{eqhelmoltz}, $w$ satisfies a classical heat equation, then for all $j\in \Z$ (still in the homogenous assumption):
$$
\|w_j\|_{L_t^\infty L^2}+ \mu 2^{2j}\|w_j\|_{L_t^1 L^2} \leq C \|w_{0,j}\|_{L^2},
$$
and, together with Proposition \ref{estimlinseuil} (we recall that $\nu_0=\min(\nu, \mu)=\min(\mu+2\lambda, \mu)>0$),
\begin{itemize}
 \item For all $j\leq\overline{j}_0$,
\begin{multline}
\|u_j\|_{L_t^\infty L^2}+ \nu_0 2^{2j}\|u_j\|_{L_t^1 L^2} +(1+\nu2^j)\left(\|q_j\|_{L_t^\infty L^2}+ \nu \min(\frac{1}{\ee^2},2^{2j})\|q_j\|_{L_t^1 L^2}\right) \leq \\
C \max(\frac{1}{M},M^2) \left((1+\nu 2^j)(1+\sqrt{p})\|q_{0,j}\|_{L^2} +(1+\frac{1}{\sqrt{p}}) \|u_{0,j}\|_{L^2}\right),\\
\end{multline}
\item For all $j>\overline{j}_0$,
\begin{multline}
\|u_j\|_{L_t^\infty L^2}+ \nu_0 2^{2j}\|u_j\|_{L_t^1 L^2} +(1+\nu2^j)\left(\|q_j\|_{L_t^\infty L^2}+ \nu \min(\frac{1}{\ee^2},2^{2j})\|q_j\|_{L_t^1 L^2}\right) \leq \\
C \max(1,M) \left((1+\nu 2^j)\|q_{0,j}\|_{L^2} +(1+\frac{1}{\sqrt{p}}) \|u_{0,j}\|_{L^2}\right),\\
\end{multline}
\end{itemize}
which implies Proposition \ref{estimlinloc}. $\blacksquare$

\section{Proof of the advected linear estimates}

\subsection{Presentation of the difficulties}

The aim of this section is to prove Theorem \ref{apriori}. A natural idea is to use Proposition \ref{estimlinloc} and put the advection terms as external forces:  we can write that for all $j\in \Z$,
\begin{multline}
 \|\ddj u\|_{L_t^\infty L^2}+\nu_0 2^{2j} \|\ddj u\|_{L_t^1 L^2}+(1+\nu 2^j)\left(\|\ddj q\|_{L_t^\infty L^2} +\nu\min(\frac{1}{\ee^2}, 2^{2j})\|\ddj q\|_{L_t^1 L^2}\right)\leq \\
C_{p,\frac{\nu^2}{4\kappa}} \Bigg[ (1+\nu 2^j)\|\ddj q_0\|_{L^2} +\|\ddj v_0\|_{L^2} + (1+\nu 2^j)\left(\|\ddj F\|_{L_t^1 L^2}+ \|\ddj (v.\nabla q)\|_{L_t^1 L^2}\right)\\
+\left(\|\ddj G\|_{L_t^1 L^2} +\|\ddj (v.\nabla u)\|_{L_t^1 L^2}\right)\Bigg].
\label{energieBF}
\end{multline}
The next step consists in multiplying by $2^{j(s-1)}$ with $-\fd+1<s<\fd+1$ and sum over all frequencies. There are three additional terms to estimate. Thanks to paraproduct and remainder laws (we refer to (\ref{estimbesov}) in the appendix) we have:
\begin{multline}
 \|v.\nabla u\|_{\dot{B}_{2,1}^{s-1}}\leq \|\dot{T}_v \nabla u\|_{\dot{B}_{2,1}^{s-1}} +\|\dot{T}_{\nabla u}v\|_{\dot{B}_{2,1}^{s-1}}+\|\dot{R}(v,\nabla u)\|_{\dot{B}_{2,1}^{s-1}}\\
\leq C\left(\|v\|_{L^\infty} \|\nabla u\|_{\dot{B}_{2,1}^{s-1}} + \|\nabla u\|_{\dot{B}_{\infty,\infty}^{s-1-\fd}}\|v\|_{\dot{B}_{2,1}^\fd}+ \|v\|_{\dot{B}_{\infty,\infty}^0} \|\nabla u\|_{\dot{B}_{2,1}^{s-1}} \right)\leq C\|v\|_{\dot{B}_{2,1}^\fd} \|u\|_{\dot{B}_{2,1}^s}\\
\leq C \|v\|_{\dot{B}_{2,1}^\fd} \|u\|_{\dot{B}_{2,1}^{s-1}}^{\frac{1}{2}} \|u\|_{\dot{B}_{2,1}^{s+1}}^{\frac{1}{2}}\leq \frac{1}{2K} \nu_0\|u\|_{\dot{B}_{2,1}^{s+1}} +C^2\frac{K}{2\nu_0} \|v\|_{\dot{B}_{2,1}^\fd}^2 \|u\|_{\dot{B}_{2,1}^{s-1}},
\end{multline}
Consequently there exists a nonnegative summable sequence $\left(c_j(t)=c_j(u,v,t)\right)_{j\in \Z}$ whose sum is $1$ such that for all $j\in \Z$,
\begin{equation}
 \|\ddj (v.\nabla u)\|_{L_t^1 L^2} \leq 2^{-j(s-1)}  \int_0^t c_j(\tau)\left(\frac{1}{2K} \nu_0\|u\|_{\dot{B}_{2,1}^{s+1}} +C^2\frac{K}{2\nu_0}\|v\|_{\dot{B}_{2,1}^\fd}^2 \|u\|_{\dot{B}_{2,1}^{s-1}}\right) d\tau.
\label{advectBF1}
\end{equation}
Similarly we obtain thanks to Proposition \ref{estimhyb1}:
$$
\|v.\nabla q\|_{\dot{B}_{2,1}^{s-1}}\leq C \|v\|_{\dot{B}_{2,1}^\fd} \|q\|_{\dot{B}_{2,1}^s} \leq C \|v\|_{\dot{B}_{2,1}^\fd} (\|q\|_{\dot{B}_{2,1}^{s-1}} +\|q\|_{\dot{B}_{2,1}^{s}})^{\frac{1}{2}} (\|q\|_{\dot{B}_{\ee}^{s+1,s-1}} +\|q\|_{\dot{B}_{\ee}^{s+2,s}})^{\frac{1}{2}},
$$
and there exists a nonnegative summable sequence whose sum is $1$, once again denoted by $\left(c_j(t)=c_j(q,v,t)\right)_{j\in \Z}$, such that for all $j\in \Z$
\begin{multline}
 \|\ddj (v.\nabla q)\|_{L_t^1 L^2} \leq 2^{-j(s-1)} \int_0^t c_j(\tau)\Bigg( \frac{1}{2K} (\nu\|q\|_{\dot{B}_{\ee}^{s+1,s-1}} +\nu^2\|q\|_{\dot{B}_{\ee}^{s+2,s}})\\
+C^2\frac{K}{2} \max(1,\frac{1}{\nu^3}) \|v\|_{\dot{B}_{2,1}^\fd}^2 (\|q\|_{\dot{B}_{2,1}^{s-1}}+\nu \|q\|_{\dot{B}_{2,1}^s})\Bigg) d\tau.
\label{advectBF2}
\end{multline}
Writing the Bony decomposition,
$$
\|v.\nabla q\|_{\dot{B}_{2,1}^s}\leq \|\dot{T}_v \nabla q\|_{\dot{B}_{2,1}^s} +\|\dot{T}_{\nabla q}v\|_{\dot{B}_{2,1}^s}+\|\dot{R}(v,\nabla q)\|_{\dot{B}_{2,1}^s}
$$
Using (\ref{estimbesov}), we get that:
$$
\|\dot{T}_{\nabla q}v\|_{\dot{B}_{2,1}^s}+\|\dot{R}(v,\nabla q)\|_{\dot{B}_{2,1}^s} \leq C \|v\|_{\dot{B}_{2,1}^{\fd+1}} \|q\|_{\dot{B}_{2,1}^s} \leq C \|v\|_{\dot{B}_{2,1}^{\fd+1}} \frac{1}{\nu}\left(\|q\|_{\dot{B}_{2,1}^{s-1}}+\nu \|q\|_{\dot{B}_{2,1}^s}\right).
$$
which implies that there exists a nonnegative summable sequence whose sum is $1$, once again denoted by $(c_j(t))_{j\in\Z}$ such that:
\begin{equation}
 \nu 2^j \|\ddj (\dot{T}_{\nabla q}v +\dot{R}(v,\nabla q))\|_{L_t^1 L^2} \leq C 2^{-j(s-1)} \int_0^t c_j(\tau) \|v\|_{\dot{B}_{2,1}^{\fd+1}} \left(\|q\|_{\dot{B}_{2,1}^{s-1}} +\nu \|q\|_{\dot{B}_{2,1}^s}\right)d\tau.
\label{advectBF3}
\end{equation}
Unfortunately, we are not able to estimate this way $\dot{T}_v \nabla q$: there are too many derivatives involved in high frequency for the density fluctuation. But if we restrict to the low frequencies when $j\leq 0$ we can deal with this term: as the frequencies of $\dot{S}_{k-1}v. \ddk \nabla q$ are localized in an annulus of size $2^k$, there exists an integer $N_1$ (only depending on the parameters chosen in the Littlewood-Paley theory) such that:
$$
\ddj (\dot{T}_v \nabla q)=\sum_{|k-j|\leq N_1} \ddj(\dot{S}_{k-1}v. \ddk \nabla q)
$$
and for all $j\leq 0$,
$$
2^{js} \|\ddj (\dot{T}_v \nabla q)\|_{L^2} \leq 2^{js} \sum_{|k-j|\leq N_1} \|v\|_{L^\infty} 2^{k}\|\ddk q\|_{L^2}
$$
There exists a nonnegative summable sequence $(d_l(t)=d_l(q,t))_{l\in\Z_-}$ , with $\|d\|_{l^1(\Z_-)}\leq 1$ and a constant $C=C(N_1)$ such that for all $l\leq N_1$,
$$
\|\ddl q(t)\|_{L^2} \leq C d_l(t) 2^{-l(s+1)} \|q(t)\|_{\dot{B}_{2,1,N_1}^{s+1}},
$$
where the low frequency summation is well-defined because:
$$
\|q\|_{\dot{B}_{2,1,N_1}^{s+1}} =\sum_{l\leq N_1} 2^{l(s+1)} \|\ddl q\|_{L^2} \leq \sqrt{\sum_{l\leq N_1} 2^{l(s+2)} \|\ddl q\|_{L^2}} \sqrt{\sum_{l\leq N_1} 2^{ls} \|\ddl q\|_{L^2}}
$$
As we consider $j\leq 0$, we have $j-N_1\leq k \leq j+N_1\leq N_1 \leq -\log_2 \ee$ if $\ee$ is small enough, so that we can write
$$
\|\ddk q\|_{L^2} \leq C d_k(t) 2^{-k(s+1)} \|q\|_{\dot{B}_\ee^{s+2,s}}^{\frac{1}{2}} \|q\|_{\dot{B}_{2,1}^s}^{\frac{1}{2}},
$$
and then
$$
2^{js} \|\ddj (\dot{T}_v \nabla q)\|_{L^2} \leq C d_j'(t) \|v\|_{L^\infty} \|q\|_{\dot{B}_\ee^{s+2,s}}^{\frac{1}{2}} \|q\|_{\dot{B}_{2,1}^s}^{\frac{1}{2}},
$$
where $d_j'(t)=\sum_{|k-j|\leq N_1} d_k 2^{(j-k)s}$ is in $l^1(\Z_-)$ as a convolution of summable sequences (and its norm is bounded and only depends on $s$). Finally, we obtain that for all $j\leq 0$,
\begin{multline}
 \nu 2^j \|\ddj (\dot{T}_v \nabla q)\|_{L_t^1 L^2} \leq C 2^{-j(s-1)} \int_0^t d_j'(\tau) \nu \|v\|_{L^\infty} \|q\|_{\dot{B}_\ee^{s+2,s}}^{\frac{1}{2}} \|q\|_{\dot{B}_{2,1}^s}^{\frac{1}{2}} d\tau\\
\leq 2^{-j(s-1)}\int_0^t d_j'(\tau) \left(\frac{1}{2K} \nu^2\|q\|_{\dot{B}_\ee^{s+2,s}}+\frac{C^2 K}{2\nu} \|v\|_{L^\infty}^2 \nu \|q\|_{\dot{B}_{2,1}^s}\right) d\tau
\label{advectBF4}
\end{multline}
For all $j\in \Z$ and $t\in I$, let us introduce:
\begin{equation}
 U_j(t)=\|\ddj u\|_{L_t^\infty L^2}+\nu_0 2^{2j} \|\ddj u\|_{L_t^1 L^2}+(1+\nu 2^j)\left(\|\ddj q\|_{L_t^\infty L^2} +\nu\min(\frac{1}{\ee^2}, 2^{2j})\|\ddj q\|_{L_t^1 L^2}\right)
\label{Uj}
\end{equation}
and
\begin{multline}
 U(t)= \|u\|_{\tilde{L}_t^{\infty} \dot{B}_{2,1}^{s-1}}+ \|q\|_{\tilde{L}_t^{\infty} \dot{B}_{2,1}^{s-1}}+ \nu\|q\|_{\tilde{L}_t^{\infty} \dot{B}_{2,1}^{s}}+ \nu_0\|u\|_{\tilde{L}_t^1 \dot{B}_{2,1}^{s+1}}+ \nu\|q\|_{\tilde{L}_t^1 \dot{B}_{\ee}^{s+1,s-1}}+ \nu^2\|q\|_{\tilde{L}_t^1 \dot{B}_{\ee}^{s+2,s}}.
\end{multline}
Using estimates \eqref{advectBF1}, \eqref{advectBF2}, \eqref{advectBF3}, \eqref{advectBF4} in \eqref{energieBF} (and the fact that $\dot{B}_{2,1}^\fd \hookrightarrow L^\infty$, see appendix), we obtain that there exists a nonnegative summable sequence whose sum is $1$, denoted by $(c_j(t))_{j\in\Z}$ such that for all $j\leq 0$,
\begin{multline}
U_j(t) \leq C_{p,\frac{\nu^2}{4\kappa}} \Bigg[U_j(0) + (1+\nu 2^j)\|\ddj F\|_{L_t^1 L^2}+ \|\ddj G\|_{L_t^1 L^2}\\
+\frac{1}{2K} 2^{-j(s-1)}\int_0^t c_j(\tau) \left(\nu_0 \|u\|_{\dot{B}_{2,1}^{s+1}}+ \nu \|q\|_{\dot{B}_{\ee}^{s+1,s-1}}+ \nu^2\|q\|_{\dot{B}_{\ee}^{s+2,s}} \right) d\tau\\
+C^2\frac{K}{2} 2^{-j(s-1)}\int_0^t c_j(\tau) \left( \big(\max(1,\frac{1}{\nu^3}) +\frac{1}{\nu_0}\big) \|v(\tau)\|_{\dot{B}_{2,1}^\fd}^2
+C\|v(\tau)\|_{\dot{B}_{2,1}^{\fd+1}}\right) U(\tau) d\tau\Bigg].
\label{energieBF2}
\end{multline}
When $j>0$ we are not able to estimate the problematic term $\dot{T}_v \nabla q$ with this method. The key idea is to get rid of every advection term using a Lagrangian change of variable as in \cite{TH1} to \cite{TH4}, \cite{Dlagrangien} and \cite{CD}. This is the main difficulty of the article and the object of the following section.

\subsection{Lagrangian change of coordinates}

As stated before, the aim of this part is to get rid of the advection terms involved in system $(LR_\ee)$. Let us first consider the localized equations (as usual we set $f_j=\ddj f$...) written in the following way: 
$$
\begin{cases}
\begin{aligned}
&\d_t q_j+ \dot{S}_{j-1}v.\nabla q_j+ \div u_j= f_j,\\
&\d_t u_j+ \dot{S}_{j-1}v.\nabla u_j -\cA u_j+ p\nabla q_j-\frac{\kappa}{\ee^2} \nabla(\phi_\ee*q_j-q_j)= g_j,\\
\end{aligned}
\end{cases}
$$
where the external force terms are defined by:
$$
f_j= F_j +\left(\dot{S}_{j-1}v.\nabla q_j-\ddj (v.\nabla q)\right)\mbox{ and }g_j= G_j +\left(\dot{S}_{j-1}v.\nabla u_j-\ddj (v.\nabla u)\right).
$$
These two terms can be estimated thanks to the following commutator estimate from \cite{Dlagrangien} (we refer to lemma $B.1$ from appendix $B$):
\begin{lem} (\cite{Dlagrangien})
\sl{
There exists a sequence $(c_j)_{j\in \Z} \in l^1(\Z)$ such that $\|c\|_{l^1(\Z)}=1$ and a constant $C=C(d,\sigma)$ such that for all $j\in \Z$,
$$
\|\dot{S}_{j-1}v.\nabla h_j-\ddj (v.\nabla h)\|_{L^2} \leq C c_j 2^{-j\sigma}\|\nabla v\|_{\dot{B}_{2,\infty}^{\fd}\cap L^\infty} \|h\|_{\dot{B}_{2,1}^\sigma}
$$
}
\label{lemmeB}
\end{lem}
\begin{rem}
 \sl{Let us emphasize that we considered, up to adding external force terms that we can control ($f_j$ and $g_j$), the advection by the low frequencies of $v$, that is exactly the quantities that we were not able to deal with.}
\end{rem}
Let us set $\psi_{j,t}$ as the flow associated to $\dot{S}_{j-1}v$:
\begin{equation}
 \begin{cases}
 \partial_t \psi_{j,t}(x)=\dot{S}_{j-1}v(t,\psi_{j,t}(x))\\
\psi_{j,0}(x)=x.
\end{cases}
\label{flotdef}
\end{equation}
we can also write:
$$
\psi_{j,t}(x)=x+\int_0^t \dot{S}_{j-1}v(\tau,\psi_{j,\tau}(x)) d\tau.
$$
Thanks to propositions \ref{p:flow} and \ref{detjacobien} from the appendix (we refer to \cite{Dlagrangien} or \cite{CD}), there exists a constant $C$ such that:
\begin{equation}
 \begin{array}{lll}
\|g\circ\psi_{j,t}\|_{L^p}&\leq& e^{CV}\|g\|_{L^p}\quad\hbox{for all function }\ g\ \hbox{ in }\ 
L^p,\\[1ex]
\|D\psi_{j,t}^\pm\|_{L^\infty}&\leq&e^{CV},\\[1ex]
\|D\psi_{j,t}^\pm-I_d\|_{L^\infty}&\leq& e^{CV}-1,\\[1ex]
\|D^k\psi_{j,t}^\pm\|_{L^\infty}&\leq& C2^{(k-1)j}\Bigl(e^{CV}-1\Bigr)\quad\hbox{for }\ 
k\geq 2,
\end{array}
\label{estimflot}
\end{equation}
where
\begin{equation}
 V(t)\overset{\mbox{def}}{=}\int_0^t \|\nabla v(\tau)\|_{L^\infty}d\tau. 
\label{defV}
\end{equation}
Let us also emphasize that in the present case, the jacobian determinant of the change of variable will play a crucial role in the obtention of uniform estimates with respect to $\ee$ (contrary to the case of lemma $2.6$ from \cite{Dbook} where it produces a term that we are not able to sum):
\begin{equation}
\begin{cases}
det(D\psi_{j,t}(x))=e^{\int_{0}^t (\div \dot{S}_{j-1}v)(\tau,\psi_{j,\tau(x)})d\tau},\\
det(D\psi_{j,t}^{-1}(x))=e^{-\int_{0}^t (\div \dot{S}_{j-1}v)(\tau, X_j(\tau,t,x))d\tau} =e^{-\int_{0}^t (\div \dot{S}_{j-1}v)(\tau,\psi_{j,\tau}\circ \psi_{j,t}^{-1}(x))d\tau},
\end{cases}
\label{detjacobienj}
\end{equation}
where $X_j(\tau,t,x))$ denotes the two parameter flow associated to $\dot{S}_{j-1}v$ (we refer to \eqref{flot2param} in the appendix).

We will now perform the announced change of variable, for a function $h$, let us define $\tilde h =h \circ \psi_{j,t} =h(t,\psi_{j,t})$. Then we have $\partial_t \tilde q_j(t,x)=(\partial_t q_j+ \dot{S}_{j-1} v.\nabla q_j)(t, \psi_{j,t}(x))$, which provides the following system:
\begin{equation}
 \begin{cases}
\begin{aligned}
&\d_t \tilde q_j+ \div \tilde u_j= \tilde f_j +R_j^1,\\
&\d_t \tilde u_j -\cA \tilde u_j+ p\nabla \tilde q_j-\frac{\kappa}{\ee^2} \nabla(\phi_\ee*\tilde q_j-\tilde q_j)= \tilde g_j +R_j^2 +R_j^3 +\kappa R_j,\\
\end{aligned}
\end{cases}
\label{systchanged}
\end{equation}
where most of the remainder terms $R_q^1,$ $R_q^2$ and $R_q^3$ are exactly the same as in \cite{CD} (with the same convention: if $f:\R^d\rightarrow \R^m$ is a differentiable function then $Df$ denotes the Jacobian matrix of $f,$ and $\nabla f$ is the transposed matrix of $Df.$):
$$\displaylines{
R_j^1(t,x):={\rm Tr}\bigl(\nabla \tilde u_j(t,x)\cdot (I_d-\nabla\psi_{j,t}^{-1}(\psi_{j,t}(x)))\bigr),\cr
R_j^2(t,x):= \nabla \tilde q_j(t,x)\cdot (\nabla\psi_{j,t}^{-1}(\psi_{j,t}(x))-I_d)}
$$
and  $R_j^3:= \mu R_j^4 + (\lambda+\mu) R_j^5$ with
$$
\displaylines{
\quad R_j^{4,i}(t,x):={\rm Tr}\Bigl( (\nabla\psi_{j,t}^{-1}(\psi_{j,t}(x))-I_d)\cdot\nabla  D \tilde u_j^i(t,x) \cdot D\psi_{j,t}^{-1}(\psi_{j,t}(x))\hfill \cr\hfill
+\nabla D\tilde u_j^i(t,x)\cdot (D\psi_{j,t}^{-1}(\psi_{j,t}(x))-I_d)\Bigr) +\nabla \tilde u_j^i(t,x)\cdot 
\Delta \psi_{j,t}^{-1}(\psi_{j,t}(x)) \cr
\quad R_j^{5,k}(t,x):={\rm Tr}\Bigl( D \tilde u_j(t,x) \cdot \partial_k D\psi_{j,t}^{-1}(\psi_{j,t}(x))\Bigr)\hfill\cr\hfill
+ \Sum_{a,b,c, b\neq a, c\neq k} \partial_{bc}^2 \tilde u_j^i(t,x)\cdot \partial_k \psi_{j,t}^{-1,c}(\psi_{j,t}(x)) \cdot \partial_a \psi_{j,t}^{-1,b}(\psi_{j,t}(x))\hfill \cr\hfill
+ \Sum_{i=1}^d \partial_{ki}^2 \tilde u_j^i(t,x)\cdot \Bigl(\partial_k \psi_{j,t}^{-1,k}(t,\psi_{j,t}(x))-I_d) \cdot \partial_i \psi_{j,t}^{-1,i}(t,\psi_{j,t}(x))+ (\partial_i \psi_{j,t}^{-1,i}(\psi_{j,t}(x))-I_d) \Bigr).}
$$
There is only one additionnal remainder term compared to \cite{CD}:
\begin{equation}
R_j=\frac{\phi_\ee*\nabla q_j-\nabla q_j}{\ee^2}\circ \psi_{j,t} -\frac{\phi_\ee*\nabla \tilde q_j-\nabla \tilde q_j}{\ee^2}.
\label{Rj}
\end{equation}
Thanks to the definition of $\phi_\ee$ and $\phi$ we obtain that:
\begin{equation}
 \frac{\phi_\ee* f-f}{\ee^2}= \frac{1}{\ee^2}\int_{\R^d} \phi(z)\left(f(x-\ee z)-f(x)\right)dz
\label{phieps}
\end{equation}
and the difficulty comes from the fact that we need to obtain bounds in $L_t^1 L^2$ that are uniform with respect to $\ee$, that go to zero when $t$ is small (estimated by $e^{CV(t)}-1$), and that do not involve too many derivatives (that is $2^{\sigma j}$) for the density fluctuation in high frequencies. Remember that formally "$\frac{\phi_\ee* f- f}{\ee^2} \rightarrow \Delta f$" so we formally have "$R_j\rightarrow \Delta q_j \circ \psi_{j,t}-\Delta(q_j \circ \psi_{j,t})$": this term has no reason to be small in $\ee$ and the best we can hope is to get uniform bounds with respect to $\ee$.\\
\\
Dealing with $R_j$ is the object of the rest of this section.

\subsection{Precisions on the capillary term}

Before dealing with $R_j$ let us briefly go back in this section to the convolution term written in (\ref{phieps}): for a function $f$,
$$
\hat{\frac{\phi_\ee* f-f}{\ee^2}}(\xi)=-\frac{1-e^{-\ee^2|\xi|^2}}{\ee^2} \hat{f}(\xi).
$$
\begin{prop}
 \sl{For any suitable function $f$ and any $s\in\R$, the following two norms are equivalent:
\begin{equation}
\|f\|_{\dot{B}_\ee^{s+2,s}} \sim \|\frac{\phi_\ee*f-f}{\ee^2}\|_{\dot{B}_{2,1}^s}.
\label{equivmin}
\end{equation}
}
\label{normhybride2}
\end{prop}
\textbf{Proof:} When we consider a frequency localization of $f$: if $j\in\Z$,
$$
\hat{\frac{\phi_\ee* \ddj f-\ddj f}{\ee^2}}(\xi)=-\frac{1-e^{-\ee^2|\xi|^2}}{\ee^2} \hat{\ddj f}(\xi),
$$
and thanks to the Plancherel formula,
$$
\|\frac{\phi_\ee* \ddj f-\ddj f}{\ee^2}\|_{L^2}^2 =C \int_{2^j \cC} \left(\frac{1-e^{-\ee^2|\xi|^2}}{\ee^2}\right)^2 |\hat{\ddj f}(\xi)|^2 d\xi.
$$
Thanks again to the fact that on $\R_+$, $g:x\mapsto 1-e^{-x}$ is increasing and $h:x\mapsto(1-e^{-x})/x$ is decreasing, if $c_0 2^j\leq |\xi|\leq C_0 2^j$, we can write:
$$
\frac{1-e^{-c_0^2 \ee^2 2^{2j}}}{\ee^2}\leq \frac{1-e^{-\ee^2|\xi|^2}}{\ee^2}\leq \frac{1-e^{-C_0^2 \ee^2 2^{2j}}}{\ee^2},
$$
and
$$
\frac{c_0^2}{C_0^2}\frac{1-e^{-C_0^2 \ee^2 2^{2j}}}{\ee^2}\leq \frac{1-e^{-\ee^2|\xi|^2}}{\ee^2}\leq \frac{C_0^2}{c_0^2}\frac{1-e^{-c_0^2 \ee^2 2^{2j}}}{\ee^2}.
$$
therefore
$$
\|\ddj \frac{\phi_\ee*f-f}{\ee^2}\|_{L^2} \sim \frac{1-e^{-c_0^2 \ee^2 2^{2j}}}{\ee^2} \|\ddl f\|_{L^2}.
$$
\begin{rem}
 \sl{Obviously we can replace $c_0$ by any positive constant.}
\end{rem}
On the other hand, it will be useful to compare $\frac{1-e^{-c_0^2 \ee^2 2^{2j}}}{\ee^2}$ and $\min(\frac{1}{\ee^2}, 2^{2j})$.
\begin{itemize}
\item If $2^j\geq \frac{1}{\ee^2}$ then $2^j \ee^2 c_0^2 \geq c_0^2$ and as function $g$ is increasing and bounded from above by $1$:
\begin{equation}
 \frac{1-e^{-c_0^2}}{\ee^2} \leq \frac{1-e^{-c_0^2 \ee^2 2^{2j}}}{\ee^2}\leq \frac{1}{\ee^2} =\min(\frac{1}{\ee^2}, 2^{2j}).
\label{equivmin1}
\end{equation}
\item If $2^j\leq \frac{1}{\ee^2}$ then $2^j \ee^2 c_0^2 \leq c_0^2$ and as function $h$ is decreasing as for all $x\geq 0$, $0<h(x)\leq 1$ we have:
$$
\frac{1-e^{-c_0^2}}{c_0^2} \leq \frac{1-e^{-c_0^2 \ee^2 2^{2j}}}{\ee^2 c_0^2 2^{2j}}\leq 1,
$$
and
\begin{equation}
 (1-e^{-c_0^2})2^{2j} \leq \frac{1-e^{-c_0^2 \ee^2 2^{2j}}}{\ee^2}\leq c_0^2 2^{2j}= c_0^2 \min(\frac{1}{\ee^2}, 2^{2j}).
\label{equivmin2}
\end{equation}
\end{itemize}
From (\ref{equivmin1}) and \eqref{equivmin2}, we deduce that for all $j\in\Z$:
$$
(1-e^{-c_0^2})\min(\frac{1}{\ee^2}, 2^{2j}) \leq \frac{1-e^{-c_0^2 \ee^2 2^{2j}}}{\ee^2}\leq \max(1,c_0^2) \min(\frac{1}{\ee^2}, 2^{2j}),
$$
and
$$
\|\frac{\phi_\ee*\ddj f-\ddj f}{\ee^2}\|_{L^2} \sim \min(\frac{1}{\ee^2}, 2^{2j}) \|\ddl f\|_{L^2}.
$$
Multiplying by $2^{js}$ and summing over $j\in \Z$, we obtain that
$$
\|\frac{\phi_\ee*f-f}{\ee^2}\|_{\dot{B}_{2,1}^s} \sim \|f\|_{\dot{B}_\ee^{s+2,s}}. \quad\blacksquare
$$
\begin{rem}
 \sl{This new formulation seems more natural than \eqref{normhybride1} as instead of a fixed frequency threshold there is a continuous transition zone between the parabolically regularized frequencies and the damped frequencies}.
\end{rem}
\begin{rem}
 \sl{Notice that the norm obtained by replacing $\hat{\phi}(\xi)=e^{-|\xi|^2}$ by $e^{-\alpha|\xi|^2}$ is equivalent (the multiplicative constants depending on $\alpha$).}
\end{rem}
\begin{rem}
 \sl{Our space is a particular case of the general hybrid Besov spaces introduced by R. Danchin (see \cite{Dinv, Dbook}): we refer to the appendix for the fact that $\dot{B}_{\ee}^{s+2,s} \sim \tilde{B}_{\ee}^{s, \frac{2}{3}}$.
}
\end{rem}

\begin{rem}
 \sl{In the $L^r$-setting, we can prove that for all $j\in\Z$:
$$
\|\frac{\phi_\ee*\ddj f-\ddj f}{\ee^2}\|_{L^r} \leq C\min(\frac{1}{\ee^2}, 2^{2j}) \|\ddl f\|_{L^r}.
$$
}
\label{Lr}
\end{rem}

\subsection{Estimates on the capillary term}
\label{sectionRj}
In this section we wish to focus on the capillary term:
$$
R_j=\frac{\phi_\ee*\nabla q_j-\nabla q_j}{\ee^2}\circ \psi_{j,t} -\frac{\phi_\ee*\nabla \tilde q_j-\nabla \tilde q_j}{\ee^2}.
$$
Thanks to
$$
\nabla \tilde q_j =\nabla( q_j\circ \psi_{j,t})=\nabla q_j \circ \psi_{j,t} \times D\psi_{j,t} =\nabla q_j \circ \psi_{j,t} \times (D\psi_{j,t}-I_d) +\nabla q_j \circ \psi_{j,t},
$$
we obtain the following decomposition: $R_j= I_j+II_j$ with
\begin{equation}
\begin{cases}
\vspace{0.2cm}
I_j= \displaystyle{\frac{\phi_\ee*g_j-g_j}{\ee^2} \mbox{ where } g_j= \nabla q_j \circ \psi_{j,t} \times (I_d-D\psi_{j,t})},\\
II_j= \displaystyle{\frac{\phi_\ee*\nabla q_j-\nabla q_j}{\ee^2}\circ \psi_{j,t} -\frac{\phi_\ee*(\nabla q_j \circ \psi_{j,t}) -\nabla q_j \circ \psi_{j,t}}{\ee^2}}.
\end{cases}
\label{IjIIj}
\end{equation}
The main difficulty of this paper is then clearly established and consists in estimating (locally in frequency) the commutator between the Lagrangian change of variable and the non-local operator:
$$
 L_\ee(f)=\frac{\phi_\ee* f-f}{\ee^2}= \frac{1}{\ee^2}\int_{\R^d} \phi(z)\left(f(x-\ee z)-f(x)\right)dz.
$$
For a function $f$, and for all $j\in\Z$ we set $f_j=\ddj f$ and:
\begin{equation}
 II_j'= II_j'(f)= \frac{\phi_\ee*f_j-f_j}{\ee^2}\circ \psi_{j,t} -\frac{\phi_\ee*(f_j \circ \psi_{j,t}) -f_j \circ \psi_{j,t}}{\ee^2}.
\end{equation}
\begin{thm}
 \sl{Let $\sigma \in \R$. There exists a constant $C=C_{\sigma, d}$ such that for all $f\in\dot{B}_\ee^{\sigma+2, \sigma}$, there exists a summable positive sequence $(c_j(f))_{j\in \Z}$ whose sum is $1$ such that for all $t$ so small that
\begin{equation}
e^{2CV}-1 \leq \frac{1}{2}.
 \label{Cond1}
\end{equation}
and for all $j\in \Z$,
$$
\|II_j'(f)\|_{L^2}\leq C e^{CV}(V+e^{2CV}-1) c_j(f) 2^{-j\sigma}\|\frac{\phi_\ee*f-f}{\ee^2}\|_{\dot{B}_{2,1}^\sigma},
$$
where $V(t)=\int_0^t \|\nabla v(\tau)\|_{L^\infty}d\tau$. 
}
\label{II'}
\end{thm}
\begin{rem}
\sl{As a by-product we obtain that under the previous assumptions, if $t$ is small enough,
$$
\frac{1}{\ee^2} \|(\phi_\ee*\ddj f)\circ \psi_{j,t} -\phi_\ee*(\ddj f\circ \psi_{j,t})\|_{L^2} \leq C e^{CV}(V+e^{2CV}-1) c_j(f) 2^{-j\sigma}\|\frac{\phi_\ee*f-f}{\ee^2}\|_{\dot{B}_{2,1}^\sigma}
$$
even if neither of the left-hand side terms are spectrally localized.} 
\end{rem}
\textbf{Proof:} Let us first rewrite the non-local operator:
\begin{equation}
 L_\ee(f)(x)=\frac{\phi_\ee* f-f}{\ee^2}(x)= \frac{1}{\ee^2} \int_{\R^d} \phi_\ee(x-y)\left(f(y)-f(x)\right)dy.
\end{equation}
In the works of T. Hmidi, S. Keraani, H. Abidi, M. Zerguine (we refer to \cite{TH1}, \cite{TH2}, \cite{TH3} and \cite{TH4}) the key idea is to express the difference $II_j'$ as an integral formulation and try to retrieve the desired Besov norm thanks to an equivalent expression of this norm involving finite differences of $f$ of order $1$ (that is expressions of the type $\tau_{-y} f - f$ where $\tau_{-y} f(x)= f(x+y)$) or order $2$.
\begin{rem}
 \sl{Due to the equivalent forms of our hybrid norm, it would be useless here to simplify the expression of $II_j'(f)$ and write $II_j'= \frac{1}{\ee^2} \left( (\phi_\ee*f_j)\circ \psi_{j,t} -\phi_\ee*(f_j \circ \psi_{j,t})\right)$.}
\end{rem}
We refer to the following results of \cite{Dbook} for these useful alternative characterizations of Besov norms (we also refer to proposition $1.37$ with a simpler proof in the specific case of Sobolev spaces):
\begin{thm}\sl{(\cite{Dbook}, $2.36$)
 Let $s \in ]0,1[$ and $p,r\in [1,\infty]$. There exists a constant $C$ such that for any $u\in \cC_h'$,
$$
C^{-1} \|u\|_{\dot{B}_{p,r}^s}\leq \|\frac{\|\tau_{-y}u -u\|_{L^p}}{|y|^s}\|_{L^r (\R^d; \frac{dy}{|y|^d})} \leq C \|u\|_{\dot{B}_{p,r}^s}.
$$}
\end{thm}
and in the case where $s=1$,
\begin{thm}\sl{(\cite{Dbook}, $2.37$)
 Let $p,r\in [1,\infty]$. There exists a constant $C$ such that for any $u\in \cC_h'$,
$$
C^{-1} \|u\|_{\dot{B}_{p,r}^1}\leq \|\frac{\|\tau_{-y}u +\tau_y u -2u\|_{L^p}}{|y|}\|_{L^r (\R^d; \frac{dy}{|y|^d})} \leq C \|u\|_{\dot{B}_{p,r}^1}.
$$}
\end{thm}
\begin{rem}
 \sl{We emphasize that in the second case, when the regularity index $s$ is an integer, we have to use finite differences of order $2$ instead of order $1$. In the present article, due to the capillary term we will have to completely rewrite these results.}
\end{rem}
For example, in \cite{TH2}, the authors need to estimate in $L^p$ the commutator $|D|^\alpha (f\circ \psi)-(|D|^\alpha f)\circ \psi$ in terms of the $\dot{B}_{p,1}^\alpha$-norm of $f$ (with $\alpha \in ]0,1[$), where $\psi$ is the flow associated to a divergence-free vectorfield. To do so they write both terms in the commutator in a unified shape and manage to get the result thanks to a non-local expression of the fractionnal derivative:
$$
|D|^\alpha f(x)= C_\alpha \int_{\R^d} \frac{f(x)-f(y)}{|x-y|^{d+\alpha}}dy,
$$
and thanks to an equivalent expression of the Besov norm, using finite differences of order $1$. Let us emphasize that in this case, the vectorfield is divergence-free therefore the Jacobian determinant of $\psi$ is $1$, which simplifies the estimates.
\\

In our case, we have to estimate $II_j'$ which is a commutator between the non-local operator $L_\ee$ and the Lagrangian change of variable. Our operator is still non-local but the vectorfield is not incompressible anymore, and more important, we will have to construct an analogous of the finite difference expression of our hybrid Besov norm. As the Jacobian determinant of $\psi_{j,t}$ is not constant we cannot use the following estimate, used in \cite{TH1}, \cite{TH2}, which is due to Vishik (see \cite{Vishik} in the divergence-free case):
\begin{equation}
\|\ddl (\ddj a \circ \psi)\|_{L^p} \leq C 2^{-|j-l|} \|\nabla \psi^{sign(j-l)}\|_{L^\infty}\cdot \|\ddj a\|_{L^p}
\label{Vishincompr}
\end{equation}
In the non measure-preserving case we can recall the following result:
\begin{lem}(\cite{Dbook} lemma $2.6$)
 \sl{There exists a constant $C>0$ such that for all global diffeomorphism of $\R^d$ $\psi$, for all $p\geq 1$, for all tempered distribution $a$, and for all $j,l\in \Z$, we have
$$
\|\ddl (\ddj a \circ \psi)\|_{L^p} \leq C 2^{-j} \|J \psi^{-1}\|_{L^\infty} \|\ddj a\|_{L^p} (\|\nabla J\psi^{-1}\|_{L^\infty} \|J \psi\|_{L^\infty} + 2^l \|\nabla \psi\|_{L^\infty}).
$$}
\label{Vishcompr}
\end{lem}
Using this result in our non measure-preserving case involves an additionnal term (due to the Jacobian determinant) that we are not able to deal with (indeed, the first term in our case is not summable in $l$). On the other hand, we cannot hope to use a better estimate than the following one, that is small when $t$ is small but consumes too many derivatives:
$$
\|f_j(\psi_{j,t}(.))-f_j(.)\|_{L^2} \leq C V(t) 2^j \|f_j\|_{L^2}
$$
To bypass this difficulty we will simply write directly the expression in the Lagrangian variable and precisely trace the Jacobian determinant. Let us go back to our estimate:
\begin{multline}
 II_j'(x)= \frac{1}{\ee^2}\Big(\int_{\R^d} \phi_\ee(\psi_{j,t}(x)-y)\left(f_j(y)-f_j(\psi_{j,t}(x))\right)dy\\
-\int_{\R^d} \phi_\ee(x-y)\left(f_j(\psi_{j,t}(y))-f_j(\psi_{j,t}(x))\right)dy \Big).
\end{multline}
Let us now study $II_j'(\psi_{j,t}^{-1}(x))$ instead of $II_j'(x)$:
\begin{multline}
 II_j'(\psi_{j,t}^{-1}(x))= \frac{1}{\ee^2}\Big(\int_{\R^d} \phi_\ee(x-y)\left(f_j(y)-f_j(x)\right)dy\\
-\int_{\R^d} \phi_\ee(\psi_{j,t}^{-1}(x)-y)\left(f_j(\psi_{j,t}(y))-f_j(x)\right)dy \Big).
\end{multline}
Performing in the second integral the change of variable $y=\psi_{j,t}^{-1}(z)$, we obtain:
\begin{multline}
 \int_{\R^d} \phi_\ee(\psi_{j,t}^{-1}(x)-y)\left(f_j(\psi_{j,t}(y))-f_j(x)\right)dy\\
= \int_{\R^d} \phi_\ee(\psi_{j,t}^{-1}(x)-\psi_{j,t}^{-1}(z))\left(f_j(z)-f_j(x)\right) |det D\psi_{j,t}^{-1}(z)|dz
\end{multline}
Thanks to \eqref{detjacobienj} we have (recall that $X_j$ is the two-parameter flow associated to $\dot{S}_{j-1}v$: $X_j(\tau,t,z)=\psi_{j,\tau}\circ \psi_{j,t}^{-1}(z)$):
$$
det(D\psi_{j,t}^{-1}(z))=e^{\displaystyle{-\int_{0}^t (\div \dot{S}_{j-1}v)(\tau, X_j(\tau,t,z))d\tau}},
$$
so that
\begin{multline}
 II_j'(\psi_{j,t}^{-1}(x))= \frac{1}{\ee^2}\Big(\int_{\R^d} \phi_\ee(x-y)\left(f_j(y)-f_j(x)\right).\\
\left[1-\frac{\phi_\ee\big(\psi_{j,t}^{-1}(x)-\psi_{j,t}^{-1}(y)\big)}{\phi_\ee(x-y)} e^{\displaystyle{-\int_{0}^t (\div \dot{S}_{j-1}v)(\tau, X_j(\tau,t,y))d\tau}} \right]dy\Big).
\end{multline}
Recall that for all $x\in \R^d$, $\phi_\ee(x)=\frac{1}{(2\pi \ee)^d}e^{-\frac{|x|^2}{4\ee^2}}$ and perform the change of variable $z=\frac{x-y}{\ee}$, we obtain:
\begin{multline}
II_j'(\psi_{j,t}^{-1}(x))= \frac{1}{(2\pi)^d\ee^2}\Bigg(\int_{\R^d} e^{\displaystyle{-\frac{|y|^2}{4}}}\left(f_j(x-\ee y)-f_j(x)\right).\\
\left[1-e^{\displaystyle{\frac{|y|^2}{4}\left(1-\frac{|\psi_{j,t}^{-1}(x)-\psi_{j,t}^{-1}(x-\ee y)|^2}{\ee^2|y|^2}\right)}} e^{\displaystyle{-\int_{0}^t (\div \dot{S}_{j-1}v)(\tau, X_j(\tau,t,x-\ee y))d\tau}}\right]dy\Bigg).
\end{multline}
As we want to estimate the $L^2$ norm of this quantity, that is a Besov norm with integer regularity index $s=0$, the finite difference of order $1$ will not be sufficient for our need, and we will have to introduce finite differences of order $2$. Indeed, using the present quantity would only involve a term in $\ee|y|$  and when estimating in low frequencies, there would be either a multiplicative coefficient $1/\ee$ or an additionnal derivative term $2^{-j}$ (that would prevent any convergence when $-j$ is large). To be able to do a correct estimate we need at least $\ee^2|y|^2$.

To do this we simply write $II_j'=\frac{1}{2}(II_j'+II_j')$, and perform the change of variable $z=-y$ in the second integral. If we set:
\begin{equation}
 \begin{cases}
A= \displaystyle{\frac{|y|^2}{4}\left(1-\frac{|\psi_{j,t}^{-1}(x)-\psi_{j,t}^{-1}(x-\ee y)|^2}{\ee^2|y|^2}\right),}\\
B=\displaystyle{-\int_{0}^t (\div \dot{S}_{j-1}v)(\tau, X_j(\tau,t,x-\ee y))d\tau,}\\
C= \displaystyle{\frac{|y|^2}{4}\left(1-\frac{|\psi_{j,t}^{-1}(x)-\psi_{j,t}^{-1}(x+\ee y)|^2}{\ee^2|y|^2}\right),}\\
D=\displaystyle{-\int_{0}^t (\div \dot{S}_{j-1}v)(\tau, X_j(\tau,t,x+\ee y))d\tau,}
 \end{cases}
\label{ABCD}
\end{equation}
then
\begin{multline}
 II_j'(\psi_{j,t}^{-1}(x))= \frac{1}{2(2\pi)^d\ee^2}\Bigg(\int_{\R^d} e^{\displaystyle{-\frac{|y|^2}{4}}}\left(f_j(x-\ee y)-f_j(x)\right) \left[1-e^A e^B\right]dy\\
+\int_{\R^d} e^{\displaystyle{-\frac{|y|^2}{4}}}\left(f_j(x+\ee y)-f_j(x)\right) \left[1-e^C e^D\right]dy \Bigg)\\
= III_j+ IV_j.
\end{multline}
where
\begin{equation}
 \begin{cases}
  III_j= \displaystyle{\frac{1}{2(2\pi)^d\ee^2} \int_{\R^d} e^{\displaystyle{-\frac{|y|^2}{4}}}\left(f_j(x-\ee y) +f_j(x+\ee y) -2f_j(x)\right) \left[1-e^A e^B\right]dy,}\\
IV_j= \displaystyle{\frac{1}{2(2\pi)^d\ee^2} \int_{\R^d} e^{\displaystyle{-\frac{|y|^2}{4}}}\left(f_j(x+\ee y)-f_j(x)\right) \left[e^A e^B-e^C e^D\right]dy.}
 \end{cases}
\label{IIIetIV}
\end{equation}
\begin{rem}
 \sl{Even if there is only a finite difference of degree $1$ in the second term, we will be able to get correct estimates thanks to the coefficient $e^A e^B-e^C e^D$.}
\end{rem}
We will make an extensive use of the following elementary consequence of the mean-value theorem:
\begin{lem}
 \sl{For any $x,y\in \R$, $|e^x-e^y|\leq |x-y|e^{\max(x,y)}$.}
\label{lemexp}
\end{lem}
Let us begin with $III_j$:
\begin{prop}
 \sl{Under the previous assumptions, there exist a positive constant $C=C_{\sigma,d}$ and a nonnegative sequence $(c_j= c_j(f))_{j\in \Z}$ whose summation is $1$, such that if $t$ is so small that $e^{2CV(t)}-1\leq \frac{1}{2}$, we have:
$$
\|III_j\|_{L^2}\leq C (V+e^{2CV}-1) e^{CV} 2^{-j\sigma} c_j \|\frac{\phi_\ee* f-f}{\ee^2}\|_{\dot{B}_{2,1}^\sigma}.
$$
\label{propIII}
}
\end{prop}
To prove this result we will successively prove the following lemmas:
\begin{lem}
 \sl{There exists a constant $C$ such that for all $j\in\Z$, $f$, and all $t$ is so small that $e^{2CV(t)}-1\leq \frac{1}{2}$,
\begin{multline}
 \|III_j\|_{L^2}\leq \frac{1}{2(2\pi)^d} (V+e^{2CV}-1) e^{CV}\\
\times \frac{1}{\ee^2} \int_{\R^d} e^{-\frac{|y|^2}{16}}\|f_j(.-\ee y) +f_j(.+\ee y) -2f_j(.)\|_{L^2} dy,
\end{multline}
where $V$ is defined in \eqref{defV}.
}
\label{lemIIIa}
\end{lem}
\begin{lem}
 \sl{For all $\sigma\in \R$ there exists a constant $C_{\sigma,d}$ such that for any $f\in \dot{B}_{2,1}^\sigma$, there exists a nonnegative summable sequence $(c_j(f))_{j\in \Z}$ with $\|c_j(f)\|_{l^1(\Z)}=1$, such that
\begin{multline}
\frac{1}{\ee^2} \int_{\R^d} e^{-\frac{|y|^2}{16}}\|f_j(.-\ee y) +f_j(.+\ee y) -2f_j(.)\|_{L^2} dy\\
 \leq C_{\sigma,d} 2^{-j\sigma} c_j(f) \|\frac{\phi_\ee* f-f}{\ee^2}\|_{\dot{B}_{2,1}^\sigma}. 
\end{multline}
}
\label{lemIIIb}
\end{lem}
\textbf{Proof of lemma \ref{lemIIIa}:} we only need to estimate with lemma \ref{lemexp} the following decomposition of the coefficient appearing in $III_j$:
$$
1-e^A e^B = 1-e^B +e^B(1-e^A).
$$
The first estimate is straightforward: as $|B|\leq C\int_0^t \|\nabla \dot{S}_{j-1} v (\tau)\|_{L^\infty} d\tau \leq CV(t)$,
\begin{equation}
 |1-e^B|\leq |B|e^{|B|} \leq CV(t) e^{CV(t)}.
\label{bout1}
\end{equation}
Before estimating the second term, as in the work of T. Hmidi and S. Keraani (\cite{TH2}) we need to give a precise bound for the variation ratio of the flow:
\begin{lem}
 \sl{Under the same assumptions, for all $x,y\in \R^d$ with $y\neq 0$, we have:
$$
\left|\frac{|\psi_{j,t}^{-1}(x)-\psi_{j,t}^{-1}(x-\ee y)|^2}{\ee^2|y|^2} -1\right|\leq e^{2CV(t)}-1.
$$
}
\label{estimtauxflot}
\end{lem}
\textbf{Proof:} For all $x\neq y \in \R^d$, thanks to \eqref{estimflot},
$$
\frac{|\psi_{j,t}^{-1}(x)-\psi_{j,t}^{-1}(y)|}{|x-y|}\leq \|D \psi_{j,t}^{-1}\|_{L^\infty} \leq e^{CV(t)}.
$$
Similarly, for all $x\neq y \in \R^d$,
$$
\frac{|\psi_{j,t}(x)-\psi_{j,t}(y)|}{|x-y|}\leq \|D \psi_{j,t}\|_{L^\infty} \leq e^{CV(t)}.
$$
and applying it at $\psi_{j,t}^{-1}(x)$ and $\psi_{j,t}^{-1}(y)$,
$$
\frac{|x-y|}{|\psi_{j,t}^{-1}(x)-\psi_{j,t}^{-1}(y)|}\leq e^{CV(t)},
$$
so that we get that for all $x$ and $y\neq 0$,
$$
-e^{-2CV(t)}(e^{2CV(t)}-1)=e^{-2CV(t)}-1 \leq \frac{|\psi_{j,t}^{-1}(x)-\psi_{j,t}^{-1}(x-\ee y)|^2}{\ee^2|y|^2} -1\leq e^{2CV(t)}-1
$$
in other words,
$$
-1\leq -e^{-2CV(t)}\leq \frac{1}{e^{2CV(t)}-1} \left(\frac{|\psi_{j,t}^{-1}(x)-\psi_{j,t}^{-1}(x-\ee y)|^2}{\ee^2|y|^2} -1\right)\leq 1
$$
that is exactly the announced result. $\blacksquare$
\\
Coming back to lemma \ref{lemIIIa} and thanks to lemma \ref{lemexp}, we immediately get that:
\begin{equation}
 e^B|1-e^A|\leq e^{|B|}e^{|A|}|A| \leq e^{CV(t)} e^{\displaystyle{\frac{|y|^2}{4}(e^{2CV(t)}-1)}} \frac{|y|^2}{4}(e^{2CV(t)}-1).
\label{bout2}
\end{equation}
Gathering \eqref{bout1} and \eqref{bout2}, we get
$$
\|1-e^A e^B\|_{L_x^\infty} \leq e^{\displaystyle{CV(t)}} \left(CV(t)+ e^{\displaystyle{\frac{|y|^2}{4}(e^{2CV(t)}-1)}} \frac{|y|^2}{4}(e^{2CV(t)}-1)\right).
$$
so that, taking the $L_x^2$-norm in the expression of $III_j$ from \eqref{IIIetIV}, we obtain:
\begin{multline}
 \|III_j\|_{L_x^2} \leq \displaystyle{\frac{1}{2(2\pi)^d\ee^2} \int_{\R^d} e^{\displaystyle{-\frac{|y|^2}{4}}} \|f_j(.-\ee y) +f_j(.+\ee y) -2f_j(.)\|_{L^2} \|1-e^A e^B\|_{L^\infty}dy,}\\
\leq \frac{1}{2(2\pi)^d\ee^2} \int_{\R^d} e^{CV(t)} \left(CV(t)e^{-\frac{|y|^2}{4}} + \frac{|y|^2}{4} e^{-\frac{|y|^2}{4}} e^{\frac{|y|^2}{4}(e^{2CV(t)}-1)} (e^{2CV(t)}-1)\right)\\
\times \|f_j(.-\ee y) +f_j(.+\ee y) -2f_j(.)\|_{L_x^2} dy.
\end{multline}
Assume that $t$ is so small that $e^{2CV(t)}-1\leq \frac{1}{2}$, then (using again that for all $x\geq 0$, $xe^{-x}\leq \frac{2}{e}e^{-x/2}$)
$$
\frac{|y|^2}{4} e^{-\frac{|y|^2}{4}} e^{\frac{|y|^2}{4}(e^{2CV(t)}-1)}\leq \frac{|y|^2}{4} e^{-\frac{|y|^2}{8}} \leq \frac{4}{e} e^{-\frac{|y|^2}{16}},
$$
which ends the proof of lemma \ref{lemIIIa}. $\blacksquare$
\begin{rem}
 \sl{Similarly, writing $L_\ee(f_j)=\frac{1}{2}(L_\ee(f_j)+L_\ee(f_j))$ and performing the change of variable $z=-y$ in the second term lead to:
$$
L_\ee(f_j)= \displaystyle{\frac{1}{2(2\pi)^d\ee^2} \int_{\R^d} e^{\displaystyle{-\frac{|y|^2}{4}}}\left(f_j(x-\ee y) +f_j(x+\ee y) -2f_j(x)\right) dy,}
$$ 
taking the $L_x^2$-norm we obtain that:
$$
\|\frac{\phi_\ee* f_j-f_j}{\ee^2}\|_{L^2} \leq \displaystyle{\frac{1}{2(2\pi)^d\ee^2} \int_{\R^d} e^{\displaystyle{-\frac{|y|^2}{4}}}\|f_j(.-\ee y) +f_j(.+\ee y) -2f_j(.)\|_{L^2} dy,}
$$
and as $\tau_{\alpha}(\ddj f)=\ddj(\tau_{\alpha}f)$ (which is localized in frequency), we immediately get that:
$$
\|\frac{\phi_\ee* f-f}{\ee^2}\|_{\dot{B}_{2,1}^\sigma} \leq \displaystyle{\frac{C}{\ee^2} \int_{\R^d} e^{\displaystyle{-\frac{|y|^2}{16}}}\|f(.-\ee y) +f(.+\ee y) -2f(.)\|_{\dot{B}_{2,1}^\sigma} dy}.
$$
On the other hand, lemma \ref{lemIIIb} gives the reverse estimate so that:
$$
\|q\|_{\dot{B}_{\ee}^{\sigma+2,\sigma}} \sim \displaystyle{\frac{1}{\ee^2} \int_{\R^d} e^{\displaystyle{-\frac{|y|^2}{16}}}\|f(.-\ee y) +f(.+\ee y) -2f(.)\|_{\dot{B}_{2,1}^\sigma} dy.}
$$
}
\end{rem}
\textbf{Proof of lemma \ref{lemIIIb}:} let us adapt the proof of \cite{Dbook} (th $2.36$ and $2.37$). As $f_j= \ddj f$, we can write that for all $x$ and $y$,
$$
f_j(x-\ee y)+f_j(x+\ee y)-2f_j(x)= \left(\tau_{-\ee y} \ddj f +\tau_{\ee y} \ddj f -2\ddj f\right)(x)
$$
and as in the proof of theorem $2.37$ from \cite{Dbook}, thanks twice to the mean value theorem, we get that:
$$
\|\tau_{-\ee y} \ddj f +\tau_{\ee y} \ddj f -2\ddj f\|_{L^2}\leq C 2^{2j} \ee^2 |y|^2 \sum_{|j'-j|\leq 1} \|\dot{\Delta}_{j'} f\|_{L^2},
$$
and thanks to the definition of the hybrid Besov norm $\|\frac{\phi_\ee *f-f}{\ee^2}\|_{\dot{B}_{2,1}^\sigma}= \|f\|_{\dot{B}_{\ee}^{\sigma+2, \sigma}}$, there exists a nonnegative summable sequence (whose summation is $1$) $(c_j)_{j\in \Z}$ such that:
$$
\|\tau_{-\ee y} \ddj f +\tau_{\ee y} \ddj f -2\ddj f\|_{L^2}\leq C 2^{2j} \ee^2 |y|^2 \sum_{|j'-j|\leq 1} \frac{2^{-j'\sigma} c_{j'}}{\min(\frac{1}{\ee^2},2^{2j'})} \|\frac{\phi_\ee *f-f}{\ee^2}\|_{\dot{B}_{2,1}^\sigma}.
$$
Moreover
$$
\frac{1}{\min(\frac{1}{\ee^2},2^{2j'})}= \frac{1}{\min(\frac{1}{\ee^2},2^{2j})} \frac{\max(\ee^2, 2^{-2j'})}{\max(\ee^2, 2^{-2j})},
$$
and thanks to the following result,
\begin{lem}(\cite{Dinv}, Proposition $5.3$)
 \sl{Let $\alpha>0$, $a,b\in\R$. Then we have
$$
\frac{\max(\alpha,2^{-a})}{\max(\alpha,2^{-b})}\leq
\begin{cases}
 1 & \mbox{if}\quad a\geq b,\\
2^{b-a} & \mbox{if}\quad a\leq b.
\end{cases}
$$}
\label{estimax}
\end{lem}
as $|j'-j|\leq 1$, we obtain that:
$$
\frac{1}{\min(\frac{1}{\ee^2},2^{2j'})}\leq 2 \frac{1}{\min(\frac{1}{\ee^2},2^{2j})},
$$
so there exists a nonnegative summable sequence of summation $1$, that we will also denote by $(c_j)_{j\in \Z}$, such that:
$$
 \|\tau_{-\ee y} \ddj f +\tau_{\ee y} \ddj f -2\ddj f\|_{L^2}\leq C_\sigma 2^{2j} \ee^2 |y|^2 2^{-j\sigma} c_j \max(\ee^2, 2^{-2j}) \|\frac{\phi_\ee *f-f}{\ee^2}\|_{\dot{B}_{2,1}^\sigma}.
$$
On the other hand, we can also roughly write:
$$
\|\tau_{-\ee y} \ddj f +\tau_{\ee y} \ddj f -2\ddj f\|_{L^2}\leq 3\|\ddj f\|_{L^2} \leq C 2^{-j\sigma} c_j \max(\ee^2, 2^{-2j}) \|\frac{\phi_\ee *f-f}{\ee^2}\|_{\dot{B}_{2,1}^\sigma}.
$$
Combining the last two estimates, if we define $j_{\ee y}\in \Z$ the unique integer satisfying $1\leq \ee |y| 2^{j_{\ee y}}<2$ (that is $\ee |y| 2^{j_{\ee y}} \sim 1$) , we get the following statement:
\begin{multline}
 \|\tau_{-\ee y} \ddj f +\tau_{\ee y} \ddj f -2\ddj f\|_{L^2}\leq C_\sigma 2^{-j\sigma} c_j \max(\ee^2, 2^{-2j}) \|\frac{\phi_\ee *f-f}{\ee^2}\|_{\dot{B}_{2,1}^\sigma}\\
\times \begin{cases}
        2^{2j} \ee^2 |y|^2 & \mbox{ if } \ee|y|2^j<1 \quad (j<j_{\ee y})\\
1 & \mbox{ if } \ee|y|2^j\geq 1 \quad (j\geq j_{\ee y}).
       \end{cases}
\label{estimdiff2cas}
\end{multline}
We can now prove lemma \ref{lemIIIb}: let us follow the proof from \cite{Dbook} and begin by splitting the integral:
\begin{multline}
\frac{1}{\ee^2} \int_{\R^d} e^{-\frac{|y|^2}{16}}\|f_j(.-\ee y) +f_j(.+\ee y) -2f_j(.)\|_{L^2} dy\\
  \leq \frac{1}{\ee^2} \left(2^{2j} \ee^2 \int_{\ee |y|2^j<1} |y|^2 e^{-\frac{|y|^2}{16}}dy +\int_{\ee |y|2^j\geq 1} e^{-\frac{|y|^2}{16}}dy\right)\\
\times C_\sigma 2^{-j\sigma} c_j \max(\ee^2, 2^{-2j}) \|\frac{\phi_\ee *f-f}{\ee^2}\|_{\dot{B}_{2,1}^\sigma}
\label{intsplit}
\end{multline}
We have:
\begin{multline}
 \int_{|y|<\frac{1}{\ee 2^j}} |y|^2 e^{-\frac{|y|^2}{16}}dy =C \int_0^{\frac{1}{\ee 2^j}} r^2 e^{-\frac{r^2}{16}} r^{d-1}dr \leq C_d \int_0^{\frac{1}{\ee 2^j}} r e^{-\frac{r^2}{32}} dr = C_d (1-e^{-\frac{1}{32\ee^2 2^{2j}}})\\
\leq C_d \min(1, \frac{1}{\ee^2 2^{2j}}).
\end{multline}
Similarly,
\begin{multline}
 \int_{|y|\geq \frac{1}{\ee 2^j}} e^{-\frac{|y|^2}{16}}dy =C \int_{\frac{1}{\ee 2^j}}^\infty e^{-\frac{r^2}{16}} r^{d-1}dr \leq C_d \int_{\frac{1}{\ee 2^j}}^\infty r e^{-\frac{r^2}{32}} dr = C_d e^{-\frac{1}{32\ee^2 2^{2j}}}.
\end{multline}
Plugging this into \eqref{intsplit}, we obtain that:
\begin{multline}
\frac{1}{\ee^2} \int_{\R^d} e^{-\frac{|y|^2}{16}}\|f_j(.-\ee y) +f_j(.+\ee y) -2f_j(.)\|_{L^2} dy\\
  \leq C_\sigma 2^{-j\sigma} c_j \max(1, \frac{1}{\ee^2 2^{2j}}) \|\frac{\phi_\ee *f-f}{\ee^2}\|_{\dot{B}_{2,1}^\sigma} C_d \left(2^{2j} \ee^2 \min(1, \frac{1}{\ee^2 2^{2j}}) +e^{-\frac{1}{32\ee^2 2^{2j}}} \right)\\
\leq C_{\sigma,d} 2^{-j\sigma} c_j \|\frac{\phi_\ee *f-f}{\ee^2}\|_{\dot{B}_{2,1}^\sigma} \left(R_1 +R_2\right),
\end{multline}
where
$$
\begin{cases}
\displaystyle{R_1= 2^{2j} \ee^2 \min(1, \frac{1}{\ee^2 2^{2j}}) \max(1, \frac{1}{\ee^2 2^{2j}}) =1,}\\
\displaystyle{R_2= \max(e^{-\frac{1}{32\ee^2 2^{2j}}}, \frac{1}{\ee^2 2^{2j}}e^{-\frac{1}{32\ee^2 2^{2j}}}) \leq \max(1, C e^{-\frac{1}{64\ee^2 2^{2j}}}) \leq C,}
\end{cases}
$$
which ends the proof of lemma \ref{lemIIIb}. $\blacksquare$

Lemmas \ref{lemIIIa} and \ref{lemIIIb} immediately imply proposition \ref{propIII}, and we can now turn to $IV_j$:
\begin{prop}
 \sl{Under the same assumptions, there exist a positive constant $C=C_{\sigma,d}$ and a nonnegative summable sequence $(c_j= c_j(f))_{j\in \Z}$ whose sum is $1$, such that if $t$ is so small that $e^{2CV(t)}-1\leq \frac{1}{2}$, we have:
$$
\|IV_j\|_{L^2}\leq C_{\sigma,d} (V+e^{CV}-1) e^{CV} 2^{-j\sigma} c_j \|\frac{\phi_\ee* f-f}{\ee^2}\|_{\dot{B}_{2,1}^\sigma}.
$$
\label{propIV}
}
\end{prop}
To prove this result we will prove the following lemmas:
\begin{lem}
 \sl{There exists a constant $C$ such that for all $j\in\Z$, $f$, and all $t$ so small that $e^{2CV(t)}-1\leq \frac{1}{2}$,
\begin{multline}
 \|IV_j\|_{L^2}\leq \frac{1}{(2\pi)^d} (V +e^{CV}-1) e^{CV} \\
\times \frac{1}{\ee^2} \int_{\R^d} \min(1, \ee 2^j |y|) e^{-\frac{|y|^2}{16}}\|f_j(.+\ee y) -f_j(.)\|_{L^2} dy,
\end{multline}
where $V$ is defined in \eqref{defV}.
}
\label{lemIVa}
\end{lem}
\begin{lem}
 \sl{For all $\sigma\in \R$ there exists a constant $C_{\sigma,d}$ such that for any $f\in \dot{B}_{2,1}^\sigma$, there exists a nonnegative summable sequence $(c_j(f))_{j\in \Z}$ with $\|c_j(f)\|_{l^1(\Z)}=1$, such that
\begin{multline}
\frac{1}{\ee^2} \int_{\R^d} \min(1, \ee 2^j |y|) e^{-\frac{|y|^2}{16}}\|f_j(.+\ee y) -f_j(.)\|_{L^2} dy\\
 \leq C_{\sigma,d} 2^{-j\sigma} c_j(f) \|\frac{\phi_\ee* f-f}{\ee^2}\|_{\dot{B}_{2,1}^\sigma}. 
\end{multline}
}
\label{lemIVb}
\end{lem}
\textbf{Proof of lemma \ref{lemIVa}}: as announced here we can only rely on a finite difference of order 1, and we need to carefully estimate the following coefficient to obtain summable sequences in low frequencies:
$$
e^A e^B - e^C e^D = e^B (e^A -e^C) + e^C (e^B-e^D).
$$
Let us begin with the second term: thanks to lemma \ref{lemexp} (we refer to \eqref{ABCD} for the expressions of the various terms involved),
\begin{multline}
 |e^B-e^D|\leq e^{\max(|D|, |B|)}|D-B| \leq e^{\max(|D|, |B|)} (|D|+|B|)\\
\leq 2e^{CV} \int_0^t \|\div \dot{S}_{j-1}(\tau)\|_{L^\infty} d\tau \leq 2V e^{CV}.
\end{multline}
and then,
\begin{equation}
 e^C |e^B-e^D| \leq 2e^{\displaystyle{\frac{|y|^2}{4}(e^{2CV(t)}-1)}} V e^{CV}.
\label{DB1}
\end{equation}
This estimate is useful only for high frequencies: indeed for the low frequency regime, as in the proof of proposition \ref{propIII}, after integration with respect to $y$, the result won't be summable when $j$ goes to $-\infty$. In order to get a suitable estimate for low frequencies, let us use once again the elementary mean value theorem:
\begin{multline}
|e^B-e^D|\leq e^{\max(|D|, |B|)}|D-B|\\
\leq e^{CV(t)} \left| \int_{0}^t \left( \div \dot{S}_{j-1}v(\tau, X_j(\tau,t,x+\ee y)) -\div \dot{S}_{j-1}v(\tau, X_j(\tau,t,x-\ee y))\right)d\tau \right|\\
\leq e^{CV(t)} \int_{0}^t 2\ee |y| \|\nabla \left(\div \dot{S}_{j-1}v(\tau, X_j(\tau,t,.))\right)\|_{L^\infty} d\tau.
\end{multline}
Obviously:
$$
\nabla \left(\div \dot{S}_{j-1}v(\tau, X_j(\tau,t,z))\right)= \nabla \div \dot{S}_{j-1}v(\tau, X_j(\tau,t,z)) \times DX_j(\tau,t,z)),
$$
so that, thanks to the Bernstein lemma (we refer to the appendix), \eqref{estimflot} and \eqref{detjacobienj}:
\begin{multline}
 \|\nabla \left(\div \dot{S}_{j-1}v(\tau, X_j(\tau,t,.))\right)\|_{L^\infty} \leq \|\nabla \div \dot{S}_{j-1}v(\tau, X_j(\tau,t,.))\|_{L^\infty} \|DX_j(\tau,t,.))\|_{L^\infty}\\
\leq \|\nabla \div \dot{S}_{j-1}v(\tau)\|_{L^\infty} e^{CV(t)} \leq 2^j \|\nabla v(\tau)\|_{L^\infty} e^{CV(t)},
\end{multline}
and
\begin{equation}
|e^B-e^D|\leq e^{2CV} 2\ee 2^j |y| \int_0^t \|\nabla v(\tau)\|_{L^\infty} d\tau \leq e^{2CV} 2\ee 2^j |y| V(t).
\label{DB2}
\end{equation}
Combining \eqref{DB1} and \eqref{DB2} we end up with the better estimate:
\begin{equation}
e^C |e^B-e^D|\leq Ce^{\displaystyle{\frac{|y|^2}{4}(e^{2CV(t)}-1)}} e^{CV} V(t) \min(1, \ee 2^j |y|),
\label{DB}
\end{equation}
which now can be used for all frequency. The same has to be done for the other term: rough estimates first imply that
\begin{equation}
 |e^A-e^C|\leq e^{\max(|A|, |C|)}2\max(|A|, |C|) \leq 2 e^{\displaystyle{\frac{|y|^2}{4}(e^{2CV(t)}-1)}} \frac{|y|^2}{4}(e^{2CV(t)}-1).
\label{AC1}
\end{equation}
And more precisely (needed later for the low frequencies convergence), we write that:
\begin{equation}
 |e^A-e^C|\leq e^{\displaystyle{\frac{|y|^2}{4}(e^{2CV(t)}-1)}} |F|,
\end{equation}
with
\begin{multline}
F= \frac{|y|^2}{4} \left(\frac{|\psi_{j,t}^{-1}(x)-\psi_{j,t}^{-1}(x+\ee y)|^2}{\ee^2|y|^2} -\frac{|\psi_{j,t}^{-1}(x)-\psi_{j,t}^{-1}(x-\ee y)|^2}{\ee^2|y|^2}\right)\\
=\frac{1}{4\ee^2} \left(|\psi_{j,t}^{-1}(x)-\psi_{j,t}^{-1}(x+\ee y)|^2 -|\psi_{j,t}^{-1}(x)-\psi_{j,t}^{-1}(x-\ee y)|^2\right).
\end{multline}
Thanks to the identity $|a|^2-|b|^2=(a+b|a-b)$, where $(.|.)$ is the usual scalar product in $\R^d$, we get:
\begin{multline}
|F| =\frac{1}{4\ee^2} \left|\left(2\psi_{j,t}^{-1}(x)-\psi_{j,t}^{-1}(x+\ee y)-\psi_{j,t}^{-1}(x-\ee y)\Big| \psi_{j,t}^{-1}(x-\ee y) -\psi_{j,t}^{-1}(x+\ee y)\right)\right|\\
\leq \frac{1}{4\ee^2} |2\psi_{j,t}^{-1}(x)-\psi_{j,t}^{-1}(x+\ee y)-\psi_{j,t}^{-1}(x-\ee y)|\times |\psi_{j,t}^{-1}(x-\ee y) -\psi_{j,t}^{-1}(x+\ee y)|.
\end{multline}
Thanks to the mean value theorem (used twice for the first factor and once for the second), we can write:
$$
|F| \leq \frac{1}{4\ee^2} (\ee |y|)^2 \|D^2 \psi_{j,t}^{-1}\|_{L^\infty} \times(\ee |y|) \|D \psi_{j,t}^{-1}\|_{L^\infty} \leq \frac{1}{4} \ee |y|^3 \|D^2 \psi_{j,t}^{-1}\|_{L^\infty} \|D \psi_{j,t}^{-1}\|_{L^\infty},
$$
and using the estimates for the flow (see \eqref{estimflot}), we obtain that
$$
 |F| \leq C \ee |y|^3 2^j e^{CV} (e^{CV}-1),
$$
and
\begin{equation}
|e^A-e^C|\leq C e^{\displaystyle{\frac{|y|^2}{4}(e^{2CV(t)}-1)}} \ee |y|^3 2^j e^{CV} (e^{CV}-1).
\label{AC2}
\end{equation}
From \eqref{AC1} and \eqref{AC2} we deduce that:
\begin{equation}
e^B |e^A -e^C| \leq C  e^{CV(t)}(e^{CV(t)}-1) \min(1, \ee 2^j |y|) |y|^2 e^{\displaystyle{\frac{|y|^2}{4}(e^{2CV(t)}-1)}}.
\label{AC}
\end{equation}
Finally, thanks to \eqref{AC} and \eqref{DB},
\begin{equation}
 \|e^A e^B - e^C e^D\|_{L_x^\infty} \leq C e^{CV(t)}(V(t)+e^{CV(t)}-1) \min(1, \ee 2^j |y|) (1+|y|^2) e^{\displaystyle{\frac{|y|^2}{4}(e^{2CV(t)}-1)}}.
\end{equation}
Taking the $L^2$-norm in \eqref{IIIetIV}, and thanks to the previous estimate, we have:
\begin{multline}
 \|IV_j\|_{L_x^2} \leq \displaystyle{\frac{1}{2(2\pi)^d\ee^2}} \int_{\R^d} e^{\displaystyle{-\frac{|y|^2}{4}}}\|f_j(.+\ee y)-f_j(.)\|_{L^2} \|e^A e^B-e^C e^D\|_{L^\infty}dy\\
\leq \displaystyle{\frac{1}{2(2\pi)^d}} C e^{CV(t)}(V+e^{CV(t)}-1)\\
\times\frac{1}{\ee^2}\int_{\R^d} e^{\displaystyle{-\frac{|y|^2}{4}}}  \min(1, \ee 2^j |y|) (1+|y|^2) e^{\displaystyle{\frac{|y|^2}{4}(e^{2CV(t)}-1)}} \|f_j(.+\ee y)-f_j(.)\|_{L^2} dy
\end{multline}
and as previously, if $t$ is small enough so that $e^{2CV(t)}-1\leq\frac{1}{2}$,
\begin{multline}
 \|IV_j\|_{L_x^2} \leq \displaystyle{\frac{1}{2(2\pi)^d}} C e^{CV(t)}(V+e^{CV(t)}-1)\\
\times\frac{1}{\ee^2}\int_{\R^d} (1+|y|^2)e^{\displaystyle{-\frac{|y|^2}{8}}}  \min(1, \ee 2^j |y|)  \|f_j(.+\ee y)-f_j(.)\|_{L^2} dy
\end{multline}
We conclude thanks to the fact that for all $x\geq 0$, $x e^{-x} \leq \frac{2}{e} e^{-\frac{x}{2}}$. $\blacksquare$
\\
The proof of lemma \ref{lemIVb} will follow the lines of the proof of theorem 2.36 from \cite{Dbook}:   
$$
f_j(x+\ee y)-f_j(x)= \left(\tau_{\ee y} \ddj f -\ddj f\right)(x)
$$
and thanks to the mean value theorem,
$$
\|\tau_{\ee y} \ddj f -\ddj f\|_{L^2}\leq C 2^j \ee |y| \sum_{|j'-j|\leq 1} \|\dot{\Delta}_{j'} f\|_{L^2},
$$
and as in the proof of lemma \ref{lemIIIb}, using the definition of the hybrid Besov norm, there exists a nonnegative summable sequence (whose summation is $1$) once more denoted by $(c_j)_{j\in \Z}$ (depending on $f$ and $t$) such that:
$$
\|\tau_{\ee y} \ddj f -\ddj f\|_{L^2}\leq C 2^j \ee |y| 2^{-j\sigma} c_j \max(\ee^2,2^{-2j}) \|\frac{\phi_\ee *f-f}{\ee^2}\|_{\dot{B}_{2,1}^\sigma}.
$$
On the other hand, a rough estimate provides:
$$
\|\tau_{\ee y} \ddj f -\ddj f\|_{L^2}\leq 2 \|\ddj f\|_{L^2} \leq C 2^{-j\sigma} c_j \max(\ee^2,2^{-2j}) \|\frac{\phi_\ee *f-f}{\ee^2}\|_{\dot{B}_{2,1}^\sigma}.
$$
so that we end up with:
\begin{multline}
 \|\tau_{\ee y} \ddj f -\ddj f\|_{L^2}\leq C_\sigma 2^{-j\sigma} c_j \max(\ee^2, 2^{-2j}) \|\frac{\phi_\ee *f-f}{\ee^2}\|_{\dot{B}_{2,1}^\sigma}\\
\times \begin{cases}
        2^j \ee |y| & \mbox{ if } \ee|y|2^j<1 \quad (j<j_{\ee y})\\
1 & \mbox{ if } \ee|y|2^j\geq 1 \quad (j\geq j_{\ee y}).
       \end{cases}
\label{estimdiff2cas2}
\end{multline}
Returning to the integral, we have:
\begin{multline}
\frac{1}{\ee^2} \int_{\R^d} \min(1, \ee 2^j |y|) e^{-\frac{|y|^2}{16}}\|f_j(.+\ee y) -f_j(.)\|_{L^2} dy\\
 \leq \max(1,\frac{1}{\ee^2 2^{2j}}) C_{\sigma,d} 2^{-j\sigma} c_j(f) \|\frac{\phi_\ee* f-f}{\ee^2}\|_{\dot{B}_{2,1}^\sigma} \int_{\R^d} \min(1, \ee^2 2^{2j} |y|^2) e^{-\frac{|y|^2}{16}} dy. 
\end{multline}
Then the same computation as in the proof of lemma \ref{lemIIIb} leads to:
\begin{multline}
\frac{1}{\ee^2} \int_{\R^d} \min(1, \ee 2^j |y|) e^{-\frac{|y|^2}{16}}\|f_j(.+\ee y) -f_j(.)\|_{L^2} dy\\
 \leq C_{\sigma,d} 2^{-j\sigma} c_j(f) \|\frac{\phi_\ee* f-f}{\ee^2}\|_{\dot{B}_{2,1}^\sigma}. 
\end{multline}
This concludes the proof of lemma \ref{lemIVb} and together with lemma \ref{lemIVa} it implies proposition \ref{propIV}. Finally, propositions \ref{propIII} and \ref{propIV}, with the first estimate from \eqref{estimflot} end the proof of theorem \ref{II'}. $\blacksquare$

\subsection{End of the proof}

Once obtained the desired estimates of the additionnal remainder term introduced by the effect of the Lagrangian change of coordinates on the non-local capillary term, we can go back to the use of the linear estimates from proposition \ref{estimlinloc} on system \eqref{systchanged}: for all $l\in \Z$, as $\tilde q_j (0)= q_j (0, \psi_{j,0}(.))=q_j(0)$,
\begin{multline}
 \|\ddl \tilde u_j\|_{L_t^\infty L^2}+\nu_0 2^{2l} \|\ddl \tilde u_j\|_{L_t^1 L^2}+(1+\nu 2^l)\left(\|\ddl \tilde q_j\|_{L_t^\infty L^2} +\nu\min(\frac{1}{\ee^2}, 2^{2l})\|\ddj \tilde q_j\|_{L_t^1 L^2}\right)\\
\leq C_{p, \frac{\nu^2}{4\kappa}} \Bigg[ (1+\nu 2^l)\|\ddl q_j (0)\|_{L^2} +\|\ddl u_j(0)\|_{L^2} + (1+\nu 2^l)\|\ddl \tilde f_j+ \ddl R_j^1\|_{L_t^1 L^2}\\
+\|\ddl \tilde g_j +\ddl R_j^2 +\ddl R_j^3 +\kappa \ddl R_j\|_{L_t^1 L^2}\Bigg].
\label{estimlinloc2}
\end{multline}
Let us recall that thanks to \eqref{energieBF2}, all we need is to estimate the high frequencies, that is $(q_j, u_j)$ for $j\geq 0$. For this, as in \cite{Dlagrangien} and \cite{CD}, let us define some $N_0\in \Z$ (that will be fixed later), and write:
$$
\|q_j\|_{L^2} =\|\tilde q_j \circ \psi_{j,t}^{-1}\|_{L^2} \leq e^{CV} \|\tilde q_j\|_{L^2} \leq e^{CV} \left( \|\dot{S}_{j-N_0}\tilde q_j\|_{L^2} +\sum_{l\geq j-N_0} \|\ddl \tilde q_j\|_{L^2}\right).
$$
Moreover, as in \cite{Dlagrangien} and \cite{CD}, we can use the version of \eqref{Vishincompr} given by lemma \ref{Vishcompr}  in the non measure-preserving case (we refer for example to \cite{Dbook} lemma $2.6$, or \cite{Dlagrangien}, lemma $A.1$) to obtain:
$$
\|\dot{S}_{j-N_0}\tilde q_j\|_{L^2}\leq C e^{CV} \left(e^{CV}-1+ 2^{-N_0} e^{CV}\right) \|q_j\|_{L^2},
$$
so that:
\begin{equation}
\|q_j\|_{L^2} \leq C e^{CV}\left((e^{CV}-1+ 2^{-N_0} e^{CV}) \|q_j\|_{L^2} +\sum_{l\geq j-N_0} \|\ddl \tilde q_j\|_{L^2}\right),
\label{estqlagrang}
\end{equation}
and for the velocity,
\begin{equation}
\|u_j\|_{L^2} \leq C e^{CV}\left((e^{CV}-1+ 2^{-N_0} e^{CV}) \|u_j\|_{L^2} +\sum_{l\geq j-N_0} \|\ddl \tilde u_j\|_{L^2}\right).
\end{equation}
Multiplying \eqref{estqlagrang} by $2^j$ and using that in the summation, $l\geq j-N_0$, we obtain:
$$
2^j\|q_j\|_{L^2}\leq C e^{CV}\left((e^{CV}-1+ 2^{-N_0} e^{CV}) 2^j\|q_j\|_{L^2} +\sum_{l\geq j-N_0} 2^l 2^{N_0} \|\ddl \tilde q_j\|_{L^2}\right).
$$
Then, going back to $U_j$ (we refer to \eqref{Uj} for the definition), we can write that for all $j\geq 0$,
\begin{multline}
 U_j(t) \leq C e^{CV}\Bigg[\left(e^{CV}-1+ 2^{-N_0} e^{CV}\right) U_j(t) + \max(1,2^{2N_0}) \sum_{l\geq j-N_0}\\
\left( \|\ddl \tilde u_j\|_{L_t^\infty L^2}+\nu_0 2^{2l} \|\ddl \tilde u_j\|_{L_t^1 L^2}+(1+\nu 2^l)\Big(\|\ddl \tilde q_j\|_{L_t^\infty L^2} +\nu\min(\frac{1}{\ee^2}, 2^{2l})\|\ddj \tilde q_j\|_{L_t^1 L^2}\Big) \right)\Bigg],
\end{multline}
where we used once again lemma \ref{estimax} to prove that for all $l\geq j-N_0$:
\begin{equation}
 \min(\frac{1}{\ee^2},2^{2j})= \min(\frac{1}{\ee^2},2^{2l}) \frac{\max(\ee^2, 2^{-2l})}{\max(\ee^2, 2^{-2j})} \leq 2^{2N_0} \min(\frac{1}{\ee^2},2^{2l}).
\end{equation}
Thanks to \eqref{estimlinloc2} we obtain that for all $j\geq 0$,
\begin{multline}
 U_j(t) \leq C e^{CV}\Bigg[\left(e^{CV}-1+ 2^{-N_0} e^{CV}\right) U_j(t) + \max(1,2^{2N_0}) C_{p, \frac{\nu^2}{4\kappa}}\times\\
\sum_{l\geq j-N_0} \Big((1+\nu 2^l)\|\ddl q_j (0)\|_{L^2} +\|\ddl u_j(0)\|_{L^2} + (1+\nu 2^l)\Big(\|\ddl \tilde f_j\|_{L_t^1 L^2}+ \|\ddl R_j^1\|_{L_t^1 L^2}\Big)\\
+\|\ddl \tilde g_j\|_{L_t^1 L^2} +\|\ddl R_j^2\|_{L_t^1 L^2} +\|\ddl R_j^3\|_{L_t^1 L^2} +\kappa \|\ddl R_j\|_{L_t^1 L^2}\big)\Bigg],
\label{estimlinloc3}
\end{multline}
Except $R_j$, all of the remainder terms are the same as those from \cite{CD}. As we will sum for $l\geq j-N_0$, we need to provide estimates involving some positive power of $2^{j-l}$. To do this we use the well-known method of Vishik (see \cite{Vishik}) which takes advantage of the Bernstein lemma: the idea is to derivate the function in order to obtain a positive power of $2^{j-l}$ which is summable over $l\geq j-N_0$. Let us detail this for example on $\tilde f_j$: thanks to \eqref{estimflot},
$$
 \|\ddl \tilde f_j\|_{L^2} \leq C 2^{-l} \|\ddl \nabla \tilde f_j\|_{L^2} \leq C 2^{-l} \|\nabla f_j \circ \psi_{j,t}\|_{L^2} \|D\psi_{j,t}\|_{L^\infty}
$$
$$
\leq C 2^{-l} e^{2CV} \|\nabla f_j\|_{L^2} \leq C 2^{j-l} e^{2CV} \|f_j\|_{L^2}.
$$ 
and thanks to lemma \ref{lemmeB} (for $\sigma=s-1$), we can write that there exists a nonnegative sequence $(c_j(t))_{j\in \Z} \in l^1(\Z)$ such that $\|c\|_{l^1(\Z)}=1$ and a constant still denoted by $C=C(d,s)$ such that for all $j\in \Z$,
\begin{equation}
\|\ddl \tilde f_j\|_{L^2} \leq C 2^{j-l} e^{CV} \left(\|F_j\|_{L^2} + c_j(t) 2^{-j(s-1)}\|\nabla v\|_{\dot{B}_{2,1}^{\fd}} \|q\|_{\dot{B}_{2,1}^{s-1}}\right).
\label{ftilde}
\end{equation}
Similarly, we obtain that (using lemma \ref{lemmeB} for $\sigma=s$ in the second case):
\begin{equation}
 \begin{cases}
  \|\ddl \tilde g_j\|_{L^2} \leq C 2^{j-l} e^{CV} \left(\|G_j\|_{L^2} + c_j(t) 2^{-j(s-1)}\|\nabla v\|_{\dot{B}_{2,1}^{\fd}} \|u\|_{\dot{B}_{2,1}^{s-1}}\right)\\
2^l \|\ddl \tilde f_j\|_{L^2} \leq C\|\ddl \nabla \tilde f_j\|_{L^2} \leq C 2^{j-l} e^{CV} \left(2^j\|F_j\|_{L^2} + c_j(t) 2^{-j(s-1)}\|\nabla v\|_{\dot{B}_{2,1}^{\fd}} \|q\|_{\dot{B}_{2,1}^s}\right)
 \end{cases}
\label{gtilde}
\end{equation}
With the same argument we obtain (we refer to \cite{CD}, section $3.2$ for details):
\begin{equation}
 \begin{cases}
\|\ddl R_j^1\|_{L^2} \leq C 2^{j-l} (e^{CV}-1)e^{CV} 2^j \|u_j\|_{L^2} \leq C 2^{j-l} (e^{CV}-1)e^{CV} 2^{2j} \|u_j\|_{L^2}\\
2^l \|\ddl R_j^1\|_{L^2} \leq C\|\ddl \nabla R_j^1\|_{L^2} \leq C 2^{j-l} (e^{CV}-1)e^{CV} 2^{2j} \|u_j\|_{L^2}\\
\|\ddl R_j^2\|_{L^2} \leq C 2^{j-l} (e^{CV}-1)e^{CV} 2^j \|q_j\|_{L^2}\\
\|\ddl R_j^3\|_{L^2} \leq C (|\lambda+\mu|+\mu) 2^{j-l} (e^{CV}-1)e^{CV} 2^{2j} \|u_j\|_{L^2}
 \end{cases}
\label{Rj123}
\end{equation}
\begin{rem}
 \sl{Note that the last estimate from the first line is valid only when $j\geq 0$. Let us recall that we already treated the low frequencies thanks to the linear estimates (we refer to \eqref{energieBF2}).
}
\end{rem}
All that remains is to estimate $\kappa R_j$: as for the other terms, a direct estimate will not be sufficient as we eventually have to sum over $l\geq j-N_0$. We then need to use the Bernstein lemma in order to be able to factor by a positive power of $2^{j-l}$. Let us recall that $R_j= I_j + II_j'(\nabla q)$ with (we refer to \eqref{IjIIj}):
$$
\begin{cases}
\vspace{0.2cm}
I_j= \displaystyle{\frac{\phi_\ee*g_j-g_j}{\ee^2} \mbox{ where } g_j= \nabla q_j \circ \psi_{j,t} \times (I_d-D\psi_{j,t})},\\
II_j'(f)= \frac{\phi_\ee*f_j-f_j}{\ee^2}\circ \psi_{j,t} -\frac{\phi_\ee*(f_j \circ \psi_{j,t}) -f_j \circ \psi_{j,t}}{\ee^2}.\end{cases}
$$
Let us recall that section \ref{sectionRj} is devoted to the following estimate (see theorem \ref{II'} for details): there exists a nonnegative summable sequence $c_j(t)$ (whose sum is $1$) such that for all $j\in \Z$ (recall that $V(t)=\int_0^t \|\nabla v(\tau)\|_{L^\infty}d\tau$),
$$
\|II_j'(f)\|_{L^2}\leq C_{\sigma, d} e^{CV}(V+e^{2CV}-1) c_j(t) 2^{-j\sigma}\|\frac{\phi_\ee*f-f}{\ee^2}\|_{\dot{B}_{2,1}^\sigma},
$$
Let us begin with $I_j$: a rough estimate gives:
\begin{multline}
 \|\ddl I_j\|_{L^2} =\|\frac{\phi_\ee*\ddl g_j-\ddl g_j}{\ee^2}\|_{L^2} \leq C \min(\frac{1}{\ee^2}, 2^{2l}) \|g_j\|_{L^2}\\
\leq C \min(\frac{1}{\ee^2}, 2^{2l}) \|\nabla q_j \circ \psi_{j,t}\|_{L^2} \|(I_d-D\psi_{j,t})\|_{L^\infty}.
\end{multline}
Thanks again to lemma \ref{estimax},
$$
\min(\frac{1}{\ee^2},2^{2l})= \min(\frac{1}{\ee^2},2^{2j}) \frac{\max(\ee^2, 2^{-2j})}{\max(\ee^2, 2^{-2l})} \leq \min(\frac{1}{\ee^2},2^{2j})
\begin{cases}
 1 & \mbox{ if } j-N_0\leq l\leq j,\\
2^{2(l-j)} & \mbox{ if } j\leq l,
\end{cases}
$$
so that (thanks to \eqref{estimflot}) as $\max(1,2^{2(l-j)})\leq 2^{2N_0} 2^{2(l-j)}$ ($N_0$ will be taken positive),
$$
 \|\ddl I_j\|_{L^2} \leq C e^{CV} (e^{CV}-1) 2^{2N_0}2^{2(l-j)} \min(\frac{1}{\ee^2}, 2^{2j}) \|\nabla q_j\|_{L^2},
$$
which implies that in our case, in order to use Vishik's trick we need to derivate three times $\ddl I_j$ because we have to absorb $2^{2(l-j)}$ and obtain a positive power of $2^{j-l}$ (which is summable over $l\geq j-N_0$):
$$
\|\ddl I_j\|_{L^2}\leq C 2^{-3l} \|\ddl \nabla^3 I_j\|_{L^2} \leq C 2^{-3l} \|\frac{\phi_\ee*\ddl \nabla^3 g_j-\ddl \nabla^3 g_j}{\ee^2}\|_{L^2}.
$$
We then compute the third derivatives of $g_j= \nabla q_j \circ \psi_{j,t} \times (I_d-D\psi_{j,t})$ and obtain that:
$$
\|\nabla^3 g_j\|_{L^2} \leq C e^{4CV}(e^{CV}-1) 2^{3j} \|\nabla q_j\|_{L^2},
$$
which allows to conclude that
\begin{equation}
\|\ddl I_j\|_{L^2}\leq C e^{4CV}(e^{CV}-1) 2^{j-l} 2^{2N_0} \min(\frac{1}{\ee^2},2^{2j}) \|\nabla q_j\|_{L^2}.
\label{estimIj}
\end{equation}
For the second term ($II_j'(\nabla q)$), we have to derivate one time and decompose $\nabla II_j$ into three parts:
\begin{multline}
 \nabla II_j= \frac{\phi_\ee*\nabla^2 q_j-\nabla^2 q_j}{\ee^2}\circ \psi_{j,t} \cdot D\psi_{j,t} -\frac{\phi_\ee*(\nabla^2 q_j \circ \psi_{j,t} \cdot D\psi_{j,t}) -(\nabla^2 q_j \circ \psi_{j,t} \cdot D\psi_{j,t})}{\ee^2}\\
= A_j +B_j +C_j,
\end{multline}
where
$$
 \begin{cases}
A_j= \displaystyle{\frac{\phi_\ee*\nabla^2 q_j-\nabla^2 q_j}{\ee^2}\circ \psi_{j,t} \cdot (D\psi_{j,t}-I_d),}\\
B_j= II_j'(\nabla^2 q),\\
C_j= -\displaystyle{\frac{\phi_\ee*(\nabla^2 q_j \circ \psi_{j,t} \cdot (D\psi_{j,t}-I_d)) -\nabla^2 q_j \circ \psi_{j,t} \cdot (D\psi_{j,t}-I_d)}{\ee^2}},
 \end{cases}
$$
and thanks to the Bernstein lemma:
\begin{equation}
\|\ddl II_j\|_{L^2}\leq C 2^{-l} \|\ddl \nabla II_j\|_{L^2} \leq C 2^{-l} \left(\|\ddl A_j\|_{L^2} +\|\ddl B_j\|_{L^2} +\|\ddl C_j\|_{L^2}\right).
\label{decompIIj}
\end{equation}
Thanks to theorem \ref{II'} with $\sigma= s-2$ and for $t$ small enough ($e^{2CV}-1\leq \frac{1}{2}$, see the proof of theorem \ref{II'}) there exists a nonnegative summable sequence $(c_j(\tau))_{j\in\Z}$ whose sum is $1$ such that
\begin{multline}
 \|\ddl B_j\|_{L^2} \leq \|B_j\|_{L^2} \leq C_{s, d} e^{CV}(V+e^{2CV}-1) c_j(\nabla^2 q, \tau) 2^{-j(s-2)}\|\frac{\phi_\ee*\nabla^2 q-\nabla^2 q}{\ee^2}\|_{\dot{B}_{2,1}^{s-2}},\\
\leq C_{s, d} e^{CV}(V+e^{2CV}-1) c_j(\nabla^2 q, \tau) 2^j 2^{-j(s-1)}\|\frac{\phi_\ee*\nabla q-\nabla q}{\ee^2}\|_{\dot{B}_{2,1}^{s-1}}.
\label{estimBj}
\end{multline}
A direct estimate gives:
\begin{multline}
 \|\ddl A_j\|_{L^2} \leq \|A_j\|_{L^2} \leq \|\frac{\phi_\ee*\nabla^2 q_j-\nabla^2 q_j}{\ee^2} \circ \psi_{j,t}\|_{L^2}\cdot \|D\psi_{j,t}-I_d\|_{L^\infty}\\
\leq C e^{CV}(e^{CV}-1) 2^j \|\frac{\phi_\ee*\nabla q_j-\nabla q_j}{\ee^2}\|_{L^2},
\label{estimAj}
\end{multline}
and finally, for the same reason as in $I_j$, rough estimates of $\ddl C_j$ will provide $\min(\frac{1}{\ee^2}, 2^{2l})\leq 2^{2(l-j)} \min(\frac{1}{\ee^2}, 2^{2j})$ and we once again have to derivate $C_j$: this time we will derivate only two more times because we already have $2^{-l}$ in factor:
$$
2^{-l}\|\ddl C_j\|_{L^2} \leq C 2^{-3l} \|\ddl \nabla^2 C_j\|_{L^2} \leq C e^{3CV} (e^{CV}-1) 2^{2N_0} 2^{j-l} \min(\frac{1}{\ee^2}, 2^{2j}) \|\nabla q_j\|_{L^2}.
$$
Plugging this together with \eqref{estimAj} and \eqref{estimBj} into \eqref{decompIIj}, and using \eqref{estimIj} allows to get the estimate on $R_j$:
\begin{equation}
\|\ddl R_j\| \leq C e^{CV}(V+e^{2CV}-1) 2^{j-l} 2^{2N_0} c_j(t) 2^{-j(s-1)}\|\frac{\phi_\ee*\nabla q-\nabla q}{\ee^2}\|_{\dot{B}_{2,1}^{s-1}}
\label{estimRj}
\end{equation}
Plugging \eqref{ftilde}, \eqref{gtilde}, \eqref{Rj123} and \eqref{estimRj} into \eqref{estimlinloc3}, then summing over $l\geq j-N_0$ implies that for all $j\geq 0$ and for $t$ small enough ($e^{2CV}-1\leq \frac{1}{2}$) we have:
\begin{multline}
 U_j(t) \leq C e^{2CV}\Bigg[\left(e^{CV}-1+ 2^{-N_0}\right) U_j(t) + 2^{5N_0} C_{p, \frac{\nu^2}{4\kappa}} \Bigg((1+\nu 2^j)\|q_j (0)\|_{L^2} +\|u_j(0)\|_{L^2}\\
+(1+\nu 2^j)\|F_j\|_{L_t^1 L^2} +\|G_j\|_{L_t^1 L^2}+(e^{CV}-1) \left((1+|\lambda+\mu|+\mu+\nu) 2^{2j} \|u_j\|_{L_t^1 L^2} +\|\nabla q_j\|_{L_t^1 L^2}\right)\\
+ \int_0^t 2^{-j(s-1)} c_j(\tau) \|\nabla v(\tau)\|_{\dot{B}_{2,1}^\fd} \big(\|q(\tau)\|_{\dot{B}_{2,1}^{s-1}} +\nu\|q(\tau)\|_{\dot{B}_{2,1}^s} +\|u(\tau)\|_{\dot{B}_{2,1}^{s-1}}\big)d\tau\\
+\frac{\kappa}{\nu^2}(V+e^{2CV}-1) \int_0^t c_j(\tau) 2^{-j(s-1)}\nu^2\|\frac{\phi_\ee*\nabla q-\nabla q}{\ee^2}\|_{\dot{B}_{2,1}^{s-1}}d\tau \Bigg)\Bigg].
\end{multline}
As $j\geq 0$, we have $1\leq \min(\frac{1}{\ee}, 2^{2j})$ and
$$
(1+|\lambda+\mu|+\mu+\nu) 2^{2j} \|u_j\|_{L_t^1 L^2} +\|\nabla q_j\|_{L_t^1 L^2} \leq  \left(\frac{1+|\lambda+\mu|+\mu+\nu}{\nu_0} +\frac{1}{\nu^2}\right) U_j(t),
$$
so that (remember that $t$ is so small that $e^{CV}-1\leq \frac{1}{2}$):
\begin{multline}
 U_j(t) \leq \frac{9}{4}C \Bigg[\left((e^{CV}-1)\left(1+2^{5N_0} C_{p, \frac{\nu^2}{4\kappa}} ( \frac{1+|\lambda+\mu|+\mu+\nu}{\nu_0} +\frac{1}{\nu^2})\right)+ 2^{-N_0}\right) U_j(t) + \\
2^{5N_0} C_{p, \frac{\nu^2}{4\kappa}} \Bigg((1+\nu 2^j)\|q_j (0)\|_{L^2} +\|u_j(0)\|_{L^2} +(1+\nu 2^j)\|F_j\|_{L_t^1 L^2} +\|G_j\|_{L_t^1 L^2}\\
+ \int_0^t 2^{-j(s-1)} c_j(\tau) \|\nabla v(\tau)\|_{\dot{B}_{2,1}^\fd} \big(\|q(\tau)\|_{\dot{B}_{2,1}^{s-1}} +\nu\|q(\tau)\|_{\dot{B}_{2,1}^s} +\|u(\tau)\|_{\dot{B}_{2,1}^{s-1}}\big)d\tau\\
+\frac{\kappa}{\nu^2} (V+e^{2CV}-1) \int_0^t c_j(\tau) 2^{-j(s-1)} \nu^2\|\frac{\phi_\ee*\nabla q-\nabla q}{\ee^2}\|_{\dot{B}_{2,1}^{s-1}}d\tau \Bigg)\Bigg].
\end{multline}
Let us fix $N_0>0$ large enough, and take $t$ small enough so that:
\begin{equation}
 \begin{cases}
\vspace{0.2cm}
\displaystyle{\frac{9}{4} C\cdot 2^{-N_0} \leq \frac{1}{8},}\\
\vspace{0.2cm}
\displaystyle{\frac{9}{4}C(e^{CV}-1)\left(1+2^{5N_0} C_{p, \frac{\nu^2}{4\kappa}} ( \frac{1+|\lambda+\mu|+\mu+\nu}{\nu_0} +\frac{1}{\nu^2})\right) \leq \frac{1}{8}}
\end{cases}
\label{Cond2}
\end{equation}
Then we obtain that for all $j\geq 0$,
\begin{multline}
 U_j(t) \leq 3 C 2^{5N_0} C_{p, \frac{\nu^2}{4\kappa}} \Bigg((1+\nu 2^j)\|q_j (0)\|_{L^2} +\|u_j(0)\|_{L^2}\\
+(1+\nu 2^j)\|F_j\|_{L_t^1 L^2} +\|G_j\|_{L_t^1 L^2} +\int_0^t 2^{-j(s-1)} c_j(\tau) \|\nabla v(\tau)\|_{\dot{B}_{2,1}^\fd} U(\tau)d\tau\\
+\frac{\kappa}{\nu^2} (V+e^{2CV}-1) \int_0^t c_j(\tau) 2^{-j(s-1)} \nu^2\|\frac{\phi_\ee*\nabla q-\nabla q}{\ee^2}\|_{\dot{B}_{2,1}^{s-1}}d\tau \Bigg).
\end{multline}
Now, if $t$ is so small that:
\begin{equation}
 3 C 2^{5N_0} C_{p, \frac{\nu^2}{4\kappa}} (V+e^{2CV}-1) \leq \frac{1}{2} \frac{\nu^2}{\kappa},
\label{Cond3}
\end{equation}
then for all $j\geq 0$,
\begin{multline}
2^{j(s-1)} U_j(t) \leq 3 C 2^{5N_0} C_{p, \frac{\nu^2}{4\kappa}} \Bigg(2^{j(s-1)} U_j(0) +(1+\nu 2^j)2^{j(s-1)}\|F_j\|_{L_t^1 L^2} +2^{j(s-1)}\|G_j\|_{L_t^1 L^2}\\
+\int_0^t c_j(\tau) \|\nabla v(\tau)\|_{\dot{B}_{2,1}^\fd} U(\tau)d\tau \Bigg) + \frac{1}{2} \nu^2\int_0^t c_j(\tau)\|\frac{\phi_\ee*\nabla q-\nabla q}{\ee^2}\|_{\dot{B}_{2,1}^{s-1}}d\tau.
\end{multline}
Going back to the low frequencies, if we take $K=(2C_{p,\frac{\nu^2}{4\kappa}})^{-1}$ in \eqref{energieBF2}, then for all $j\leq 0$
\begin{multline}
2^{j(s-1)} U_j(t) \leq C_{p,\frac{\nu^2}{4\kappa}} \Bigg[2^{j(s-1)} U_j(0) + (1+\nu 2^j)2^{j(s-1)}\|\ddj F\|_{L_t^1 L^2}+ 2^{j(s-1)}\|\ddj G\|_{L_t^1 L^2}\\
+2^{-j(s-1)}\int_0^t c_j(\tau) \left( \big(\max(1,\frac{1}{\nu^3}) +\frac{1}{\nu_0}\big) \|v(\tau)\|_{\dot{B}_{2,1}^\fd}^2
+\|v(\tau)\|_{\dot{B}_{2,1}^{\fd+1}}\right) U(\tau) d\tau\Bigg]\\
+\frac{1}{2} 2^{-j(s-1)}\int_0^t c_j(\tau) \left(\nu_0 \|u\|_{\dot{B}_{2,1}^{s+1}}+ \nu \|q\|_{\dot{B}_{\ee}^{s+1,s-1}}+ \nu^2\|q\|_{\dot{B}_{\ee}^{s+2,s}} \right) d\tau.
\end{multline}
Summing over $j\in \Z$ gives that
\begin{multline}
 U(t) \leq \frac{U(t)}{2} +C_{p,\frac{\nu^2}{4\kappa}} \Bigg(U(0)+ \|F\|_{L_t^1 \dot{B}_{2,1}^{s-1}} +\nu\|F\|_{L_t^1 \dot{B}_{2,1}^s} +\|G\|_{L_t^1 \dot{B}_{2,1}^{s-1}}\\
+(\frac{1+|\lambda+\mu|+\mu+\nu}{\nu_0}+\max(1, \frac{1}{\nu^3})) \int_0^t W'(\tau) U(\tau)\Bigg)
\end{multline}
where
\begin{equation}
 V(t)\overset{\mbox{def}}{=}\int_0^t \|\nabla v(\tau)\|_{L^\infty}d\tau \leq W(t)\overset{\mbox{def}}{=}\int_0^t (\|\nabla v(\tau)\|_{\dot{B_{2,1}^\fd}} +\|v(\tau)\|_{\dot{B_{2,1}^\fd}}^2) d\tau. 
\label{defW}
\end{equation}
and thanks to the Gronwall lemma, we obtain that for $t$ small enough (satisfying conditions \eqref{Cond1}, \eqref{Cond2} and \eqref{Cond3}),
\begin{multline}
 U(t) \leq 2C_{p,\frac{\nu^2}{4\kappa}} \left(U(0)+ \|F\|_{L_t^1 \dot{B}_{2,1}^{s-1}} +\nu\|F\|_{L_t^1 \dot{B}_{2,1}^s} +\|G\|_{L_t^1 \dot{B}_{2,1}^{s-1}}\right)\\
\times e^{\displaystyle{2 C_{p,\frac{\nu^2}{4\kappa}} (\frac{1+|\lambda+\mu|+\mu+\nu}{\nu_0}+\max(1, \frac{1}{\nu^3}))W(t)}}
\end{multline}
Using conditions \eqref{Cond1}, \eqref{Cond2} and \eqref{Cond3} gives the result for small times $t$ (it does not depend on the initial data, only on $V(t)$ and such that \eqref{Cond2} is satisfied), then we globalize the result (for example as in \cite{TH3}, end of section 4) for all time $t$ by subdividing $[0,t]$ into intervals $[T_{i+1}, T_i]$ where we have $e^{2C\int_{T_i}^{T_{i+1}} V(\tau)d\tau}-1 \leq \frac{1}{2}$ as well as the analogous of \eqref{Cond2} (we also refer to \cite{Dbook} chapter 10 for a connexity argument). $\blacksquare$

\subsection{Extension of the results}

In this short section, we mention the following extension of the linear estimates from proposition \ref{estimlinloc} in Besov spaces constructed on $L^r$ spaces with $r\neq 2$. This high-frequency result can be obtained from  the Fourier expression of the solution of the linear system as done in \cite{CD} (see section $2.1$).

\begin{prop} \sl{Under the same assumptions as in proposition \ref{estimlinseuil} (adapted to $r\neq 2$), there exists a constant such that for all  $j>\overline{j}_0$,
\begin{multline}
\|v_j\|_{L_t^\infty L^r}+ \nu 2^{2j}\|v_j\|_{L_t^1 L^r} +(1+\nu2^j)\left(\|q_j\|_{L_t^\infty L^r}+ \frac{\nu}{\ee^2}\|q_j\|_{L_t^1 L^r}\right) \leq \\
C \max(1,M) \left((1+\nu 2^j)\|q_{0,j}\|_{L^r} +(1+\frac{1}{\sqrt{p}}) \|v_{0,j}\|_{L^r}\right),\\
\end{multline}
}
\end{prop}

Using the Lagrangian change of variable, remark \ref{Lr} and the methods from \cite{CD}, we obtain the analogous in the $L^r$-setting of theorem \ref{estimapriori} for high frequencies, which is similar to the key proposition 6 from \cite{CD} that leads to the main result of \cite{CD}. Let us recall that in this case, the regularity index of the Besov spaces is $\frac{d}{r}-1$ which is negative when $d<r$, that allows initial data with large modulus provided that they have fast enough oscillations.

\section{Appendix}

The first part is devoted to a quick presentation of the Littlewood-Paley theory and specific properties for hybrid Besov norms used in this paper. The second section to general considerations on flows.

\subsection{Besov spaces}

\subsubsection{Littlewood-Paley theory}

As usual, the Fourier transform of $u$ with respect to the space variable will be denoted by $\mathcal{F}(u)$ or $\hat{u}$. 
In this section we will briefly state (as in \cite{CD}) classical definitions and properties concerning the homogeneous dyadic decomposition with respect to the Fourier variable. We will recall some classical results and we refer to \cite{Dbook} (Chapter 2) for proofs (and more general properties).

To build the Littlewood-Paley decomposition, we need to fix a smooth radial function $\chi$ supported in (for example) the ball $B(0,\frac{4}{3})$, equal to 1 in a neighborhood of $B(0,\frac{3}{4})$ and such that $r\mapsto \chi(r.e_r)$ is nonincreasing over $\R_+$. So that if we define $\varphi(\xi)=\chi(\xi/2)-\chi(\xi)$, then $\varphi$ is compactly supported in the annulus $\{\xi\in \R^d, c_0=\frac{3}{4}\leq |\xi|\leq C_0=\frac{8}{3}\}$ and we have that,
\begin{equation}
 \forall \xi\in \R^d\setminus\{0\}, \quad \sum_{l\in\Z} \varphi(2^{-l}\xi)=1.
\label{LPxi}
\end{equation}
Then we can define the \textit{dyadic blocks} $(\ddl)_{l\in \Z}$ by $\ddl:= \varphi(2^{-l}D)$ (that is $\hat{\ddl u}=\varphi(2^{-l}\xi)\hat{u}(\xi)$) so that, formally, we have
\begin{equation}
u=\Sum_l \ddl u
\label{LPsomme} 
\end{equation}
As (\ref{LPxi}) is satisfied for $\xi\neq 0$, the previous formal equality holds true for tempered distributions \textit{modulo polynomials}. A way to avoid working modulo polynomials is to consider the set $\cS_h'$ of tempered distributions $u$ such that
$$
\lim_{l\rightarrow -\infty} \|\dot{S}_l u\|_{L^\infty}=0,
$$
where $\dot{S}_l$ stands for the low frequency cut-off defined by $\dot{S}_l:= \chi(2^{-l}D)$. If $u\in \cS_h'$, (\ref{LPsomme}) is true and we can write that $\dot{S}_l u=\Sum_{k\leq l-1} \ddq u$. We can now define the homogeneous Besov spaces used in this article:
\begin{defi}
\label{LPbesov}
 For $s\in\R$ and  
$1\leq p,r\leq\infty,$ we set
$$
\|u\|_{\dot B^s_{p,r}}:=\bigg(\sum_{l} 2^{rls}
\|\Delta_l  u\|^r_{L^p}\bigg)^{\frac{1}{r}}\ \text{ if }\ r<\infty
\quad\text{and}\quad
\|u\|_{\dot B^s_{p,\infty}}:=\sup_{l} 2^{ls}
\|\Delta_l  u\|_{L^p}.
$$
We then define the space $\dot B^s_{p,r}$ as the subset of  distributions $u\in {\cS}'_h$ such that $\|u\|_{\dot B^s_{p,r}}$ is finite.
\end{defi}
Once more, we refer to \cite{Dbook} (chapter $2$) for properties of the inhomogeneous and homogeneous Besov spaces. Among these properties, let us mention:
\begin{itemize}
\item for any $p\in[1,\infty]$ we have the following chain of continuous embeddings:
$$
\dot B^0_{p,1}\hookrightarrow L^p\hookrightarrow \dot B^0_{p,\infty};
$$
\item if $p<\infty$ then 
  $\dot{B}^{\frac dp}_{p,1}$ is an algebra continuously embedded in the set of continuous 
  functions decaying to $0$ at infinity, in particular we make in this paper an extensive use of the injection $\dot{B}_{2,1}^\fd \hookrightarrow L^{\infty}$;
    \item for any  smooth homogeneous  of degree $m$ function $F$ on $\R^d\setminus\{0\}$
the operator $F(D)$ maps  $\dot B^s_{p,r}$ in $\dot B^{s-m}_{p,r}.$ This implies that the gradient operator maps $\dot B^s_{p,r}$ in $\dot B^{s-1}_{p,r}.$  
  \end{itemize}
We refer to \cite{Dbook} (lemma 2.1) for the following result describing how derivatives act on spectrally localized functions:
\begin{lem}(Bernstein lemma)
\label{lpfond}
{\sl
Let  $0<r<R.$   A
constant~$C$ exists so that, for any nonnegative integer~$k$, any couple~$(p,q)$ 
in~$[1,\infty]^2$ with  $q\geq p\geq 1$ 
and any function $u$ of~$L^p$,  we  have for all $\lambda>0,$
$$
\displaylines{
{\rm Supp}\, \widehat u \subset   B(0,\lambda R)
\Longrightarrow
\|D^k u\|_{L^q} \leq
 C^{k+1}\lambda^{k+N(\frac{1}{p}-\frac{1}{q})}\|u\|_{L^p};\cr
{\rm Supp}\, \widehat u \subset \{\xi\in\R^N\,/\, r\lambda\leq|\xi|\leq R\lambda\}
\Longrightarrow C^{-k-1}\lambda^k\|u\|_{L^p}
\leq
\|D^k u\|_{L^p}
\leq
C^{k+1}  \lambda^k\|u\|_{L^p}.
}$$
}
\end{lem}
This implies the following embedding result:
\begin{prop}\label{LP:embed}
\sl{For all $s\in\R,$ $1\leq p_1\leq p_2\leq\infty$ and $1\leq r_1\leq r_2\leq\infty,$
  the space $\dot B^{s}_{p_1,r_1}$ is continuously embedded in 
  the space $\dot B^{s-d(\frac1{p_1}-\frac1{p_2})}_{p_2,r_2}.$}
\end{prop}

In this paper, we mainly work with functions or distributions depending on both the time variable $t$ and the space variable $x.$ We denote by $\cC(I;X)$ the set of continuous functions on $I$ with values in $X.$ For $p\in[1,\infty]$, the notation $L^p(I;X)$ stands for the set of measurable functions on  $I$ with values in $X$ such that $t\mapsto \|f(t)\|_X$ belongs to $L^p(I)$.

In the case where $I=[0,T],$  the space $L^p([0,T];X)$ (resp. $\cC([0,T];X)$) will also be denoted by $L_T^p X$ (resp. $\cC_T X$). Finally, if $I=\R^+$ we alternately use the notation $L^p X.$

The Littlewood-Paley decomposition enables us to work with spectrally localized (hence smooth) functions rather than with rough objects. We naturally obtain bounds for each dyadic block in spaces of type $L^\rho_T L^p.$  Going from those type of bounds to estimates in  $L^\rho_T \dot B^s_{p,r}$ requires to perform a summation in $\ell^r(\Z).$ When doing so however, we \emph{do not} bound the $L^\rho_T \dot B^s_{p,r}$ norm for the time integration has been performed \emph{before} the $\ell^r$ summation.
This leads to the following notation (after J.-Y. Chemin and N. Lerner in \cite{CL}):

\begin{defi}\label{d:espacestilde}
For $T>0,$ $s\in\R$ and  $1\leq r,\rho\leq\infty,$
 we set
$$
\|u\|_{\tilde L_T^\rho \dot B^s_{p,r}}:=
\bigl\Vert2^{js}\|\ddq u\|_{L_T^\rho L^p}\bigr\Vert_{\ell^r(\Z)}.
$$
\end{defi}
One can then define the space $\tilde L^\rho_T \dot B^s_{p,r}$ as the set of  tempered distributions $u$ over $(0,T)\times \R^d$ such that $\lim_{q\rightarrow-\infty}\dot S_q u=0$ in $L^\rho([0,T];L^\infty(\R^d))$ and $\|u\|_{\tilde L_T^\rho \dot B^s_{p,r}}<\infty.$ The letter $T$ is omitted for functions defined over $\R^+.$ 
The spaces $\tilde L^\rho_T \dot B^s_{p,r}$ may be compared with the spaces  $L_T^\rho \dot B^s_{p,r}$ through the Minkowski inequality: we have
$$
\|u\|_{\tilde L_T^\rho \dot B^s_{p,r}}
\leq\|u\|_{L_T^\rho \dot B^s_{p,r}}\ \text{ if }\ r\geq\rho\quad\hbox{and}\quad
\|u\|_{\tilde L_T^\rho \dot B^s_{p,r}}\geq
\|u\|_{L_T^\rho \dot B^s_{p,r}}\ \text{ if }\ r\leq\rho.
$$
All the properties of continuity for the product and composition which are true in Besov spaces remain true in the above  spaces. The time exponent just behaves according to H\"older's inequality. 
\medbreak
Let us now recall a few nonlinear estimates in Besov spaces. Formally, any product of two distributions $u$ and $v$ may be decomposed into 
\begin{equation}\label{eq:bony}
uv=T_uv+T_vu+R(u,v), \mbox{ where}
\end{equation}
$$
T_uv:=\sum_l\dot S_{l-1}u\ddl v,\quad
T_vu:=\sum_l \dot S_{l-1}v\ddl u\ \hbox{ and }\ 
R(u,v):=\sum_l\sum_{|l'-l|\leq1}\ddl u\,\dot\Delta_{l'}v.
$$
The above operator $T$ is called ``paraproduct'' whereas $R$ is called ``remainder''. The decomposition \eqref{eq:bony} has been introduced by J.-M. Bony in \cite{BJM}.

In this article we will frequently use the following estimates (we refer to \cite{Dbook} section 2.6, \cite{Dinv}, \cite{Has1} for general statements, more properties of continuity for the paraproduct and remainder operators, sometimes adapted to $\tilde L_T^\rho \dot B^s_{p,r}$ spaces): under the same assumptions there exists a constant $C>0$ such that:
\begin{equation}
\|\dot{T}_u v\|_{\dot{B}_{2,1}^s}\leq C \|u\|_{L^\infty} \|v\|_{\dot{B}_{2,1}^s}\leq C \|u\|_{\dot{B}_{2,1}^\fd} \|v\|_{\dot{B}_{2,1}^s},
\label{estimbesov}
\end{equation}
$$
\|\dot{T}_u v\|_{\dot{B}_{2,1}^{s+t}}\leq C\|u\|_{\dot{B}_{\infty,\infty}^t} \|v\|_{\dot{B}_{2,1}^s} \leq C\|u\|_{\dot{B}_{2,1}^{t+\fd}} \|v\|_{\dot{B}_{2,1}^s} \quad (t<0),
$$
$$\|\dot{R}(u,v)\|_{\dot{B}_{2,1}^{s_1+s_2}} \leq C\|u\|_{\dot{B}_{\infty,\infty}^{s_1}} \|v\|_{\dot{B}_{2,1}^{s_2}} \leq C\|u\|_{\dot{B}_{2,1}^{s_1+\fd}} \|v\|_{\dot{B}_{2,1}^{s_2}} \quad (s_1+s_2>0),
$$
$$
 \|\dot{R}(u,v)\|_{\dot{B}_{2,1}^{s_1+s_2-\fd}} \leq C\|\dot{R}(u,v)\|_{\dot{B}_{1,1}^{s_1+s_2}} \leq C\|u\|_{\dot{B}_{2,1}^{s_1}} \|v\|_{\dot{B}_{2,1}^{s_2}} \quad (s_1+s_2>0).
$$

\subsubsection{Complements for hybrid Besov spaces}

As explained, in the compressible Navier-Stokes system, the density fluctuation has two distinct behaviours in some low and high frequencies, separated by a frequency threshold. This leads to the definition of the hybrid Besov spaces. Let us begin with the spaces that are introduced by R. Danchin in \cite{Dinv} or \cite{Dbook} (we will use these spaces only in the appendix to prove estimates with the Hybrid norms introduced in (\ref{normhybride1})):

\begin{defi}
 \sl{
For $\alpha>0$, $r\in [0, \infty]$ and $s\in \R$ we denote
$$
\|u\|_{\tilde{B}_\alpha^{s,r}} \overset{def}{=} \Sum_{l\in \Z} 2^{ls} \max(\alpha, 2^{-l})^{1-\frac{2}{r}}\|\ddl u\|_{L^2}
$$}
\end{defi}
For example with $r\in \{1,\infty\}$:
$$
\|u\|_{\tilde{B}_\alpha^{s,\infty}}= \Sum_{l\leq \log_2(\frac{1}{\alpha})} 2^{l(s-1)} \|\ddl u\|_{L^2}+ \Sum_{l> \log_2(\frac{1}{\alpha})} \alpha 2^{ls} \|\ddl u\|_{L^2}, \mbox{ and}
$$
$$
\|u\|_{\tilde{B}_\alpha^{s,1}}= \Sum_{l\leq \log_2(\frac{1}{\alpha})} 2^{l(s+1)} \|\ddl u\|_{L^2}+ \Sum_{l> \log_2(\frac{1}{\alpha})} \frac{1}{\alpha} 2^{ls} \|\ddl u\|_{L^2},
$$
\begin{rem}
\sl{As stated in \cite{Dbook} we have the equivalence
$$
\frac{1}{2} \left( \|u\|_{\dot{B}_{2,1}^{s-1}}+ \alpha \|u\|_{\dot{B}_{2,1}^s}\right)\leq \|u\|_{\tilde{B}_\alpha^{s,\infty}} \leq \|u\|_{\dot{B}_{2,1}^{s-1}}+ \alpha \|u\|_{\dot{B}_{2,1}^s}.
$$
} 
\end{rem}
We refer to \eqref{normhybride1} and Proposition \ref{normhybride2} for the precise expression of the hybrid norm used in the present article. Let us just mention that this particuliar hybrid norm is accurate for our problem, but it is also related to the hybrid norms introduced by R. Danchin:
$$
\|.\|_{\dot{B}_\ee^{s+2,s}}=\|.\|_{\tilde{B}_\ee^{s,\frac{2}{3}}}
$$
Let us now state the following result (Proposition $5$ from \cite{CH})
\begin{prop}
 \sl{
Let $s\in \R$, $\alpha>0$. For all $q\in \Tilde{B}_{\alpha}^{s,\infty}\cap \Tilde{B}_{\alpha}^{s,1}$, we have
$$
\|q\|_{\dot{B}_{2,1}^s}^2 \leq \|q\|_{\Tilde{B}_{\alpha}^{s,\infty}} \|q\|_{\Tilde{B}_{\alpha}^{s,1}}
$$
}
\end{prop}
\begin{rem}
\sl{
For all $q\in \dot{B}_{2,1}^{s-1}\cap \dot{B}_{2,1}^{s}= \Tilde{B}_1^{s,\infty}$ we have
$$
\|q\|_{\Tilde{B}_1^{s,\infty}}\leq \|q\|_{\dot{B}_{2,1}^{s-1}} +\|q\|_{\dot{B}_{2,1}^{s}}
$$
and when $\ee>0$ is small enough, for all $q\in \dot{B}_{\ee}^{s+1,s}$, we have
$$
\|q\|_{\Tilde{B}_1^{s,1}}\leq \|q\|_{\tilde{B}_\ee^{s,1}} \leq \|q\|_{\dot{B}_{\ee}^{s+1,s}} \leq \|q\|_{\dot{B}_{\ee}^{s+1,s-1}} +\|q\|_{\dot{B}_{\ee}^{s+2,s}},
$$
}
\end{rem}
so we can use the hybrid norms introduced in (\ref{normhybride1}) and we will in fact use the following results:
\begin{prop}
 \sl{Let $s\in \R$. There exists a constant $C>0$ such that for all $0<\ee<1$, and all $q\in \dot{B}_{2,1}^{s-1}\cap \dot{B}_{2,1}^{s}\cap \dot{B}_{\ee}^{s+1,s-1} \cap \dot{B}_{\ee}^{s+2,s} $, we have
$$
\|q\|_{\dot{B}_{2,1}^s}^2 \leq C (\|q\|_{\dot{B}_{2,1}^{s-1}} +\|q\|_{\dot{B}_{2,1}^{s}}) (\|q\|_{\dot{B}_{\ee}^{s+1,s-1}} +\|q\|_{\dot{B}_{\ee}^{s+2,s}})
$$}
\label{estimhyb1}
\end{prop}

\subsection{Estimates for the flow of a smooth vector-field}
In this section, we recall classical estimates for the flow 
of a smooth vector-field with bounded spatial derivatives. We refer to \cite{Dlagrangien} or \cite{CD} for more details. We also refer to \cite{TH1} for the incompressible Navier-Stokes case.
\begin{prop}
\label{p:flow}
Let $v$ be a smooth globally Lipschitz time dependent  vector-field. Let $W(t) :=\int_0^t\|\nabla v(t')\|_{L^\infty}\,dt'.$  Let $\psi_t$ 
satisfy
$$
\psi_t(x)=x+\int_0^t v(t',\psi_{t'}(x))\,dt'.
$$
Then for all $t\in\R,$ the flow $\psi_t$ is a smooth diffeomorphism over $\R^d$ and one has  if $t\geq0,$
$$
\|D\psi_t^{\pm1}\|_{L^\infty}\leq e^{W(t)},
$$
$$
\|D\psi_t^{\pm1}-I_d\|_{L^\infty}\leq  e^{W(t)}-1,
$$
$$
\|D^2\psi_t^{\pm1}\|_{L^\infty}\leq e^{W(t)} \int_0^t\|D^2v(t')\|_{L^\infty} e^{W(t')} dt',
$$
$$
\|D^3\psi_t^{\pm1}\|_{L^\infty}\leq e^{W(t)} \int_0^t\|D^3v(t')\|_{L^\infty} e^{2W(t')} dt' +3\biggl(e^{V(t)}\int_0^t \|D^2v(t')\|_{L^\infty}e^{W(t')} dt'\biggr)^{2}.
$$
\end{prop}
As in \cite{CD} we also introduce the two-parameter flow $(t,t',x)\mapsto X(t,t',x)$ which is (uniquely) defined by
\begin{equation}
\label{flot2param}
X(t,t',x)=x+\int_{t'}^tv\bigl(t'',X(t'',t',x)\bigr)\,dt''.
\end{equation}
Uniqueness for Ordinary Differential Equations entails that
$$
X(t,t'',X(t'',t',x))=X(t,t',x).
$$ 
Hence 
$\psi_t=X(t,0,\cdot)$ and $\psi_t^{-1}=X(0,t,\cdot).$ 
\begin{prop}
 \sl{Under the previous notations, the jacobian determinant of $X$ satisfies:
\begin{equation}
 det(DX(t,t',x))=e^{\int_{t'}^t (\div v)(\tau, X(\tau,t',x))d\tau},
\end{equation}
and
$$
\begin{cases}
det(D\psi_t(x))=e^{\int_{0}^t (\div v)(\tau,\psi_\tau(x))d\tau},\\
det(D\psi_t^{-1}(x))=e^{-\int_{0}^t (\div v)(\tau, X(\tau,t,x))d\tau} =e^{-\int_{0}^t (\div v)(\tau,\psi_\tau\circ \psi_t^{-1}(x))d\tau}.
\end{cases}
$$
}
\label{detjacobien}
\end{prop}
\textbf{Proof:} differentiating \eqref{flot2param} with respect to $x,$ one gets by virtue of the chain rule,
\begin{equation}
\label{flow4} 
DX(t,t',x)=I_d+\int_{t'}^t  Dv(\tau,X(\tau,t',x))\cdot DX(\tau,t',x)\,d\tau.
\end{equation}
This immediately implies that:
$$
\partial_t (DX)(t,t',x)=Dv(t,X(t,t',x))\cdot DX(t,t',x),
$$
and
$$
\partial_t det(DX(t,t',x))=tr\left(Dv(t,X(t,t',x))\right)\cdot det(DX(t,t',x)),
$$
so that we obtain the result. $\blacksquare$.
\\

The authors wish to thank Rapha\"el Danchin, Taoufik Hmidi, Miguel Rodrigues and the anonymous referees for useful remarks and suggestions.

\end{document}